\documentclass[10pt,a4paper]{article}
 \usepackage{graphics}
 \usepackage{graphicx}
 \usepackage{epsfig}

\usepackage{amssymb}
\usepackage{amsthm}
\usepackage{amsmath}
\usepackage[left=2.5cm,right=2.5cm,bottom=3.5cm,top=3.5cm]{geometry}
\usepackage{empheq}
\usepackage{epstopdf}
\usepackage{changepage}
\usepackage{rotating}
\usepackage{adjustbox}
\usepackage{graphbox}
\usepackage{color}
\usepackage{xcolor}
\usepackage{dsfont}
\usepackage{float}
\usepackage{resizegather}
\usepackage{multicol}
\usepackage{booktabs} 
\usepackage{colortbl} 
\usepackage{multirow} 
\usepackage{longtable,tabu} 
\usepackage{hhline}   
\usepackage{anyfontsize} 
\usepackage{everysel}
\usepackage{psfrag}
\usepackage{titlesec}
\usepackage{bbm}
\usepackage{bm}
\usepackage{pgf,tikz}
\usetikzlibrary{arrows}

\titleformat*{\section}{\large\bfseries}
\titleformat*{\subsection}{\bfseries}
\titleformat*{\subsubsection}{\bfseries}

\allowdisplaybreaks 

\definecolor{cr1}{RGB}{200,0,0}
\definecolor{cr2}{RGB}{0,0,200}
\definecolor{cr12}{RGB}{100,0,100}

\newcommand{\diff}[2]{\frac{\partial #1}{\partial #2}}

\newcommand{\q}{\mathbf{q}}
\newcommand{\uu}{\mathbf{u}}

\begin{document}

\begin{center}
	\textbf{ \Large{High order ADER-DG schemes for the simulation of linear seismic waves induced by nonlinear dispersive free-surface water waves} }
	
	\vspace{0.5cm}
	{C. Bassi$^{(b)}$, S. Busto$^{(a)}$, M. Dumbser$^{(a)}$\footnote{Corresponding author.\\ \hspace*{0.25cm} Email addresses: caterina.bassi@polimi.it (C. BAssi), saray.busto@unitn.it (S. Busto), michael.dumbser@unitn.it (M. Dumbser)}}
	
	\vspace{0.2cm}
	{\small
		\textit{$^{(a)}$ Laboratory of Applied Mathematics, DICAM, University of Trento, via Mesiano 77, 38123 Trento, Italy}
		
		\textit{$^{(b)}$ MOX--Modelling and Scientific Computing, Dipartimento di Matematica, Politecnico di Milano, Piazza Leonardo da Vinci 32, 20133 Milano, Italy}
	}
\end{center}

\vspace{0.4cm}

\hrule
\vspace{0.4cm}

\noindent \textbf{Abstract}

\vspace{0.1cm}
In this paper, we propose a unified and high order accurate fully-discrete one-step ADER Discontinuous Galerkin method for the simulation of linear seismic waves in the sea bottom that are generated by the propagation of free surface water waves. 
In particular, a hyperbolic reformulation of the Serre-Green-Naghdi model for nonlinear dispersive free surface flows is coupled with a first order velocity-stress formulation for linear elastic wave propagation in the sea bottom. To this end, Cartesian non-conforming meshes are defined in the solid and fluid domains and the coupling is achieved by an appropriate time-dependent pressure boundary condition in the three-dimensional domain for the elastic wave propagation, where the pressure is a combination of hydrostatic and non-hydrostatic pressure in the water column above the sea bottom. The use of a first order hyperbolic reformulation of the nonlinear dispersive free surface flow model leads to a straightforward coupling with the linear seismic wave equations, which are also written in first order hyperbolic form.  
It furthermore allows the use of explicit time integrators with a rather generous CFL-type time step restriction associated with the dispersive water waves, compared to numerical schemes applied to classical dispersive models that contain higher order derivatives and typically require implicit solvers. Since the two systems that describe the seismic waves and the free surface water waves are written in the same form of a first order hyperbolic system they can also be efficiently solved in a unique numerical framework. In this paper we choose the family of arbitrary high order accurate discontinuous Galerkin finite element schemes, which have already shown to be suitable for the numerical simulation of wave propagation problems.  
The developed methodology is carefully assessed by first considering several benchmarks for each system separately, i.e. in the framework of linear elasticity and non-hydrostatic free surface flows, showing a good  agreement with exact and numerical reference solutions. Finally, also coupled test cases are addressed.  
Throughout this paper we assume the elastic deformations in the solid to be sufficiently small so that their influence on the free surface water waves can be neglected.

\vspace{0.2cm}

\vspace{0.4cm}

\hrule
\section{Introduction} \label{sec:introduction}
The physical phenomenon we are interested in is the generation of seismic waves in the sea bottom due to the propagation of free surface water waves on the sea surface in near coastal regions.
In view of the physical characteristics of the two materials involved, the wave speeds in the solid medium and the free surface wavespeed in the fluid can differ by up to two orders of magnitude. To simulate such a complex situation, we propose the use of a high order accurate fully-discrete one-step ADER-DG scheme on non conforming meshes, which solve a coupled set of two first order hyperbolic systems. The first model is a hyperbolic reformulation \cite{EDC19} of the Serre-Green-Naghdi model for the description of non-hydrostatic free-surface flows, while the second model consists in a classical first order hyperbolic model of linear elasticity in velocity-stress formulation for the description of linear seismic waves propagation \cite{Virieux1984,BedfordDrumheller,LeVeque:2002a}.   

The most straightforward way to model water wave propagation could be to consider a fully three-dimensional free surface flow model, like for example those developed in  a series of papers by Casulli et al. \cite{Casulli1999,CasulliZanolli2002,BrugnanoCasulli,Casulli2009,CasulliStelling2011,Casulli2014}, where an  efficient semi-implicit method for fully three-dimensional hydrostatic and non-hydrostatic free surface flows has been proposed. Since, however, we are interested in the propagation of free surface sea waves near the coast, where the typical horizontal length scales are far larger than the typical vertical length scales, also simplified shallow water-type equations (SWE) may be a suitable alternative, substantially reducing computational complexity and increasing the efficiency of the final methodology compared to fully three-dimensional models. 
A wide variety of phenomena associated to the propagation of water waves can be successfully described employing the classical shallow water equations \cite{saintvenant:1871}. Nevertheless, the SWE are unable to reproduce non-hydrostatic effects and propagation of solitary and dispersive waves. For this reason, we need to go beyond the classical shallow water equations, looking for more sophisticated dispersive systems. Starting from the pioneering work \cite{boussinesq:1872}, in which a first 1D Boussinesq-type model is derived under the assumptions of weak dispersion, weak non-linearity and flat bottom topography, different dispersive systems have been proposed in the literature. Among them we recall the Peregrine system, \cite{Per67}, where the model in \cite{boussinesq:1872} is extended to two dimensions and to non flat bottom topography (still maintaining the weak nonlinearity hypothesis), the Serre model, \cite{serre:1953}, where a 1D fully non-linear approach is presented for flat bottom topography (still keeping the weak dispersion assumption), the Serre-Green-Naghdi (SGN) model, \cite{SGN2,SRT87}, where an extension of \cite{serre:1953} to two dimensions for arbitrary bottom is provided and \cite{cienfuegos:2006}, where a derivation of the model in \cite{SGN2} is built, using asymptotic expansion and irrotationality. 
Notice that, in addition to the models just mentioned, more advanced models with improved dispersion characteristics (see \cite{MMS91,MS92,nwogu:1993,madsen:2003,madsen:2010})  and models which include additional physical phenomena with respect to the classical formulation for inviscid flows (see the review \cite{kirby:2016}) have been also proposed in the literature. 

An important distinguishing feature of dispersive models is that they contain higher order space and space-time derivatives. Consequently, their numerical discretization becomes particularly difficult and a severe time step restriction arises when explicit time integration schemes are employed. A possible solution to this major drawback is the introduction of augmented first order systems, such as in \cite{JSMarie,fernandeznieto:2018}. Also in this case, however, the hyperbolicity of the SW equations is lost, thus leading to the necessity of solving an elliptic equation for the dispersive part at each time-step. A completely different approach is instead adopted in \cite{EDC19}, where a hyperbolic approximation of the non-hydrostatic system in \cite{JSMarie} and of the SGN model with mild bottom approximation are proposed. The hyperbolic approximation proves to have a dispersion relation which is very similar to the dispersion relations associated to the original non-hyperbolic models and allows to obtain very accurate numerical results. At the same time hyperbolicity allows to easily implement the model in the context of high-order finite volume (FV) and discontinuous Galerkin (DG) schemes and to realize efficient numerical simulations also in multiple space dimensions. Due to these evident advantages, we have decided to employ this approach for the water waves part in our coupled simulations. Notice that the introduction of hyperbolic reformulations is not a novelty and comes from the pioneering work by Cattaneo \cite{cattaneo:1958}, where the second order terms in the heat equation are rewritten as relaxation terms, while a hyperbolic approximation has been proposed for the first time in the context of dispersive systems in \cite{HGN} (precisely for the Serre model \cite{serre:1953}). Further developments of \cite{EDC19} are presented in \cite{BBBD20}, where an hyperbolic reformulation of SGN model for arbitrary bathymetry in introduced. 
Moreover, in \cite{EL20}, an even more general formulation has been proposed, including further dispersive Boussinesq-type systems, as the models of Yamazaki et al., \cite{YKC09}, and Peregrine, \cite{Per67,MMS91}, which are built neglecting some convective terms in the equation for the averaged vertical velocity and in an auxiliary equation accounting for the spatial variation of the mean horizontal velocity. For further first order hyperbolic reformulations of dispersive and dissipative systems, see also \cite{Dhaouadi2018} and \cite{PeshRom2014,GPRmodel}.
Concerning the propagation of seismic waves, we put ourselves in the context of linear elasticity and we adopt a first order velocity-stress formulation, which has the advantage of being hyperbolic, see e.g. \cite{LeVeque:2002a,gij1}.

A key idea of the present paper is the mathematical description of both wave propagation problems, i.e. the linear seismic waves as well as the nonlinear dispersive water waves, at the aid of first order hyperbolic systems. For this reason, we can easily couple both models with each other and solve them numerically in the time domain by employing high order accurate discontinuous Galerkin (DG) finite element methods in a straightforward manner. We decide, in particular, to use an explicit high order ADER-DG scheme \cite{gij1}, for which the usual CFL condition holds, i.e. with $\Delta t$ proportional to $\Delta x$, thus avoiding higher powers of $\Delta x$ that would be typical for explicit time discretizations of Boussinesq-type equations with higher order spatial and temporal derivatives. 

The DG method has been introduced for the first time in \cite{reed}, for the solution of a neutron transport equation. It has been then extended to time-dependent, multidimensional and nonlinear hyperbolic problems in \cite{chavent:1989,cbs1,cbs2,cbs3,cbs0}. Besides, in \cite{BassiRebay2,CBS-convection-diffusion} the Local Discontinuous Galerkin (LDG method) has been proposed in order to solve convection-diffusion equations. This approach involves rewriting the second order equations as an augmented non-hyperbolic first order system and the subsequent discretization of this augmented system with the DG method.
The DG method has been applied for the first time to equations containing higher order derivatives in \cite{YanShu}, where the LDG method is used for the resolution of linear dispersive Korteveg-de-Vries (KdV) equations, containing up to third order spatial derivatives. Extensions to linear equations with derivatives up to fifth order and nonlinear dispersive equations are instead presented in \cite{YanShub,LevyShuYan}. Applications of the DG method to the solution of nonlinear Boussinesq-type dispersive equations have been introduced in \cite{eskilsson:2006a,eskilsson:2006b,engsig:2008}. Notice that the severe time step restrictions due to higher-order spatial derivatives are overcome in \cite{DumbserFacchini}, where a fully implicit  space-time DG method is applied to both linear third order KdV equations and nonlinear Boussinesq-type systems. 

High order of accuracy in space is straightforward to obtain in the DG framework, while attaining high order in time is still a very active field of research. A successful approach consists in using the already mentioned space-time DG methods, \cite{Sander2012,Rhebergen2013,spacetimedg1,spacetimedg2,3DSIINS,3DSICNS,KlaijVanDerVegt,TD18LE,BTBD20}. An alternative, that will be followed in this work, are the ADER-DG schemes, first put forward in \cite{DumbserEnauxToro} in the context of FV methods and generalized in \cite{Dumbser2008} to the unified $P_{N}P_{M}$ framework for arbitrary high order accurate one-step FV and DG schemes. 
Classical ADER methods (Arbitrary high order DErivative Riemann problem) have been proposed by Millington et al., \cite{mill}, and Toro et al., \cite{TMN01}, in the framework of finite volume methods. The methodology is based on the resolution, at the cell interfaces, of a generalized Riemann problem with piecewise polynomial initial conditions, built using a nonlinear reconstruction (e.g. ENO or WENO methods) that circumvents Godunov's theorem. 
Then, space-time integration on an appropriate control volume is performed using a Taylor series expansion in time, where time derivatives are replaced by spatial derivatives following the Cauchy-Kovalevskaya procedure. Further developments of classical ADER schemes, including their extension to the DG framework, can be found in \cite{toro4,titarevtoro,dumbser_jsc,gij5,gassner:2011a,BTVC16,BFTVC17,DMT2019,BBDFSVC20} and references therein. 
The main inconvenient of classical ADER methods is the need of the cumbersome Cauchy-Kovalevskaya procedure. The novel ADER-DG methodology presented in \cite{DumbserEnauxToro,Dumbser2008} avoids that step by employing a new element-local space-time DG predictor, which leads to more efficient algorithms. Since then, ADER-DG has been used to solve many different models, such as compressible flows \cite{ADERNSE}, hyperbolic systems in general relativity \cite{ADERCCZ4,ADERGRMHD,dumbser2020glm} or non conservative hyperbolic systems for geophysical flows \cite{ADERNC}. Moreover, the ADER-DG method has also proven to be well suited for the simulation of seismic waves problems, see \cite{gij1,gij2,gij3,gij4,gij5}. For other high order discontinuous Galerkin finite element
schemes applied to elastic wave equations, see e.g. \cite{GroteDG, Antonietti1, Antonietti2, Antonietti3, Antonietti4}. 
Finally, in \cite{EDC19} the ADER-DG method has also been successfully applied to the solution of hyperbolic reformulations of dispersive models. Due to these considerations, the ADER-DG method appears to be a suitable choice for the discretization of the coupled system proposed in the present work.

The paper is organized as follows. In Section \ref{sec:governing_equations}, the mathematical models employed for the description of both seismic and free surface water waves are recalled. Moreover, the boundary conditions to be set on the interface between the fluid and solid domains are defined. Section~\ref{sec:method} is devoted to the description of the high order one-step ADER-DG scheme on Cartesian grids. In Section \ref{sec:num_res}, the numerical method is validated. First, some classical benchmarks are solved independently for each of the two systems of equations. Then, two test cases of the coupled problem are presented. Finally, in Section \ref{sec:conclusions}, we draft some conclusions and perspectives.

\section{Governing equations}\label{sec:governing_equations}
As already mentioned within the introduction, we aim at simulating the effect of surface water waves on the generation and propagation of seismic waves on the sea bottom. To this end, we couple the two systems that are recalled in this section: a non-hydrostatic dispersive shallow water-type model for the propagation of the free surface water waves and a linear elasticity model for the seismic waves propagating in the solid domain below the sea floor.  

\subsection{Non-hydrostatic free surface flows}\label{sec:hsgn}
During the last decades many non-hydrostatic models aiming at characterising non-hydrostatic free surface flows have been successfully developed. 
Taking into account the features of the flow to be modelled, we will employ the  Serre-Green-Naghdi (SGN), \cite{SGN1,SGN2} model, and a dispersive model recently proposed by Sainte-Marie et al. (SM) in \cite{JSMarie}. Both models are used in combination with the so-called mild bottom approximation. It is important to remark that both systems can be rewritten in a one-parameter dependent unified formulation, \cite{EDC19}, as 
\begin{gather} 
\diff{ h}{ t}+  \nabla \cdot(h \uu) = 0, \label{eq:sgn_depth}\\
\diff{ h \uu}{ t} + \nabla \cdot \left( h \uu\otimes\uu\right)  + \nabla \cdot \left( \frac{1}{2}gh^2\mathbf{I}+h p \mathbf{I}\right)  +(gh + \gamma p)\nabla z_b=0, \label{eq:sgn_velocity}\\
\diff{h w}{t} + \nabla \cdot(h \uu  w) = \gamma  p, \label{eq:sgn_va_velocity}\\
w +\frac{1}{2} h \nabla \cdot \uu -\uu \cdot \nabla z_b = 0, \label{eq:sgn_w}
\end{gather}
where we have denoted $h$ the water depth, $\uu=(u_{x},u_{y})$ the horizontal velocity vector of the fluid, $w$ the auxiliary variable of the averaged vertical flow velocity, $p$ the depth-averaged non-hydrostatic correction for the pressure, $g=9.81$ the gravity acceleration and $z_b$ the vertical coordinate of the bottom bathymetry. Moreover, the gradient and divergence operators considered refer to the horizontal plane, $\nabla=\left(\partial_{x},\partial_{y}\right)^{T}$, neglecting the vertical variable and $\mathbf{I}$ is the identity matrix of dimension two. Substituting $\gamma=\frac{3}{2}$, the former system would provide the SGN model whereas setting $\gamma=2$ leads to the SM model.

\begin{figure}[h]
	\centering
	\includegraphics[width=0.7\linewidth]{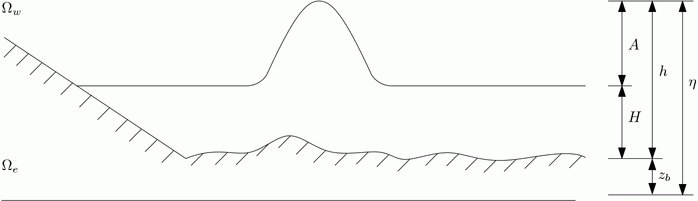}
	\begin{minipage}{0.85\linewidth}
		\caption{Two-dimensional computational domain for the simulation of the water wave propagation. $h$ denotes the water depth, $H$ is the still water depth, $A= h-H$ corresponds to the surface elevation with respect to the still water depth, $\eta$ is the total surface elevation, and $z_{b}$ gives the bottom bathymetry.}	\label{fig:plotswvariablebed}
	\end{minipage}
\end{figure}

Following \cite{EDC19}, the first order unified reformulation of the SGN and SM models reads 
\begin{gather} 
\diff{ h}{ t}+  \nabla \cdot(h \uu) = 0, \label{eq:hgn_depth}
\\
\diff{ h \uu}{ t} + \nabla \cdot \left( h \uu\otimes\uu\right)  + \nabla \cdot \left( \frac{1}{2}gh^2\mathbf{I}+h p \mathbf{I}\right)  +(gh + \gamma p)\nabla z_b=0, \label{eq:hgn_velocity}\\
\diff{h w}{t} + \nabla \cdot(h \uu  w) = \gamma  p, \label{eq:hgn_va_velocity}\\
\diff{h p}{t} + \nabla \cdot (h  \uu  p) + c^2 \left(2 w + h\nabla\cdot  \uu - 2 \uu\cdot\nabla z_b\right) = 0. \label{eq:hgn_pressure}
\end{gather}
with $c$ the artificial sound speed, $c=\alpha\sqrt{g H_{0}}$, $H_{0}$ the average still water depth and $\alpha>0$. 
The system \eqref{eq:hgn_depth}-\eqref{eq:hgn_pressure} represents the conservation of mass and momentum and is furthermore augmented by two PDEs for the auxiliary variables $p$ and $w$, with appropriate source terms, that allow the system 
to relax, for $c^{2} \rightarrow \infty$, towards the original SGN or SM models, depending on the value of the parameter $\gamma$. 
Moreover, assuming steady bathymetry, the system can be completed with the following equation: 
\begin{equation}
\partial_t z_{b} = 0. \label{eq:hgn_bottom}
\end{equation}
It is important to note that the above system is depth averaged and, thus, the vector of spatial coordinates $\mathbf{x} \in \Omega_w \subset \mathds{R}^2$, where $\Omega_w$ is the two-dimensional computational domain used for the simulation of the water wave propagation (see Figure \ref{fig:plotswvariablebed}). The fact that the equations are depth averaged, jointly with the use of a two-dimensional domain, makes this model particularly interesting concerning its computational efficiency, compared to a more complete three-dimensional non-hydrostatic formulation.

\subsection{Linear elastic wave propagation}\label{sec:linearelasticity}
Assuming the sea bottom to be a homogeneous isotropic elastic material under small deformations, it can be modelled using the linear elasticity equations \cite{Gur81,Ber05,BedfordDrumheller,LeVeque:2002a}. 
The classical formulation of linear elasticity is a second order vector wave equation, but it can be rewritten as a first order hyperbolic system in velocity-stress formulation, which leads to an easier and more direct coupling with the model governing the dispersive water waves shown in the previous section. Accordingly, the final hyperbolic system for the linear elastic wave propagation reads 
\begin{eqnarray}
\frac{\partial \bm{\sigma} }{\partial t}  - \lambda \left( \nabla \cdot \mathbf{v} \right) \mathbf{I} - \mu \left( \nabla \mathbf{v} + \nabla \mathbf{v}^T \right) = 0, \label{eqn:hooke} \\
\diff{\rho\mathbf{v}}{t}-\nabla \cdot \bm{\sigma} = 0,	\label{eqn:momentum} 
\end{eqnarray}  
where the first equation, \eqref{eqn:hooke}, is the Hooke law expressed in terms of the two Lam\'e constants $\lambda$ and $\mu$ and the second equation, 
\eqref{eqn:momentum}, represents the conservation of momentum. Here, the time is again denoted by $t$ and the spatial coordinate is  $\mathbf{x} \in \Omega_e \subset \mathds{R}^3$, 
with $\Omega_e$ a 3D computational domain used for the simulation of the seismic wave propagation.     
Furthermore, in the above system, $\bm{\sigma} = \bm{\sigma}^T$ is the symmetric stress tensor, $\tilde{\bm{\sigma}}=\frac{\partial \bm{\sigma}}{\partial t}$ corresponds to the linearized part of the first Piola-Kirchhoff stress tensor, $\mathbf{u}$ is the vector of the velocity field and $\rho$ is the density.

\subsection{Coupling of the models}
To simulate the effect of free surface waves on the generation and propagation of seismic waves on the sea bottom, the non-hydrostatic and the linear elastic wave models need to be coupled. Throughout this paper we assume that the coupling is only done in a one-way manner, i.e. the free surface waves lead to a time-dependent pressure on the sea floor, which generates low frequency and low amplitude linear elastic waves in the solid medium below the sea bottom. Instead, we assume that the elastic deformations of the solid are very small and therefore do \textit{not} couple back to the surface water waves via a time-dependent bottom geometry.  
The non-hydrostatic description of the free surface water waves propagation therefore provides appropriate boundary conditions for the stress tensor on the upper boundary of the domain $\Omega_e$, which will be denoted by $\partial \Omega_{e,u}$ in the following.  In particular, we assume zero shear stress and continuity of the normal stress on $\partial \Omega_{e,u}$, 
\begin{equation}
\sigma_{xz}=\sigma_{yz}=0, \qquad  \sigma_{zz} = \rho g h. \label{eq:boundaryconditions}
\end{equation}
As already stated before, due to their very low amplitude, the linear elastic waves are not coupled back to the non-hydrostatic shallow water solver, although one might imagine a fully coupled system by considering, in \eqref{eq:hgn_depth}-\eqref{eq:hgn_bottom}, $\partial_t z_{b} = \mathbf{v} \cdot \mathbf{n}_z$, where $\mathbf{n}_z$ is the unit normal vector in $z$ direction and $\mathbf{v}$ is the velocity of the solid.

\subsection{Unified writing of the models}
To provide a unified description of the numerical scheme for the two models considered, we rewrite the systems of equations to be solved in the general form
\begin{equation}
\partial_t \mathbf{Q}  + \nabla \cdot
\mathbf{F}(\mathbf{Q})+ \boldsymbol{\mathbf{B}}(\mathbf{Q}) \cdot\nabla
\mathbf{Q}=\mathbf{S}\left( \mathbf{Q} \right), \label{eq:general_pde}
\end{equation}  
where $\mathbf{Q}=\mathbf{Q}\left(\mathbf{x},t\right)$ is the vector of unknowns, $\mathbf{F}(\mathbf{Q})$ is the nonlinear flux tensor, $\boldsymbol{\mathbf{B}}(\mathbf{Q}) \cdot\nabla
\mathbf{Q}$ is a genuinely non-conservative term, and $\mathbf{S}\left( \mathbf{Q} \right)$ is an algebraic source term. Therefore, for the hyperbolic non-hydrostatic model \eqref{eq:hgn_depth}-\eqref{eq:hgn_bottom}, we have
\begin{eqnarray}
\mathbf{Q} = \begin{pmatrix}
h        \\
h u_{1}  \\
h u_{2}  \\
h w      \\
h p      \\
z_{b} 
\end{pmatrix}, \quad
\mathbf{F}(\mathbf{Q}) = 
\begin{pmatrix}
h u_{1} &  h u_{2}\\
h u_{1}^{2} +\frac{1}{2} g h^{2} + h p & h u_{1} u_{2}\\
h u_{1} u_{2} & h u_{2}^{2} +\frac{1}{2} g h^{2} + h p\\
h w u_{1} & h w u_{2}\\
h u_{1} \left( p + c^{2} \right) & h u_{2} \left( p + c^{2} \right)  \\
0  & 0
\end{pmatrix}, \notag\\
\boldsymbol{\mathbf{B}}(\mathbf{Q}) \cdot\nabla
\mathbf{Q} = \begin{pmatrix}
0  \\
\left( g h + \gamma p\right)\partial_{x} z_{b}  \\
\left( g h + \gamma p\right)\partial_{y} z_{b}  \\
0  \\
c^{2}\left[ - u_{1} \partial_{x}h  - u_{2} \partial_{y}h - 2 \left( u_{1} \partial_{x}z_{b} + u_{2} \partial_{y}z_{b}\right)  \right]   \\
0 
\end{pmatrix}, \quad
\mathbf{S}\left( \mathbf{Q} \right) = \begin{pmatrix}
0 \\
0 \\
0 \\
\gamma p\\
-2 c^{2} w\\
0 
\end{pmatrix}.
\end{eqnarray}
On the other hand, the definition of $\mathbf{Q} =\left(\sigma_{xx}, \sigma_{yy}, \sigma_{zz}, \sigma_{xy}, \sigma_{yz}, \sigma_{xz}, u_{1}, u_{2}, u_{3}\right)^{T} $ and
\begin{equation}
\boldsymbol{\mathbf{B}}(\mathbf{Q}) \cdot\nabla
\mathbf{Q} = \begin{pmatrix}
- (\lambda+2\mu)\partial_{x}u_{1} - \lambda        \partial_{y}u_{2} - \lambda        \partial_{z}u_{3} \\
- \lambda       \partial_{x}u_{1} - (\lambda+2\mu) \partial_{y}u_{2} - \lambda        \partial_{z}u_{3} \\ 
- \lambda       \partial_{x}u_{1} - \lambda        \partial_{y}u_{2} - (\lambda+2\mu) \partial_{z}u_{3} \\    
- \mu \partial_{x}u_{2}  - \mu \partial_{y}u_{1} \\    
- \mu \partial_{y}u_{3}  - \mu \partial_{z}u_{2} \\
- \mu \partial_{x}u_{3}  - \mu \partial_{z}u_{1} \\
- \frac{1}{\rho} \left( \partial_{x}\sigma_{xx} + \partial_{y}\sigma_{xy} + \partial_{z}\sigma_{xz} \right) \\
- \frac{1}{\rho} \left( \partial_{x}\sigma_{xy} + \partial_{y}\sigma_{yy} + \partial_{z}\sigma_{yz} \right) \\
- \frac{1}{\rho} \left( \partial_{x}\sigma_{xz} + \partial_{y}\sigma_{yz} + \partial_{z}\sigma_{zz} \right) 
\end{pmatrix}
\end{equation}
leads to the linear elasticity model \eqref{eqn:hooke}-\eqref{eqn:momentum}.

\section{Numerical discretization}\label{sec:method}
The high order accurate fully-discrete one-step ADER discontinuous Galerkin methodology (ADER-DG), \cite{Dumbser2008,Dumbser2014}, is used in order to discretize the two models considered, namely the hyperbolic reformulation of the SGN equations (HSGN), presented in Section \ref{sec:hsgn}, and the linear elasticity system recalled in Section \ref{sec:linearelasticity}. ADER-DG methods fall into the framework of the general  $P_{N}P_{M}$ schemes proposed in \cite{Dumbser2008}, that extend the local predictor ADER methodology presented in \cite{DumbserEnauxToro} for FV also to the DG framework. More precisely, we focus on the pure DG case, where $N=M$, which has shown to be appropriate to solve linear and non-linear hyperbolic conservation laws. In this section, we provide a brief summary of the method on Cartesian grids. Further developments of this methodology, including the use of unstructured mesh and Cartesian grids with adaptative mesh refinement (AMR) employed to solve a wide variety of hyperbolic PDEs, in both the Eulerian and the Lagrangian framework, can be found, for instance, in \cite{ADERNSE,GPRmodelMHD,Dumbser2014,Zanotti2015,DFTBW18,GBCKSD2019,ADERCCZ4,BCDGP20} and references therein.

Before describing the numerical method to be employed, we define a discretization of the computational domain $\Omega$ using a Cartesian grid made of elements of the form $\Omega_{i}=\left[x_{i}-\frac{1}{2}\Delta x,x_{i}+\frac{1}{2}\Delta x\right]\times \left[y_{i}-\frac{1}{2}\Delta y,y_{i}+\frac{1}{2}\Delta y\right]\times \left[z_{i}-\frac{1}{2}\Delta z,z_{i}+\frac{1}{2}\Delta z\right]$, 
with $\mathbf{x}_{i}=\left(x_{i},y_{i},z_{i}\right) $ the barycentre of cell $\Omega_{i}$ 
and $\Delta x$, $\Delta y$, $\Delta z$ the cell size on each spatial coordinate direction. 
Next, following the classical DG approach, we assume that the space of discrete solutions of \eqref{eq:general_pde} is generated by  spatial basis functions $\phi_k$ constructed as the tensor product of piecewise polynomials up to degree $N$. In particular, we consider the orthogonal Lagrange interpolation polynomials passing through the Gauss-Legendre quadrature points of a $N+1$ Gauss quadrature formula. Then, within each element $\Omega_i$ the discrete solution of the system can be written as
\begin{equation}
\mathbf{u}_h(\mathbf{x},t^n) = \phi_l (\mathbf{x})\;
\hat{\mathbf{u}}^n_l,   \quad \mathbf{x} \in
\Omega_i, \label{eq:disc_sol}
\end{equation}
where the classical Einstein summation convection is employed, $\hat{\mathbf{u}}^n_l$ denote the degrees of freedom of the solution, and $l$ is a multidimensional 
index referring to the one-dimensional basis functions, $\phi_{l_{m}}$, on a reference 																												
element $\Omega_{\mathrm{ref}}=\left[0,1\right]$, to be used in the tensor product. 
The reference coordinates 
$0 \leq \xi, \eta, \zeta \leq 1$ are obtained via the transformations 
$x = x_{i}-\frac{1}{2}\Delta x + \xi \Delta x$, 
$y = y_{i}-\frac{1}{2}\Delta y + \eta \Delta y$, and $z = z_{i}-\frac{1}{2}\Delta z + \zeta \Delta z$, respectively.
Since the chosen basis functions are not time dependent, the direct use of a classical DG approach would result in a low order scheme in time. To attain high order of accuracy in time, we use the ADER-DG methodology which can be divided into two main steps:
\begin{itemize}	
	\item Local space-time predictor computation. System \eqref{eq:general_pde} is solved ``in the small'' using a locally implicit space-time discontinuous Galerkin scheme on each element, which neglects iterations between neighbour cells.
	
	\item Explicit update of the solution using a one-step corrector. The space-time predictor is used into the weak formulation of \eqref{eq:general_pde} which takes into account the fluxes between cells and provides the solution of the system at the new time instant. 
\end{itemize}

We come now to further detail each step.

\subsection{Local space-time predictor}\label{sec:predictor}
To determine the local space-time predictor solution, $\q_h(\mathbf{x},t)$, which will 
lead to a high order scheme in space and time avoiding the cumbersome Cauchy-Kovalewskaya 
procedure used in the original ADER methods \cite{Toro2001,Toro2002,toro3,TT04}, we employ the weak formulation 
in space-time proposed in \cite{DumbserEnauxToro,Dumbser2008}. 
Let us consider the space-time test functions, $\theta_{k}=\theta_{k}(\mathbf{x},t) $, 
built as the product of the already introduced nodal spatial basis functions 
and an additional one dimensional basis function for the time dependency. with 
the additional transformation for the reference time $\tau$ given by $t = t^n + \tau \Delta t$.
Then,  multiplying \eqref{eq:general_pde} by $\theta_{k}$ and integrating over the space-time 
control volume, $\Omega_{i}\times\left[t^{n},t^{n+1}\right]$, yields
\begin{equation}
\int \limits_{t^n}^{t^{n+1}} \!\! \int\limits_{\Omega_i  }
\theta_k \, \partial_t \q_h \,d\mathbf{x} \, dt + \int
\limits_{t^n}^{t^{n+1}} \!\! \int\limits_{\Omega_i  } \theta_k \,
\nabla \cdot \mathbf{F}(\q_h) \,d\mathbf{x}\,dt  
+ \int
\limits_{t^n}^{t^{n+1}} \!\! \int\limits_{\Omega_i^{\circ} } \theta_k
\mathbf{\mathbf{B}}(\q_h ) \cdot \nabla \q_h \,d\mathbf{x}\,dt =
\int\limits_{\Omega_i } \theta_k
\mathbf{S}(\q_h )   \,d\mathbf{x}\,dt. 
\label{eq:predictor1}
\end{equation}
Within an implicit space-time DG method, \cite{spacetimedg1,spacetimedg2}, the weak formulation \eqref{eq:predictor1} would be now integrated by parts in space and time to provide the solution at the new step. However, we are just interested in obtaining a local approximation of the predictor so we only integrate by parts in time, neglecting the interaction between neighbours, 
\begin{eqnarray}
\int\limits_{\Omega_i }
\theta_k(\mathbf{x},t^{n+1}) \q_h(\mathbf{x},t^{n+1}) \,d\mathbf{x}  -  
\int\limits_{\Omega_i }
\theta_k(\mathbf{x},t^{n}) \mathbf{u}_h(\mathbf{x},t^{n}) \,d\mathbf{x}     
-\int \limits_{t^n}^{t^{n+1}} \!\! \int\limits_{\Omega_i  }
\partial_t \theta_k \,  \q_h \,d\mathbf{x} \, dt + 
\nonumber \\ 
\int
\limits_{t^n}^{t^{n+1}} \!\! \int\limits_{\Omega_i } \theta_k \,
\nabla \cdot \mathbf{F}(\q_h) \,d\mathbf{x}\,dt  
+ \int
\limits_{t^n}^{t^{n+1}} \!\! \int\limits_{\Omega_i^{\circ} } \theta_k
\mathbf{\mathbf{B}}(\q_h ) \cdot \nabla \q_h \,d\mathbf{x}\,dt =
\int\limits_{\Omega_i } \theta_k
\mathbf{S}(\q_h )   \,d\mathbf{x}\,dt.
\label{eq:predictor}
\end{eqnarray}
Moreover, in \eqref{eq:predictor}, we have taken into account that the predictor at time $t^{n}$ is given by the degrees of freedom of the solution at the previous time step, $\mathbf{u}_h(\mathbf{x},t^{n})$, thus respecting the causality principle. This fact can be seen as the realization of an upwinding approximation in time. Therefore, the above nonlinear system has as only unknown the degrees of freedom of the space-time expansion, 
\begin{equation}
\q_h(\boldsymbol{x},t) = \theta_k(\boldsymbol{x},t) \, \hat{\q}_k
\label{eq.stdof}
\end{equation}
and can be solved locally at each cell using a discrete Picard iteration procedure. 
Since the Picard iteration matrix is nilpotent, it will converge in at most $N+1$ iterations, as it has been proven in \cite{Jackson} for homogeneous linear conservation laws. The convergence of the Picard iteration for nonlinear systems of conservation laws was proven in \cite{BCDGP20}. The solution of \eqref{eq:predictor} constitutes the only (element-local) implicit step on the whole ADER-DG algorithm.

\subsection{Fully discrete one-step ADER-DG schemes}
The solution $\mathbf{q}_{h}$, obtained at the predictor step, does not account for the neighbouring flux contributions, so it can not be used as the solution of the PDE system at time $t^{n+1}$. To correct this issue we employ an explicit one-step DG approach.
We first multiply the governing PDE system \eqref{eq:general_pde}, by the test functions $\phi_k$ and we then integrate over a space-time control volume $\Omega_i \times [t^{n},t^{n+1}]$, obtaining the following weak problem 
\begin{equation}
\int \limits_{t^n}^{t^{n+1}} \int\limits_{\Omega_i } \phi_k \left(
\partial_t \mathbf{Q}  + \nabla \cdot
\mathbf{F}(\mathbf{Q})+ \boldsymbol{\mathbf{B}}(\mathbf{Q}) \cdot\nabla
\mathbf{Q}\right) \,d\mathbf{x}\,dt = 
\int \limits_{t^n}^{t^{n+1}} \int\limits_{\Omega_i } \phi_k \, \mathbf{S}\left( \mathbf{Q} \right).\label{eq:weakPDE}
\end{equation}
Taking into account \eqref{eq:disc_sol}, integrating the flux divergence term by 
parts in space and the time derivative by parts in time yields
\begin{gather}
\left( \, \int\limits_{\Omega_i} \phi_k \phi_l \, d\mathbf{x}\right)
\left( \hat{\mathbf{u}}_l^{n+1} - \hat{\mathbf{u}}_l^{n} \, \right) +
\int \limits_{t^n}^{t^{n+1}} \int\limits_{\partial
	\Omega_i } 
\phi_k \mathcal{G}\left(\q_h^-, \q_h^+ \right) \cdot
\mathbf{n} \, dS \, dt 
+ \int \limits_{t^n}^{t^{n+1}} \!\!
\int\limits_{\partial \Omega_i } \!\! \phi_k \mathcal{D}\left(
\q_h^-,\q_h^+ \right) \cdot \mathbf{n} \, dS \, dt \nonumber \\  -
\int \limits_{t^n}^{t^{n+1}} \!\! \int\limits_{\Omega_i } \!\!\! \nabla
\phi_k \cdot \mathbf{F}(\q_h) \,d\mathbf{x}\,dt + \int
\limits_{t^n}^{t^{n+1}} \!\! \int\limits_{\Omega_i^{\circ} } \phi_k
\mathbf{\mathbf{B}}(\mathbf{q}_h) \cdot \nabla \mathbf{q}_h
\,d\mathbf{x}\,dt = 
\int\limits_{\Omega_i  } \phi_k
\mathbf{\mathbf{S}}(\mathbf{q}_h)  
\,d\mathbf{x}\,dt,  
\label{eq:ADER-DG}
\end{gather}
where $\mathbf{n}$ denotes the outward unit normal at the cell boundary, 
$\partial \Omega_{i}$, and $\q_h=\q_h(\mathbf{x},t)$ is the local space-time predictor already introduced in Section \ref{sec:predictor}. 
Since we are using a discontinuous Galerkin scheme, the basis functions are allowed to jump across cell interfaces. To account for the discontinuities arising in the second term of \eqref{eq:ADER-DG}, we make use of a Riemann solver at the element interfaces, see e.g. \cite{toro-book}. 
The initial condition for the numerical flux function is then given by the right and left states $\q_h^-,\, \q_h^+$ computed at the predictor step which yields the order of accuracy sought. In particular, within this work, we consider the Rusanov numerical flux function, hence
\begin{equation}
\mathcal{G}\left(\q_h^-, \q_h^+ \right) \cdot \mathbf{n} = \frac{1}{2}
\left( \mathbf{F}(\q_h^+) + \mathbf{F}(\q_h^-) \right) \cdot \mathbf{n}
- \frac{1}{2} s_{\max} \left( \q_h^+ - \q_h^- \right).
\label{eq.rusanov} 
\end{equation} 
Furthermore, we also need to develop a proper discretization of the non conservative products at the element boundaries. 
To this end, we consider the works \cite{Pares2004,Castro2006,Pares2006,MunozPares}, 
based on the theory of Dal Maso, Le Floch and Murat \cite{DLMtheory} on nonconservative 
hyperbolic PDE systems, and their extension to higher order DG schemes in \cite{Rhebergen2008,ADERNC}. 
Within this framework it is usual to build also so-called well-balanced schemes  \cite{Bermudez1994,GarciaNavarro1,LeVequeWB}. 
Accordingly, the third term in \eqref{eq:ADER-DG} is approximated with a path integral in phase-space between the two extrapolated values related to the face, $\q_h^-$ and $\q_h^+$, 
\begin{equation}
\mathcal{D}\left( \q_h^-,\q_h^+ \right) \cdot \mathbf{n} =  \frac{1}{2}
\left(\int \limits_{0}^{1} \mathbf{\mathcal{B}} \left(
\bm{\psi}(\q_h^-,\q_h^+,s) \right)\cdot\mathbf{n} \, ds
\right)\cdot\left(\q_h^+ - \q_h^-\right), \label{eq:PC}
\end{equation}
where we have used the linear segment path
\begin{equation}
\bm{\psi} = \bm{\psi}(\q_h^-, \q_h^+, s) = \q_h^- + s \left( \q_h^+ -
\q_h^- \right)\,, \qquad s \in [0,1].
\end{equation}
As it can be seen in \cite{MPT13,ACD18}, different paths could have been chosen to perform the former integral attending to special features of non conservative and source terms. Nevertheless, for the systems addressed in this work, the easiest straight line path already provides good results.

From the computational point of view, it is important to remark that, as a consequence of the nodal tensor-product basis employed, the scheme can be written in a dimension by dimension fashion and integral operators are decomposed as the product of one-dimensional operators \cite{BCDGP20}.
The resulting explicit one-step ADER-DG scheme is conditionally stable with stability condition
\begin{equation}
\Delta t < \frac{\mathrm{CFL }\, h_{\mathrm{min}}}{\left(2N+1\right)\left|\lambda_{\mathrm{max}}\right|},
\label{eq:cfl}
\end{equation}
with $\textnormal{CFL}<1/d$, where $\left|\lambda_{\mathrm{max}}\right|$ is the maximum eigenvalue of the system, $h_{\mathrm{min}}$ denotes the minimum characteristic mesh size and $d$ is the number of space dimensions. 
The final time step for the coupled problem is the lowest one between the timesteps associated to the two considered models.

\subsection{Boundary conditions}
\begin{figure}[h]
	\centering
	\includegraphics[width=0.45\linewidth]{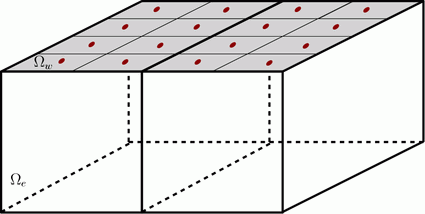}
	\begin{minipage}{0.85\linewidth}
		\caption{Sketch of the nonconforming mesh at the interface, $\partial \Omega_{e,w}$. The 2D cells (grey) filling the fluid domain, $\Omega_w$, are built as a subgrid of the faces of the 3D elements discretizing the solid domain, $\Omega_e$.}	\label{fig:coupling}
	\end{minipage}
\end{figure}

A major challenge concerning the coupling of the two models is the large discrepancy, in the spatial length 
scales, between the elastic waves, with a typical propagation speed of about $3000-4000$m/s and the short free surface water waves, whose propagation speed is between $10$m/s and $30$m/s. The corresponding difference in the expected wave lengths is huge, with wave lengths of the order of kilometres in solid media and well below $10^{2}$ m in water. To address this problem, we employ two non-conforming meshes with different spatial resolutions. An initial 3D mesh made of rectangular bricks is first designed to cover $\Omega_e$. Then, on $\Omega_w$, a much finer grid is built as a refinement of the 2D mesh made up by the faces of the mesh designed in $\Omega_e$, lying on the boundary $\partial \Omega_{e,w}$ (see Figure \ref{fig:coupling}). 
Interpolation of the discrete solution from $\Omega_w$ onto the boundary $\partial \Omega_{e,w}$ is
carried out using appropriate high order Gaussian quadrature formulas, \cite{stroud}, so that the boundary condition \eqref{eq:boundaryconditions} can be imposed on the solid domain. 

Further boundary conditions for the solid media include periodic and free surface boundaries. Regarding the last ones, the exact Riemann solver of Godunov can be employed, \cite{gij1,gij2,TDCRWB19}. Accordingly, the opposite value for the incoming normal stress to the boundary is set.
Finally, on the fluid domain we consider either periodic or Dirichlet boundary conditions.

\section{Numerical test problems}\label{sec:num_res}
This section is devoted to the assessment of the developed methodology. 
Initially, we address the two mathematical models independently, analysing the solution obtained 
for classical benchmarks of linear elasticity and non-hydrostatic flows.
Once the numerical method is validated, 
we present two showcases of the coupled problem, reporting the seismic wave propagation 
generated by a soliton and by a sinusoidal wave train on the water surface.

\subsection{Linear elastic wave problems}
Linear elasticity is a well established research field, so many numerical tests can be found to validate the proposed methodology. In what follows, we will first validate the numerical method using a p-s-wave test, whose exact periodic solution is known. Then, we focus on a classical benchmark of seismic wave propagation, the so-called Lamb problem, and on a stiff inclusion test, which accounts for large material parameter variations. A final wave propagation test in three dimensions is also included.

\subsubsection{Numerical convergence test}

Following \cite{gij1,TD18LE}, a $p-$ and $s-$wave test case is employed to verify the accuracy of the ADER-DG scheme.
We consider the computational domain  $\Omega = [-1.5,1.5] \times [-1.5,1.5]\in\mathbb{R}^{2}$ with periodic boundary conditions in $x$ and $y$ directions and we define the initial condition
\begin{equation}
\mathbf{Q}\left(\mathbf{x},0\right) = \alpha \mathbf{r}_{p} \sin\left(2\pi \mathbf{n}\cdot \mathbf{x}\right) +\alpha \mathbf{r}_{s} \sin\left(2\pi \mathbf{n}\cdot \mathbf{x}\right),
\end{equation}
with $\alpha= 0.1$, $\mathbf{n}=\left(n_{x},n_{y}\right)=\left(1,1\right)$, $\mathbf{r}_{p},\, \mathbf{r}_{s}$ the eigenvectors associated to the $p-$ and $s-$waves,
\begin{equation}
\mathbf{r}_{p} = \left(\lambda+2\mu n_x^{2},\lambda+2\mu n_y^{2}, 2\mu n_{x} n_{y}, -c_{p}n_{x}, -c_{p}n_{y}\right), \quad
\mathbf{r}_{s} = \left(-2\mu n_x n_y,2\mu n_x n_y, \mu\left(n_{x}^{2}-n_{y}^{2}\right), c_{s}n_{y}, -c_{s}n_{x}\right),
\end{equation}
$c_{p}=\sqrt{\left(\lambda+2\mu\right)/\rho}$ the $p-$wave speed and $c_{s}=\sqrt{\mu/\rho}$ the $s-$wave speed. Setting the material parameters to $\lambda=2$, $\mu=1$, $\rho=1$, leads to propagation velocities, $c_{p}=2$ and $c_{s}=1$.
The former initial condition generates a sinusoidal $p-$wave travelling along the diagonal direction, $\mathbf{n}=\left(1,1\right)$, and another sinusoidal $s-$wave moving in the opposite direction. 
Taking $t_{\mathrm{end}}=3\sqrt{2}$, the solution coincides with the initial datum and a convergence analysis can be performed. 
Table \ref{tab:LE_accuracy} shows the $L^{2}$ errors, $\epsilon_{L^{2}}$, and the convergence order, $\mathcal{O}$, obtained for $N\in\left\lbrace 3,4,5,6,7 \right\rbrace$, $M=N$. The spatial grids were built using the same number of elements in $x$ and $y$ directions, $N_{y}=N_{x}$, and the time step is computed at each iteration according to the CFL condition. All variables attain the optimal order of convergence sought.

\begin{table}
	\hspace*{-1.5cm}
	\begin{tabular}{|c|c||c|c|c|c|c||c|c|c|c|c||c|}
		\hline
		$N$ & $N_{x}$ & $\epsilon_{L^{2}}\!\left(\!u\!\right) $ & $\epsilon_{L^{2}}\!\left(\!v\!\right) $ & $\epsilon_{L^{2}}\!\left(\!\sigma_{x\!x}\!\right) $& $\epsilon_{L^{2}}\!\left(\!\sigma_{y\!y}\!\right) $& $\epsilon_{L^{2}}\left(\!\sigma_{x\!y}\!\right) $& $\mathcal{O}\!\left(\!u\!\right) $ & $\mathcal{O}\!\left(\!v\!\right) $ & $\mathcal{O}\!\left(\!\sigma_{x\!x}\!\right) $ & $\mathcal{O}\!\left(\!\sigma_{y\!y}\!\right) $ & $\mathcal{O}\!\left(\!\sigma_{x\!y}\!\right) $ & Teor.\\
		\hline\hline
		\multirow{5}{*}{$3$} &$10$ & $1.18E\!-\!02$ & $1.66E\!-\!02$ & $2.70E\!-\!02$ & $3.55E\!-\!02$ & $1.07E\!-\!02$  & - & - & - &  - &  - &\multirow{5}{*}{$4$}  \\	\hhline{|~|-||-|-|-|-|-||-|-|-|-|-||~|}
		&$15$ & $2.24E\!-\!03$ & $3.29E\!-\!03$ & $5.17E\!-\!03$ & $6.57E\!-\!03$ & $1.98E\!-\!03$ & $4.10$ & $3.99$ & $4.08$ & $4.16$ & $4.16$ & \\	\hhline{|~|-||-|-|-|-|-||-|-|-|-|-||~|}
		&$20$ & $7.02E\!-\!04$ & $1.07E\!-\!03$ & $1.63E\!-\!03$ & $2.06E\!-\!03$ & $6.18E\!-\!04$ & $4.03$ & $3.92$ & $4.00$ & $4.04$ & $4.05$ & \\	\hhline{|~|-||-|-|-|-|-||-|-|-|-|-||~|}
		&$30$ & $1.38E\!-\!04$ & $2.15E\!-\!04$ & $3.26E\!-\!04$ & $4.07E\!-\!04$ & $1.21E\!-\!04$ & $4.02$ & $3.95$ & $3.98$ & $4.00$ & $4.02$ & \\	\hhline{|~|-||-|-|-|-|-||-|-|-|-|-||~|}
		&$40$ & $4.34E\!-\!05$ & $6.85E\!-\!05$ & $1.05E\!-\!04$ & $1.30E\!-\!04$ & $3.81E\!-\!05$ & $4.01$ & $3.97$ & $3.94$ & $3.97$ & $4.01$ & \\	\hline\hline
		\multirow{5}{*}{$4$}& $10$ & $7.67E\!-\!04$ & $1.15E\!-\!03$ & $1.78E\!-\!03$ & $2.20E\!-\!03$ & $6.73E\!-\!04$  & - & - & - &  - &  - &\multirow{5}{*}{$5$} \\	\hhline{|~|-||-|-|-|-|-||-|-|-|-|-||~|}
		& $15 $ & $1.09E\!-\!04$ & $1.53E\!-\!04$ & $2.54E\!-\!04$ & $3.14E\!-\!04$ & $9.52E\!-\!05$ & $4.81$ & $4.99$ & $4.79$ & $4.80$ & $4.82$ &   \\	\hhline{|~|-||-|-|-|-|-||-|-|-|-|-||~|}
		& $20$ & $2.66E\!-\!05$ & $3.61E\!-\!05$ & $6.29E\!-\!05$ & $7.74E\!-\!05$ & $2.32E\!-\!05$ & $4.90$ & $5.01$ & $4.85$ & $4.87$ & $4.90$ &   \\	\hhline{|~|-||-|-|-|-|-||-|-|-|-|-||~|}
		& $30$ & $3.58E\!-\!06$ & $4.72E\!-\!06$ & $8.74E\!-\!06$ & $1.06E\!-\!05$ & $3.12E\!-\!06$ & $4.94$ & $5.01$ & $4.87$ & $4.89$ & $4.95$ &   \\	\hhline{|~|-||-|-|-|-|-||-|-|-|-|-||~|}
		& $40$ & $8.56E\!-\!07$ & $1.12E\!-\!06$ & $2.15E\!-\!06$ & $2.60E\!-\!06$ & $7.45E\!-\!07$ & $4.98$ & $5.01$ & $4.87$ & $4.90$ & $4.98$ &   \\	\hline\hline
		\multirow{5}{*}{$5$}& $10$ & $4.61E\!-\!05$ & $6.38E\!-\!05$ & $1.09E\!-\!04$ & $1.35E\!-\!04$ & $4.14E\!-\!05$ & - & - & - &  - &  - & \multirow{5}{*}{$6$} \\	\hhline{|~|-||-|-|-|-|-||-|-|-|-|-||~|}
		& $15 $ & $4.04E\!-\!06$ & $6.00E\!-\!06$ & $9.84E\!-\!06$ & $1.21E\!-\!05$ & $3.64E\!-\!06$ & $6.00$ & $5.83$ & $5.92$ & $5.94$ & $5.99$ &    \\	\hhline{|~|-||-|-|-|-|-||-|-|-|-|-||~|}
		& $20 $ & $7.16E\!-\!07$ & $1.10E\!-\!06$ & $1.81E\!-\!06$ & $2.21E\!-\!06$ & $6.46E\!-\!07$ & $6.02$ & $5.90$ & $5.88$ & $5.92$ & $6.01$ &    \\	\hhline{|~|-||-|-|-|-|-||-|-|-|-|-||~|}
		& $30 $ & $6.25E\!-\!08$ & $9.84E\!-\!08$ & $1.71E\!-\!07$ & $2.04E\!-\!07$ & $5.64E\!-\!08$ & $6.01$ & $5.95$ & $5.82$ & $5.88$ & $6.01$ &    \\	\hhline{|~|-||-|-|-|-|-||-|-|-|-|-||~|}
		& $40$ & $1.11E\!-\!08$ & $1.76E\!-\!08$ & $3.25E\!-\!08$ & $3.80E\!-\!08$ & $1.00E\!-\!08$ & $6.00$ & $5.97$ & $5.77$ & $5.84$ & $4.04$ &    \\	\hline\hline
		\multirow{5}{*}{$6$}& $10$ & $2.72E\!-\!06$ & $4.10E\!-\!06$ & $6.53E\!-\!06$ & $8.01E\!-\!06$ & $2.46E\!-\!06$ & - & - & - &  - &  - & \multirow{5}{*}{$7$} \\	\hhline{|~|-||-|-|-|-|-||-|-|-|-|-||~|}
		& $15 $ & $1.67E\!-\!07$ & $2.36E\!-\!07$ & $4.18E\!-\!07$ & $5.05E\!-\!07$ & $1.50E\!-\!07$ & $6.88$ & $7.05$ & $6.78$ & $6.82$ & $6.90$ &   \\	\hhline{|~|-||-|-|-|-|-||-|-|-|-|-||~|}
		& $20 $ & $2.28E\!-\!08$ & $3.11E\!-\!08$ & $5.97E\!-\!08$ & $7.11E\!-\!08$ & $2.04E\!-\!08$ & $6.93$ & $7.04$ & $6.77$ & $6.82$ & $6.94$ &   \\	\hhline{|~|-||-|-|-|-|-||-|-|-|-|-||~|}
		& $30 $ & $1.36E\!-\!09$ & $1.80E\!-\!09$ & $3.87E\!-\!09$ & $4.50E\!-\!09$ & $1.21E\!-\!09$ & $6.96$ & $7.03$ & $6.75$ & $6.81$ & $6.97$ & \\	\hhline{|~|-||-|-|-|-|-||-|-|-|-|-||~|}
		& $40 $ & $1.82E\!-\!10$ & $2.38E\!-\!10$ & $5.55E\!-\!10$ & $6.35E\!-\!10$ & $1.62E\!-\!10$ & $6.99$ & $7.03$ & $6.75$ & $6.81$ & $6.99$ & \\	\hline\hline
		\multirow{5}{*}{$7$}& $5$ & $3.31E\!-\!05 $ & $4.50E\!-\!05 $& $7.76E\!-\!05 $ & $9.58E\!-\!05 $ & $3.07E\!-\!05 $ & - & - & - &  - &  - &  \multirow{5}{*}{$8$} \\	\hhline{|~|-||-|-|-|-|-||-|-|-|-|-||~|}
		& $10 $ & $1.40E\!-\!07 $ & $1.96E\!-\!07 $& $3.47E\!-\!07 $ & $4.23E\!-\!07 $ & $1.29E\!-\!07 $ & $7.88 $ & $7.84 $ & $7.80 $ & $7.82 $ & $7.90 $ &   \\	\hhline{|~|-||-|-|-|-|-||-|-|-|-|-||~|}
		& $15 $ & $5.42E\!-\!09 $ & $8.10E\!-\!09 $& $1.45E\!-\!08 $ & $1.73E\!-\!08 $ & $5.01E\!-\!09 $ & $8.02 $ & $7.86 $ & $7.83 $ & $7.88 $ & $8.01 $ &   \\	\hhline{|~|-||-|-|-|-|-||-|-|-|-|-||~|}
		& $20 $ & $5.38E\!-\!10 $ & $8.28E\!-\!10 $& $1.56E\!-\!09 $ & $1.82E\!-\!09 $ & $4.98E\!-\!10 $ & $8.03$ & $7.92$ & $7.75$ & $7.82$ & $8.03$ &   \\	\hhline{|~|-||-|-|-|-|-||-|-|-|-|-||~|}
		& $25 $ & $9.02E\!-\!11 $ & $1.41E\!-\!10 $& $2.81E\!-\!10 $ & $3.22E\!-\!10 $ & $8.34E\!-\!11 $ & $8.00$ & $7.93$ & $7.69$ & $7.77$ & $8.00$ &   \\	\hline
	\end{tabular}
	\caption{Numerical convergence test for linear elasticity. $L^{2}$ errors and convergence rates obtained for $N\in\left\lbrace 3,4,5,6,7\right\rbrace$ at $t_{\mathrm{end}}=3\sqrt{2}$.}
	\label{tab:LE_accuracy}
\end{table}

\subsubsection{2D Lamb's problem}
\begin{figure}
	\centering
	\includegraphics[width=0.7\linewidth]{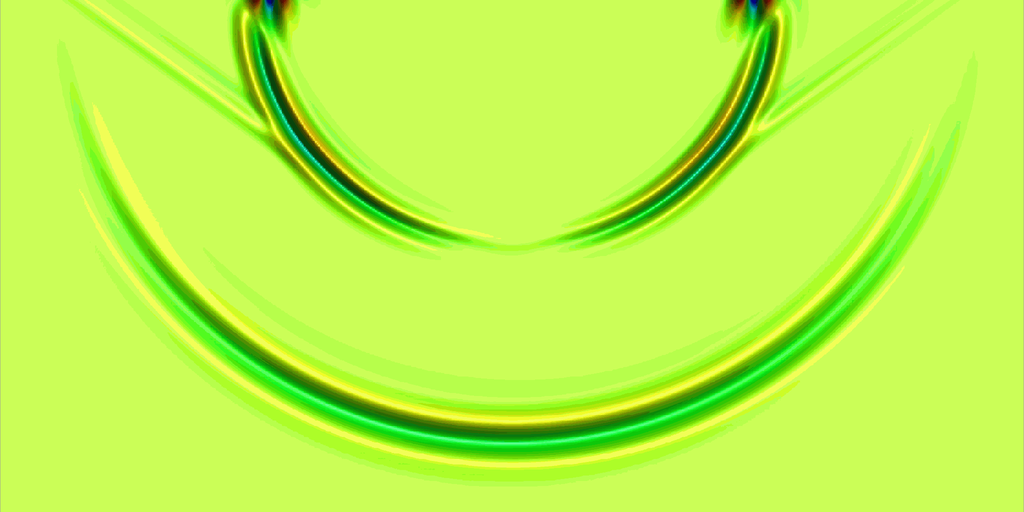}
	
	\vspace{0.05\linewidth}
	\includegraphics[width=0.7\linewidth]{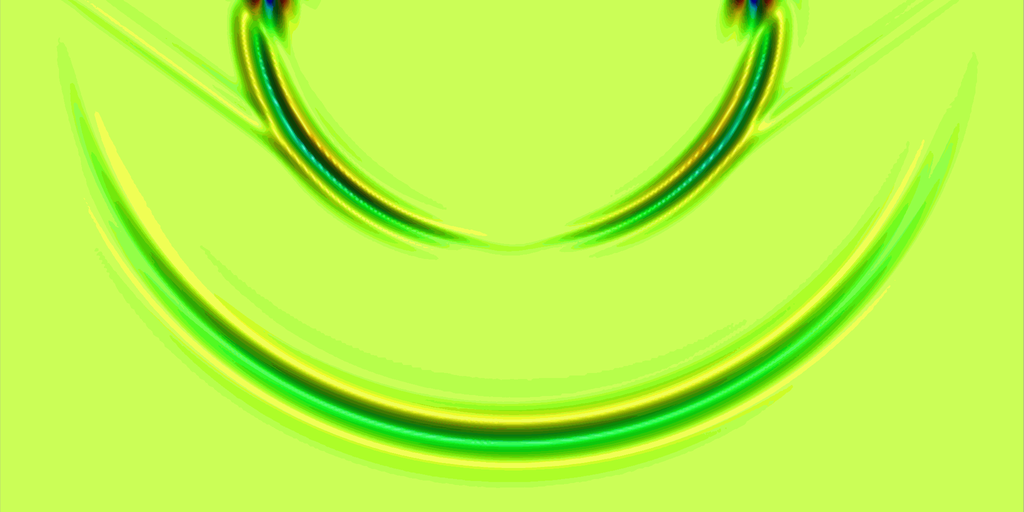}
	
	\caption{2D Lamb's problem. Vertical velocity contour plot at $t=0.6$. Top: ADER-DG $\mathcal{O}4$ run on a mesh made of $200\times100$ elements. Bottom: reference solution computed with a $P_{4}P_{4}$ ADER-DG scheme on an unstructured mesh of $N_{i}=180276$ triangles.}
	\label{fig:2DLamb_verticalvelocity}
\end{figure}

\begin{figure}
	\centering
	\includegraphics[width=0.3\linewidth]{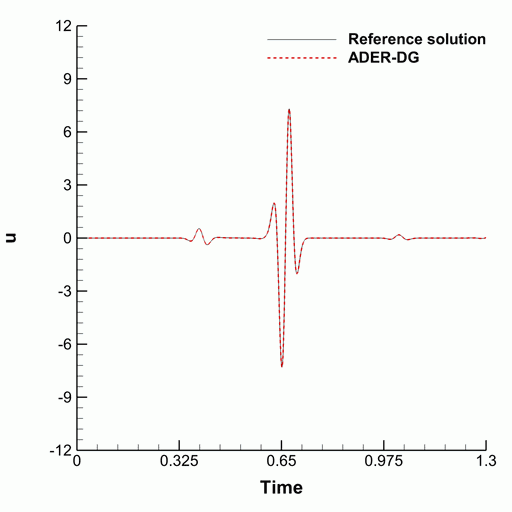}\hspace{0.05\linewidth}
	\includegraphics[width=0.3\linewidth]{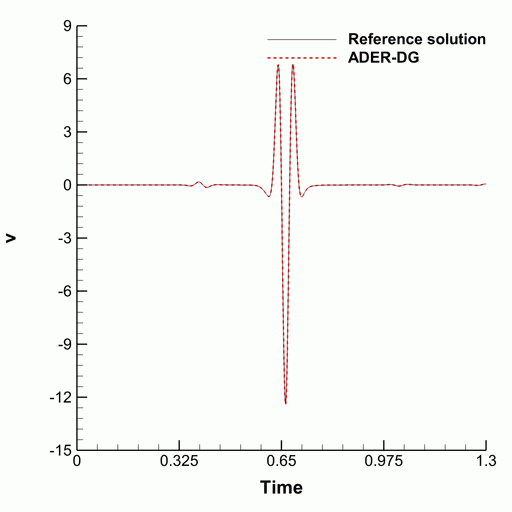}
	
	\includegraphics[width=0.3\linewidth]{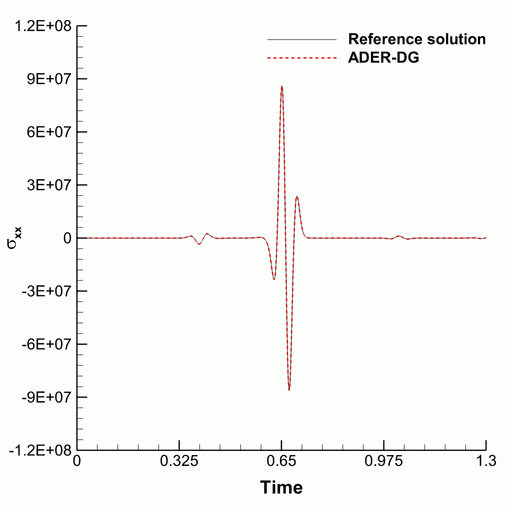}\hspace{0.02\linewidth}	
	\includegraphics[width=0.3\linewidth]{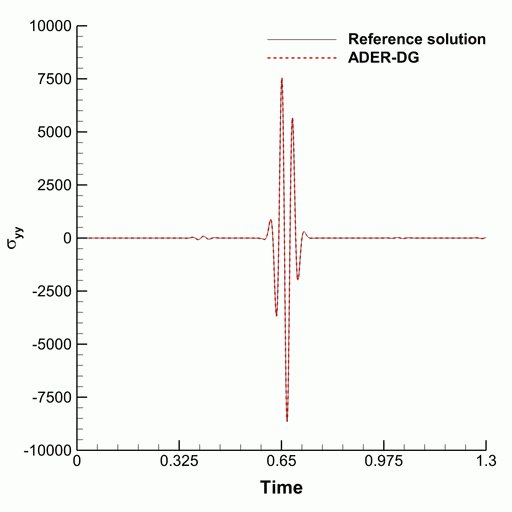}\hspace{0.02\linewidth}
	\includegraphics[width=0.3\linewidth]{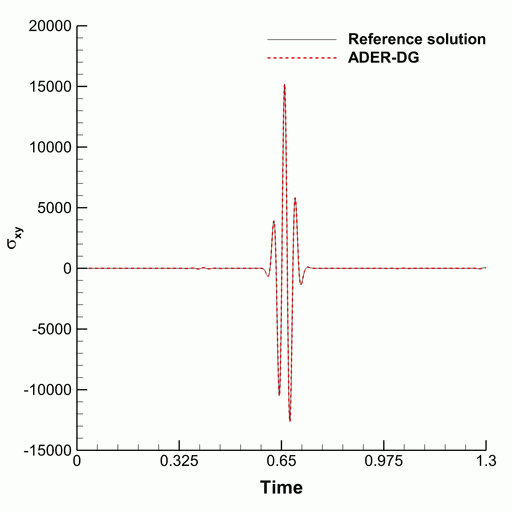}
	
	\caption{2D Lamb's problem. Seismogram at receiver $\mathbf{x}_{1}= \left(990,0\right)$ obtained using ADER-DG $\mathcal{O}4$ on a mesh made of $200\times100$ elements and reference solution, computed with a $P_{4}P_{4}$ ADER-DG scheme on an unstructured mesh of $N_{i}=180276$ triangles.}
	\label{fig:2DLamb_pickpoint990}
\end{figure}

The Lamb's problem is a well known benchmark used to test numerical methods for linear elastic waves and has first been put forward in \cite{Lamb1904}. In this paper we consider one of its variants, already analysed in \cite{gij1,GPRmodel}, where a rectangular domain $\Omega = \left[-2000,2000\right]\times \left[-2000,0\right]$ and a point source of the form 
\begin{equation}
\mathbf{f}_{\mathbf{v}} \left(\mathbf{x},t\right) = \rho_{s} a_{1} \left(\frac{1}{2} + a_{2}\left(t-t_{D} \right)^{2} \right) \exp \left(a_{2}\left(t-t_{D} \right)^{2} \right) \delta\left(\mathbf{x}- \mathbf{x}_{s}\right) \mathbf{e}_{y} \label{eq:lamb_source}
\end{equation}
in the momentum equation, \eqref{eqn:momentum}, are chosen. 
We locate the source near the free surface using the Dirac delta distribution at $\mathbf{x}_{s}=\left(0,-1\right)$ and we set the related parameters to $\rho_{s}=2200$, $a_{1}=-2000$, $a_{2}= -\left(\pi f_{c}\right)^{2}$, $f_{c}=14.5$, $t_{D}= 0.08$.  The considered homogeneous material has density $\rho=2200$ and Lam\'e constants $\lambda = 7509672500$, $\mu = 7509163750$, so the propagation velocities are $c_{p}=3200$ and $c_{s}=1847.5$. The domain is discretized using a Cartesian grid made of $200\times100$ elements and the simulation is run for $N=3$. 
The vertical velocity contour plot at time $t=0.6$ is depicted in Figure \ref{fig:2DLamb_verticalvelocity}. We observe a good agreement with the reference solution that has been computed by the \texttt{SeisSol} 
community code\footnote{\textcolor{blue}{www.seissol.org}} using the $P_{4}P_{4}$ ADER-DG scheme, presented in \cite{gij1,gij2,SeisSol1,SeisSol2}, on an unstructured mesh made of $N_{i}=180276$ triangles. In Figure \ref{fig:2DLamb_pickpoint990} we observe an almost perfect agreement between the two simulations in the seismogram obtained at a receiver located in $\mathbf{x}_{1}= \left(990,0\right)$. 

\subsubsection{Stiff inclusion}
\begin{figure}
	\centering
	\includegraphics[width=0.9\linewidth,trim= 10 0 10 0,clip]{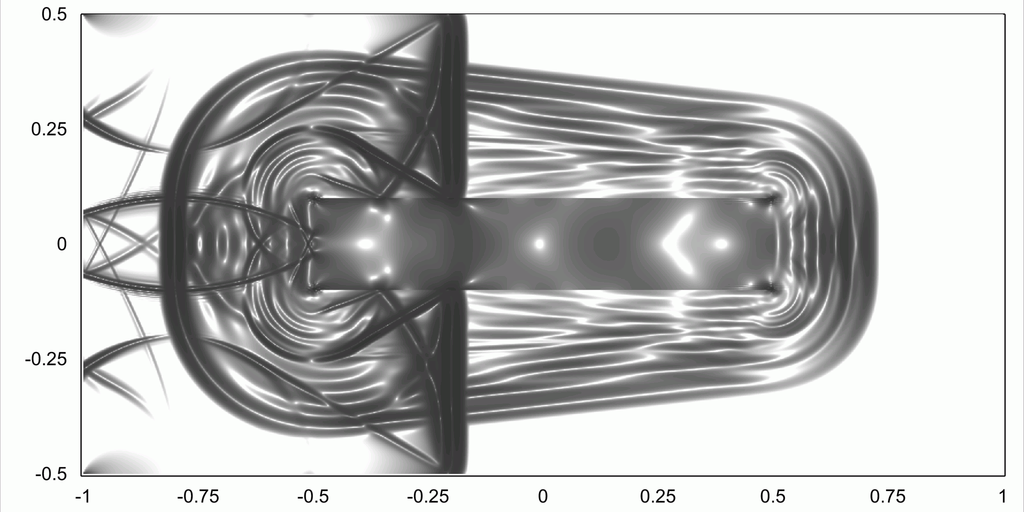}
	
	\vspace{0.05\linewidth}
	\includegraphics[width=0.9\linewidth,trim= 10 0 10 0,clip]{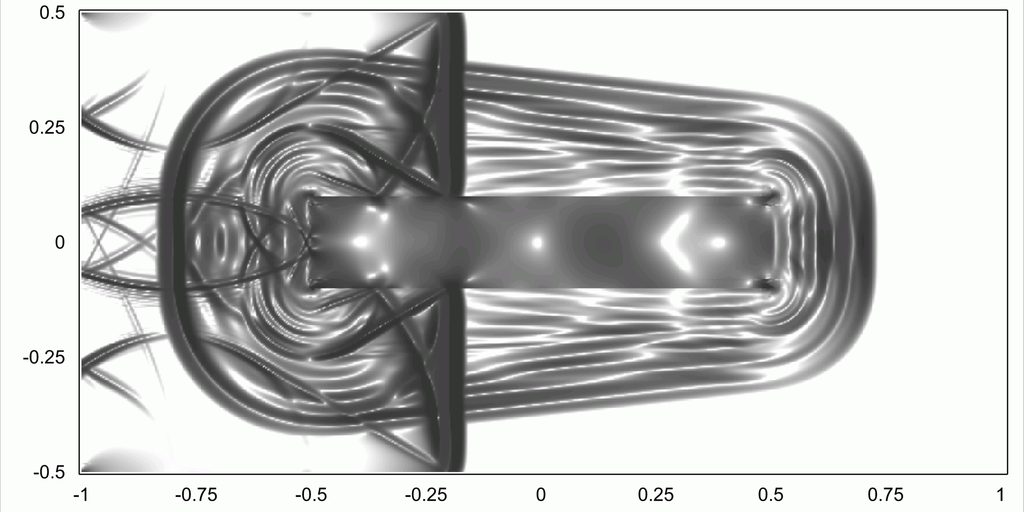}
	
	\caption{Stiff inclusion. Component $\sigma_{x x}$ of the wavefield at time $t=0.3$. Top: ADER-DG $\mathcal{O}4$ on a mesh made of $400\times 100$ elements. Bottom ADER-DG $\mathcal{O}7$ on a mesh made of $100\times 50$ elements.} 
	\label{fig:SI_t30}
\end{figure}

To study the behaviour of the method under large jumps of material parameters, we consider the stiff inclusion benchmark, \cite{LeVeque:2002a,gij1}. 
The computational domain, $\Omega=[-1,1]\times[-0.5,0.5]$,  is divided into two regions with different materials. The outer material properties are $\lambda_{\mathrm{out}}=2$, $\mu_{\mathrm{out}}=1$ and $\rho_{\mathrm{out}}=1$, whereas the inner material, placed in $\Omega_{\mathrm{in}}=[-0.5,0.5]\times[-0.1,0.1]$, has $\lambda_{\mathrm{in}}=200$, $\mu_{\mathrm{in}}=100$ and $\rho_{\mathrm{in}}=1$. The initial field is characterized by a p-wave of the form
\begin{equation}
\mathbf{Q}\left(\mathbf{x},0\right) =  \mathbf{r}_{p} \exp\left[ -\frac{1}{2 \sigma^{2}} \left( \mathbf{n}\cdot \left( \mathbf{x}-\mathbf{x}_{0}\right)\right) ^{2} \right],
\end{equation}
with $\mathbf{x}_{0}=\left(-0.08,0 \right) $ the initial wave location, $\mathbf{n}=\left(1,0\right) $ and, $\sigma=0.01$ the standard deviation. Free surface boundary conditions are considered, so that the normal and shear stresses vanish at the boundary. Consequently, surface waves will develop from the beginning of the simulation. At time $t=0.15$ the planar wave reaches the stiffer material, where elastic waves propagate ten times faster that in the outer media. Then, the reflection of waves inside the inclusion yield to its vibration, which results in small amplitude waves propagating into the outer domain. 
The $\sigma_{x x}$ pattern generated at time $t=0.3$ is plotted in Figure~\ref{fig:SI_t30} for two different simulations.
The solution obtained for the $\mathcal{O}4$ ADER-DG scheme on a mesh made of $400\times 100$ elements agrees well with the one obtained using the  $\mathcal{O}7$ ADER-DG scheme on a much coarser mesh, $100\times 50$ elements. To compare the results with data available in the bibliography one may refer to \cite{LeVeque:2002a} as well as to \cite{gij1}. Overall, one can note a very good qualitative agreement of the wavefield with the numerical reference solutions available in the literature.  

\subsubsection{3D wave propagation}
\begin{figure}
	\centering
	\includegraphics[width=0.5\linewidth]{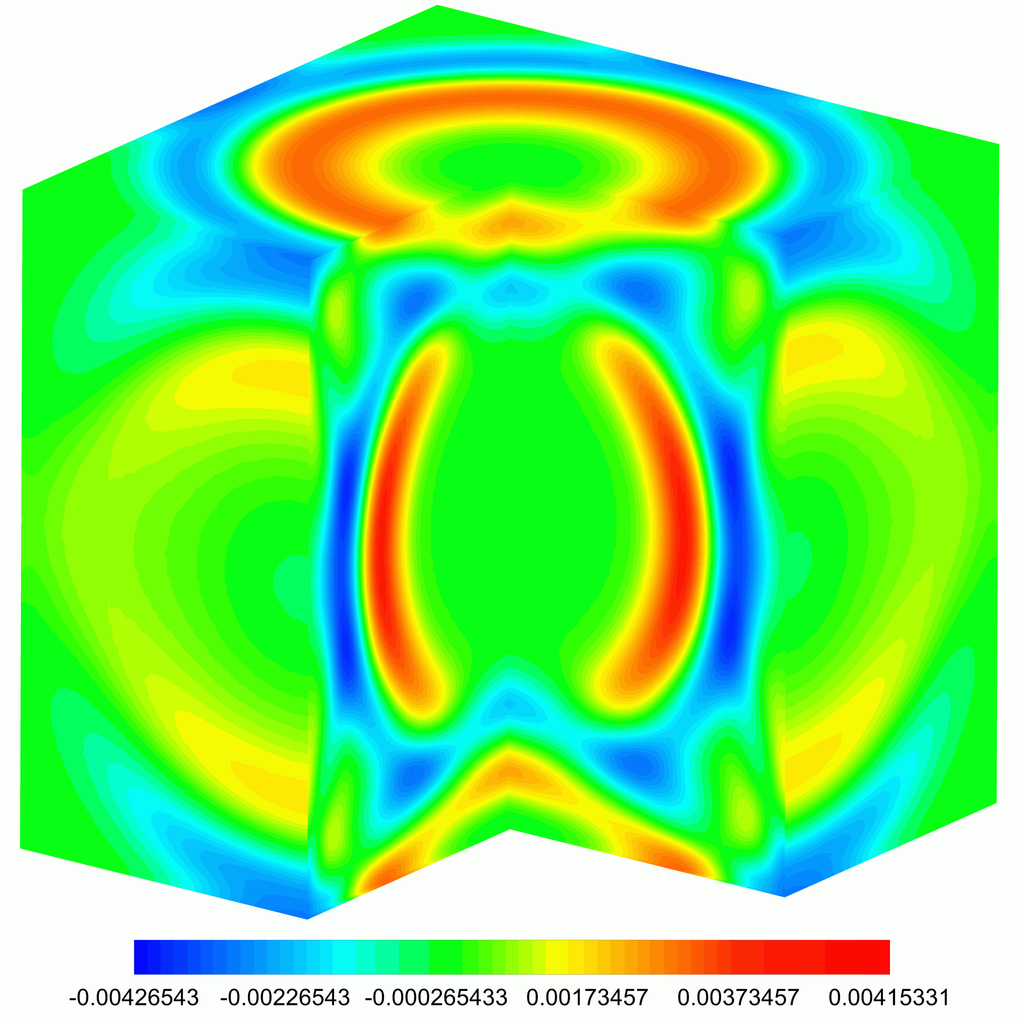}
	
	\caption{3D wave propagation. Vertical velocity contour plot at $t=2$ obtained on a mesh made of $30\times 30 \times 30$ elements using the fourth order ADER-DG scheme.}
	\label{fig:3DExp_w_t2}
\end{figure}

\begin{figure}
	\centering
	\includegraphics[width=0.32\linewidth]{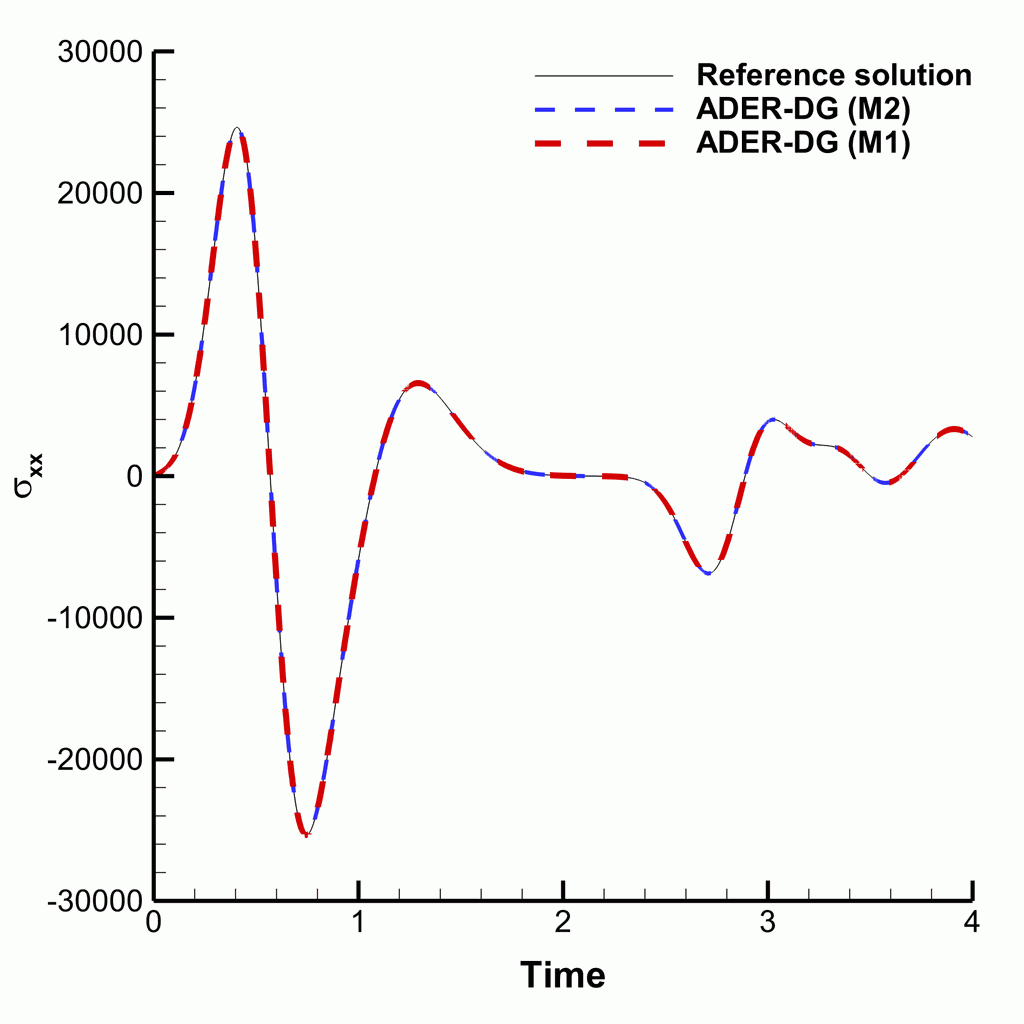}\hfill
	\includegraphics[width=0.32\linewidth]{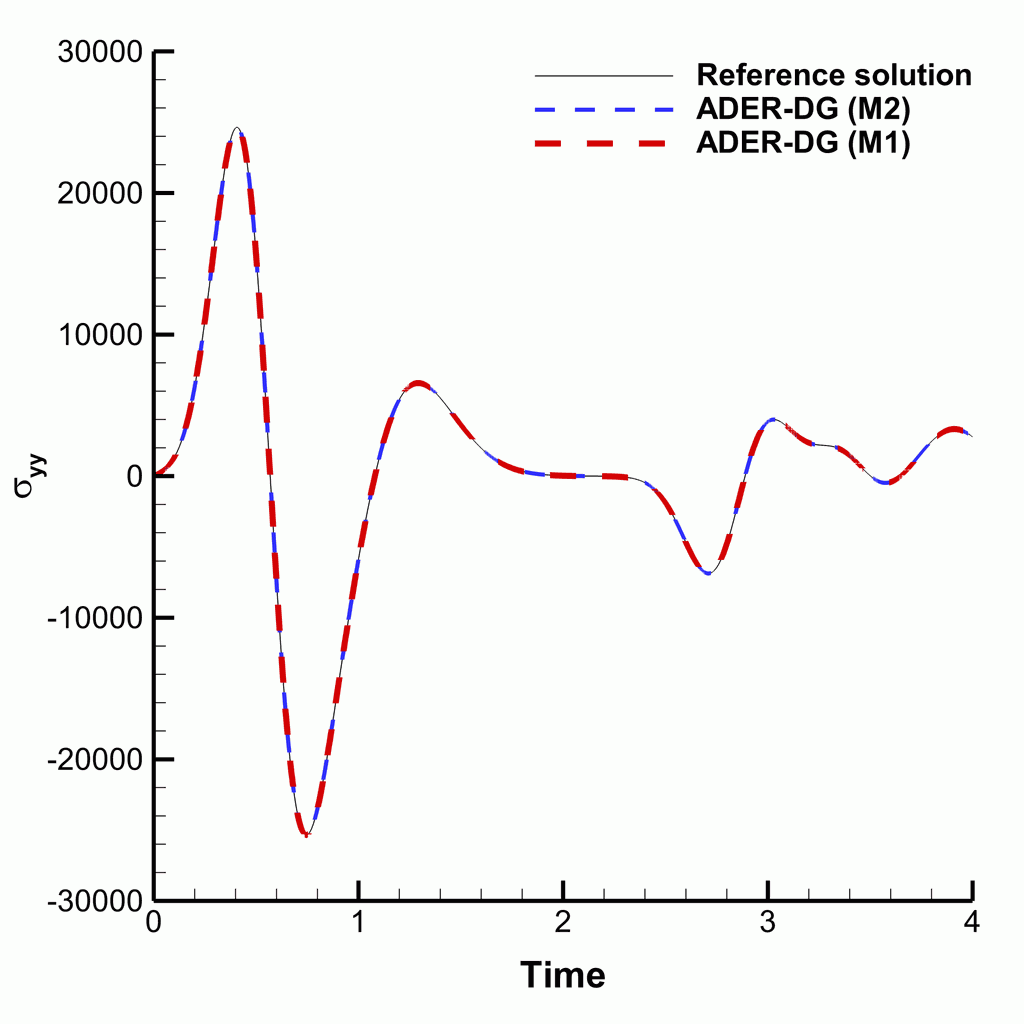}\hfill
	\includegraphics[width=0.32\linewidth]{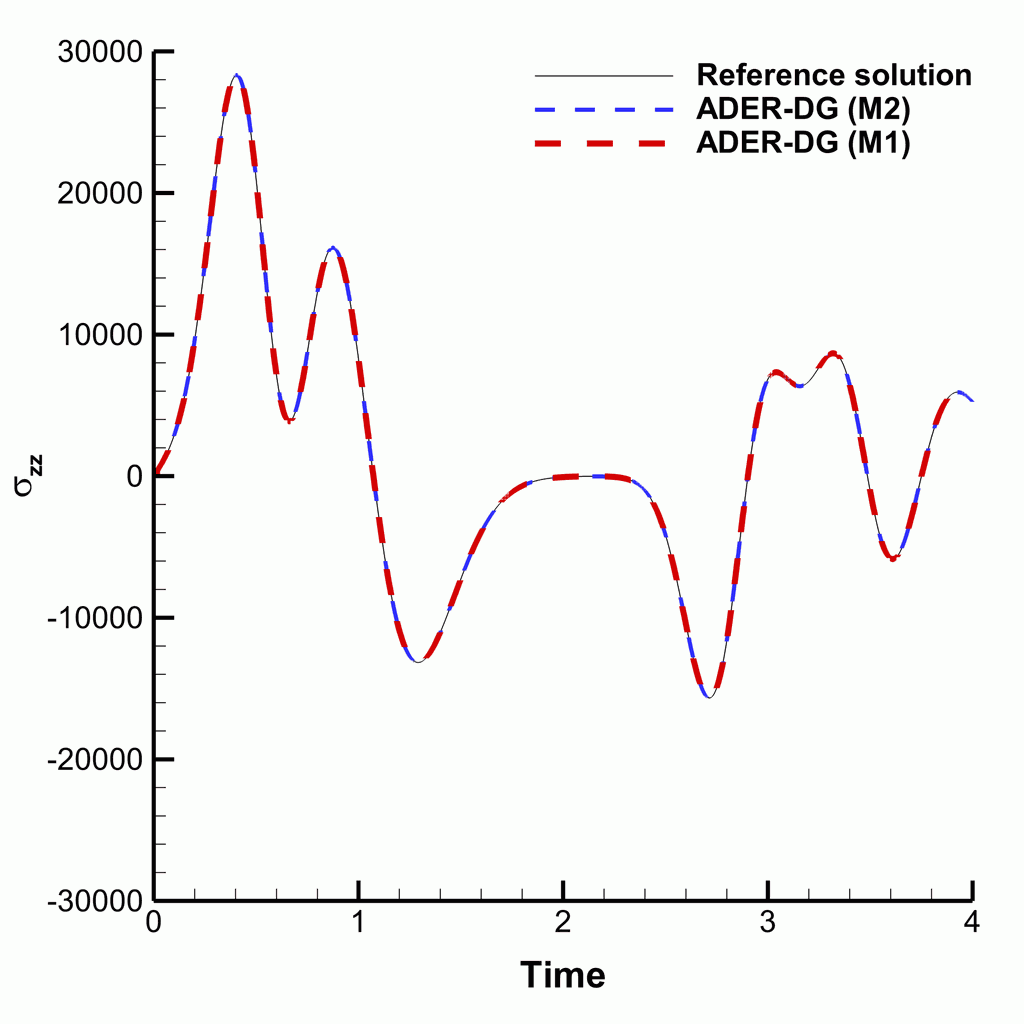}\\	
	\includegraphics[width=0.32\linewidth]{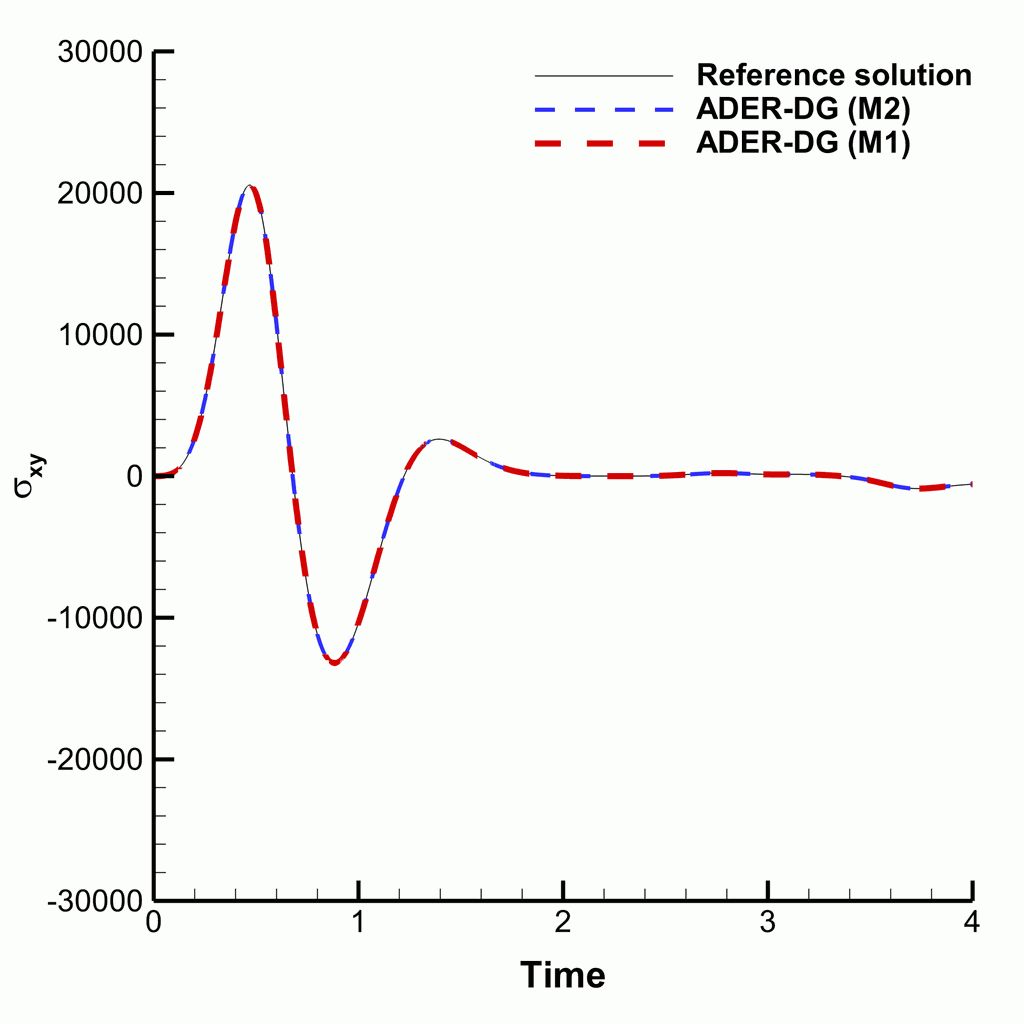}\hfill
	\includegraphics[width=0.32\linewidth]{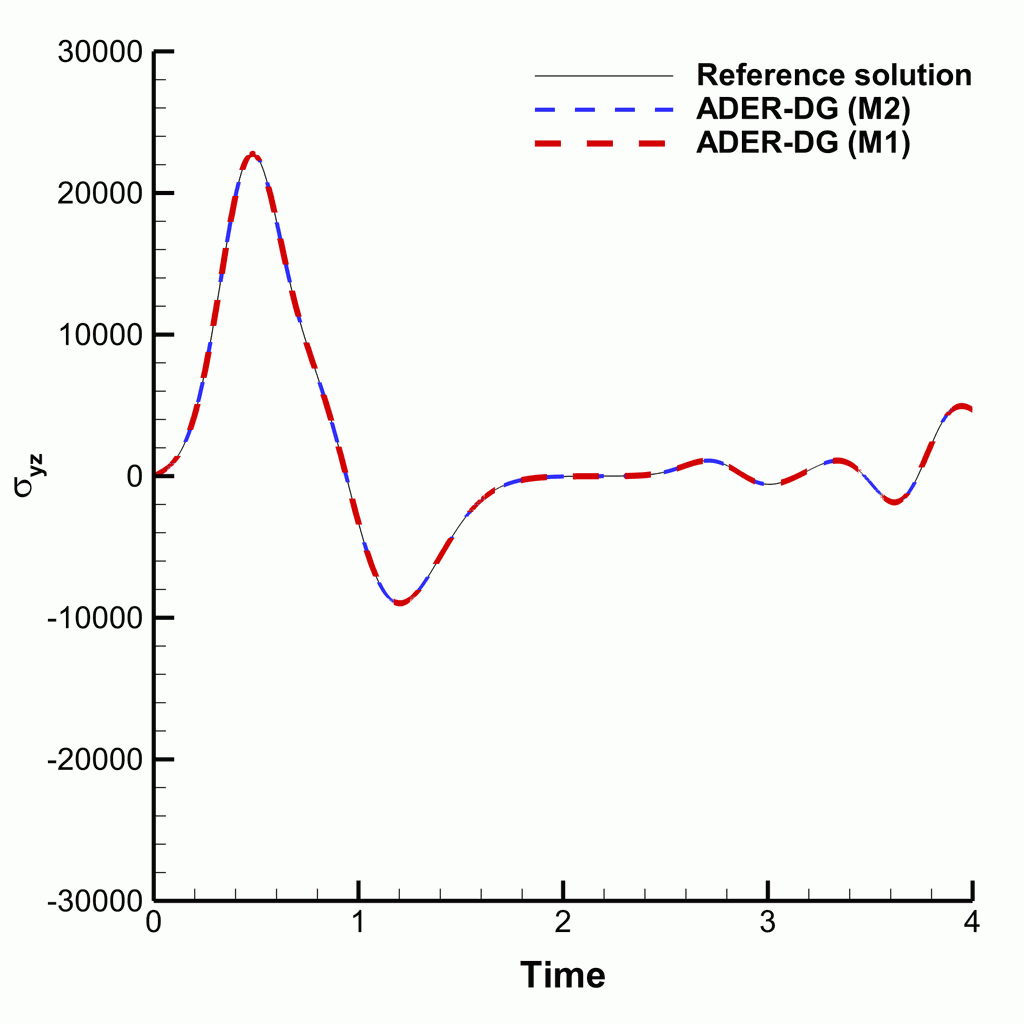}\hfill
	\includegraphics[width=0.32\linewidth]{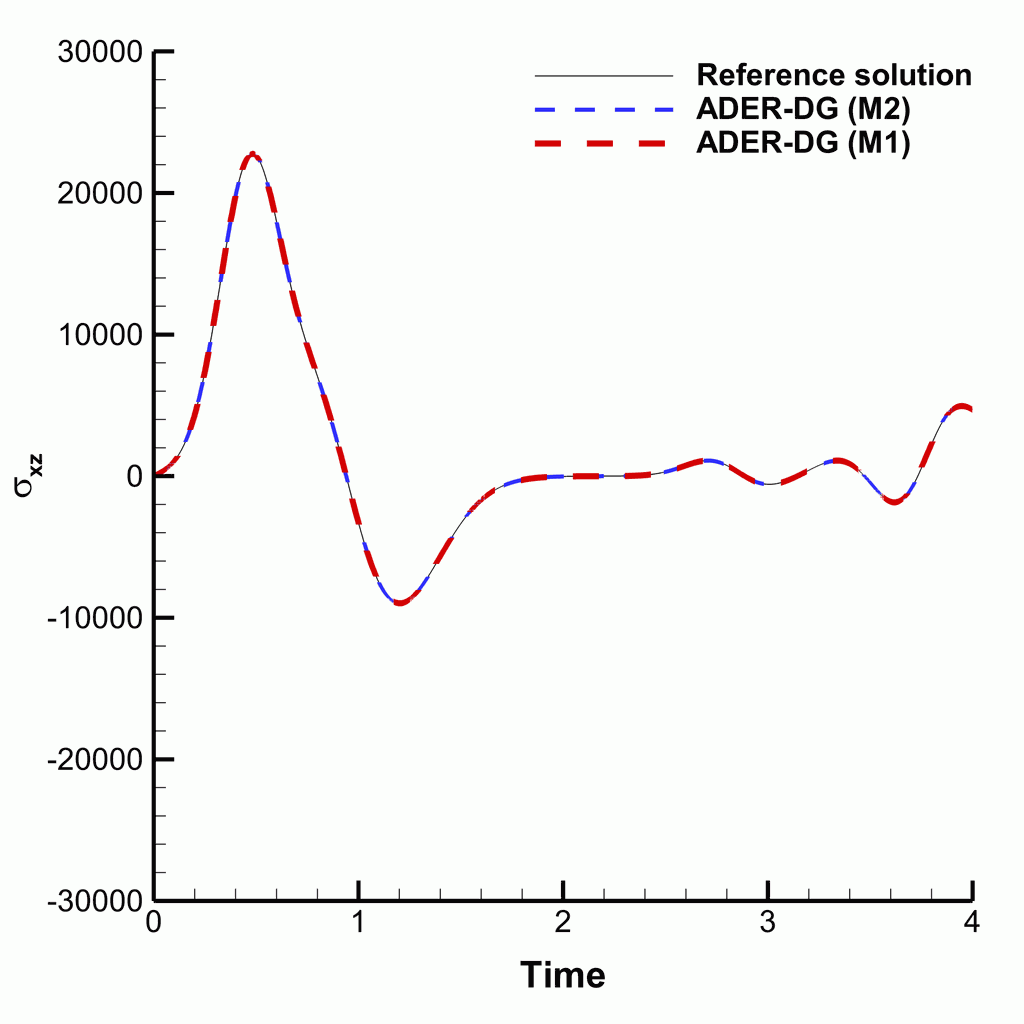}\\	
	\includegraphics[width=0.32\linewidth]{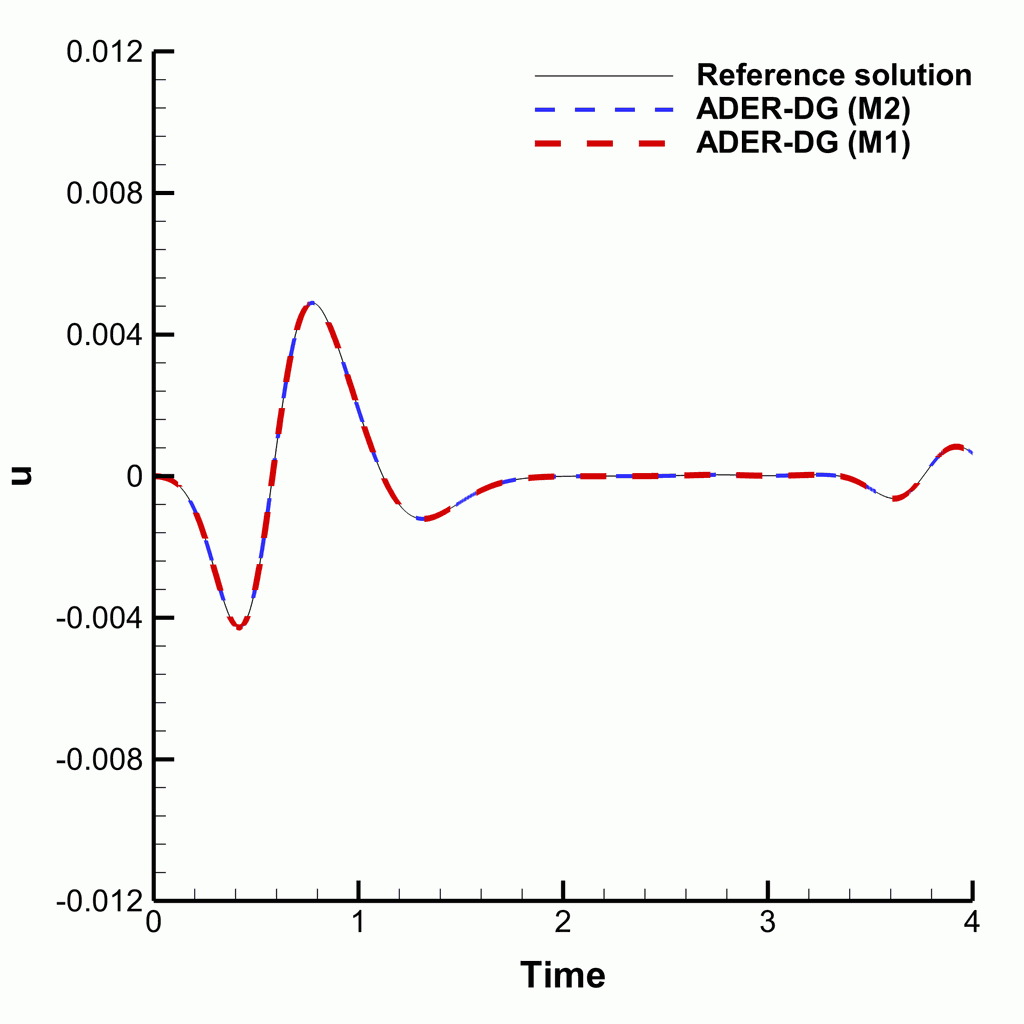}\hfill
	\includegraphics[width=0.32\linewidth]{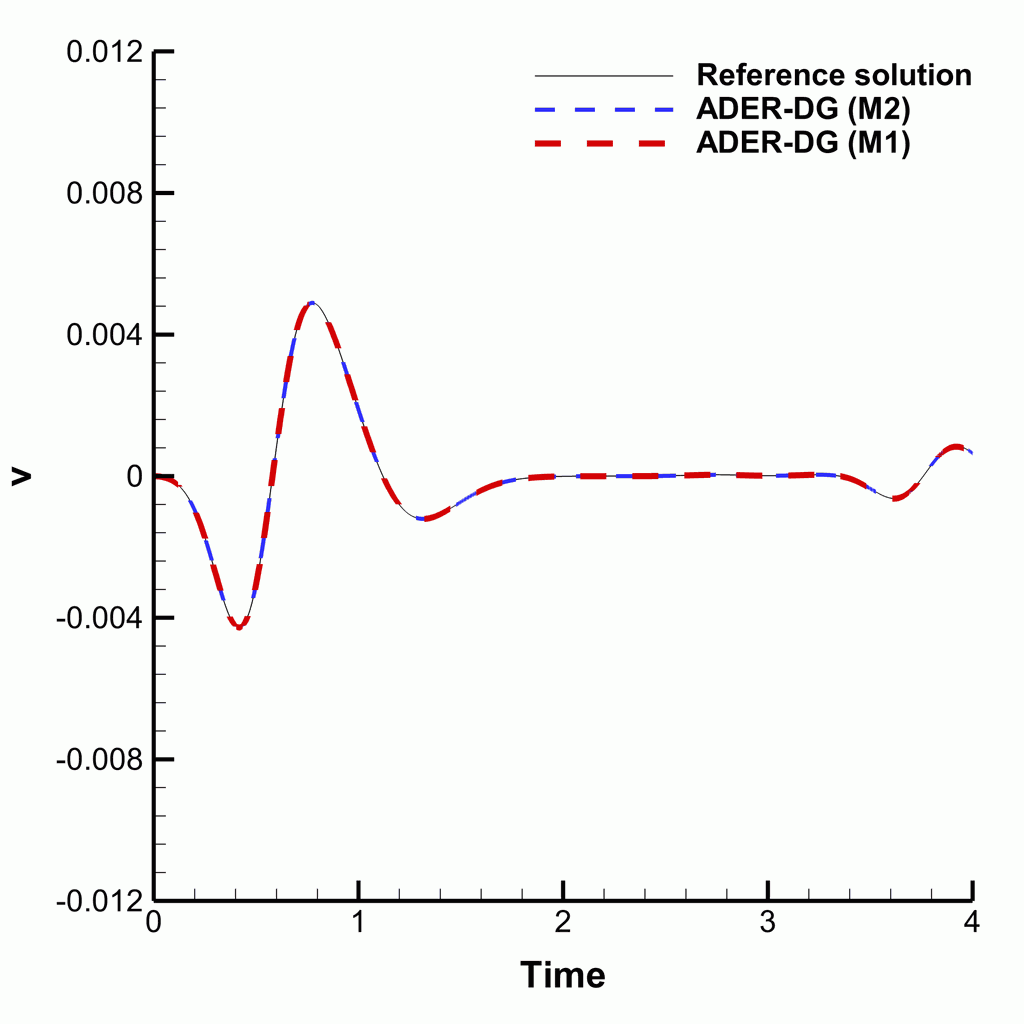}\hfill
	\includegraphics[width=0.32\linewidth]{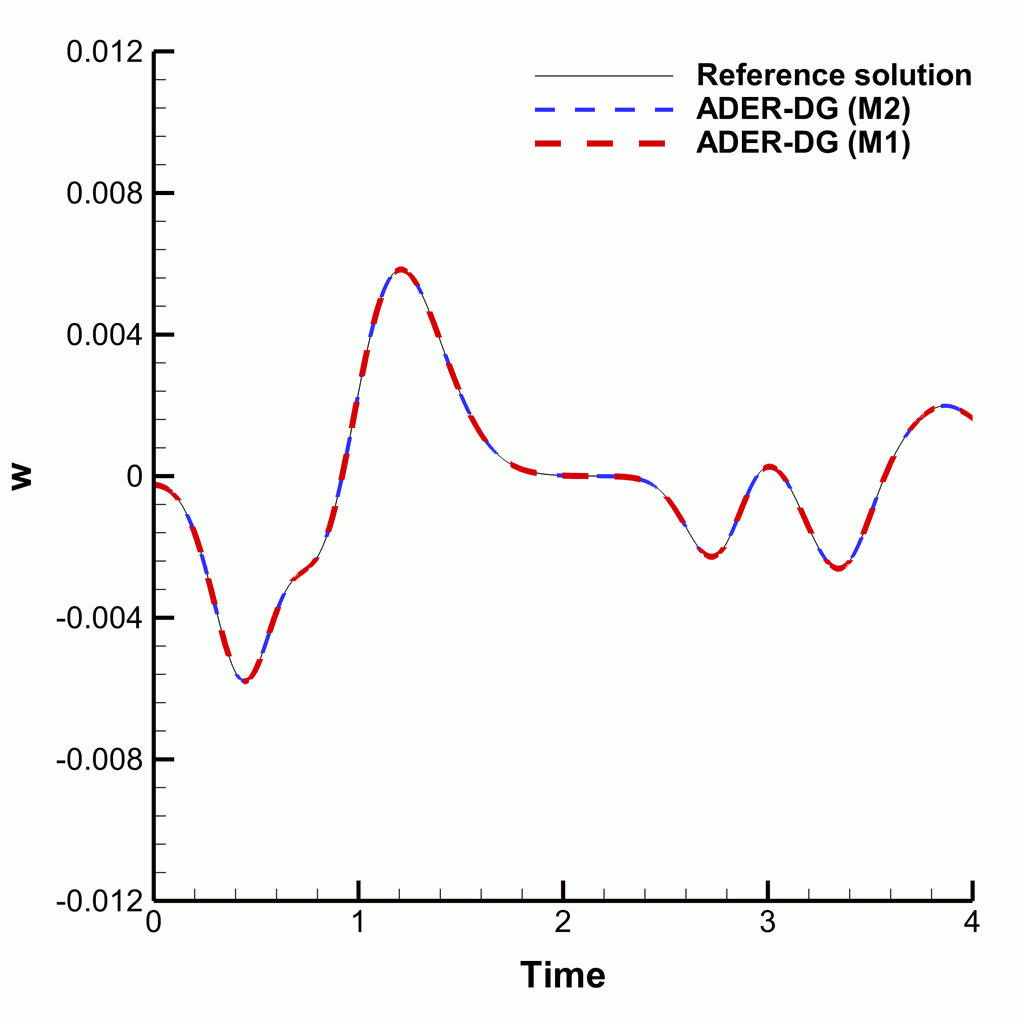}
	
	\caption{3D wave propagation. Seismogram at receiver $\mathbf{x}_{1}=\left(1000,1000,1000\right)$ obtained using ADER-DG $\mathcal{O}4$ on two different meshes, M1 made of $30\times 30 \times 30$ (red long dashed line) and M2 using $60\times 60 \times 60$ elements (navy short dashed line), and reference solution computed with a $P_{3}P_{3}$ ADER-DG scheme on an unstructured grid of $195301$ tetrahedra (black line).}
	\label{fig:3DExp_timeevolution_xyz1000}
\end{figure}

\begin{figure}
	\centering
	\includegraphics[width=0.32\linewidth]{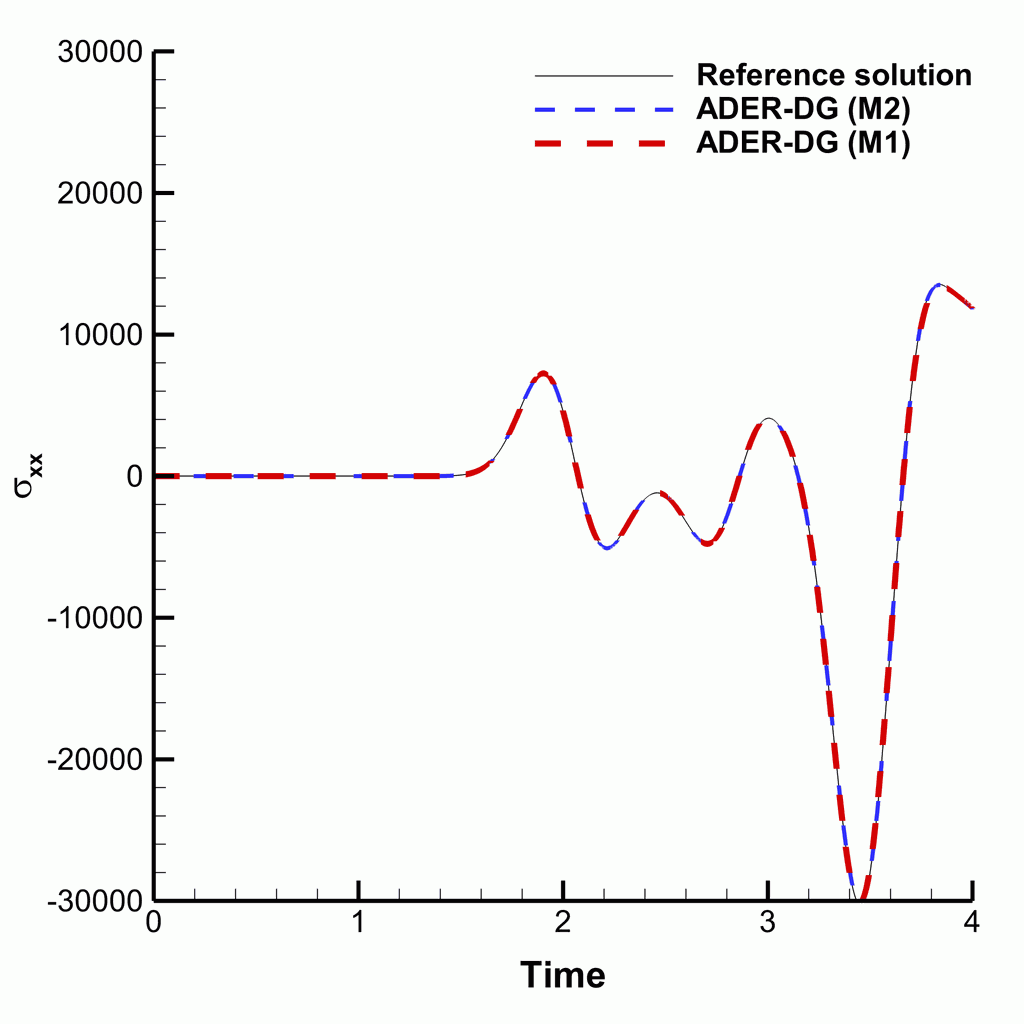}\hfill
	\includegraphics[width=0.32\linewidth]{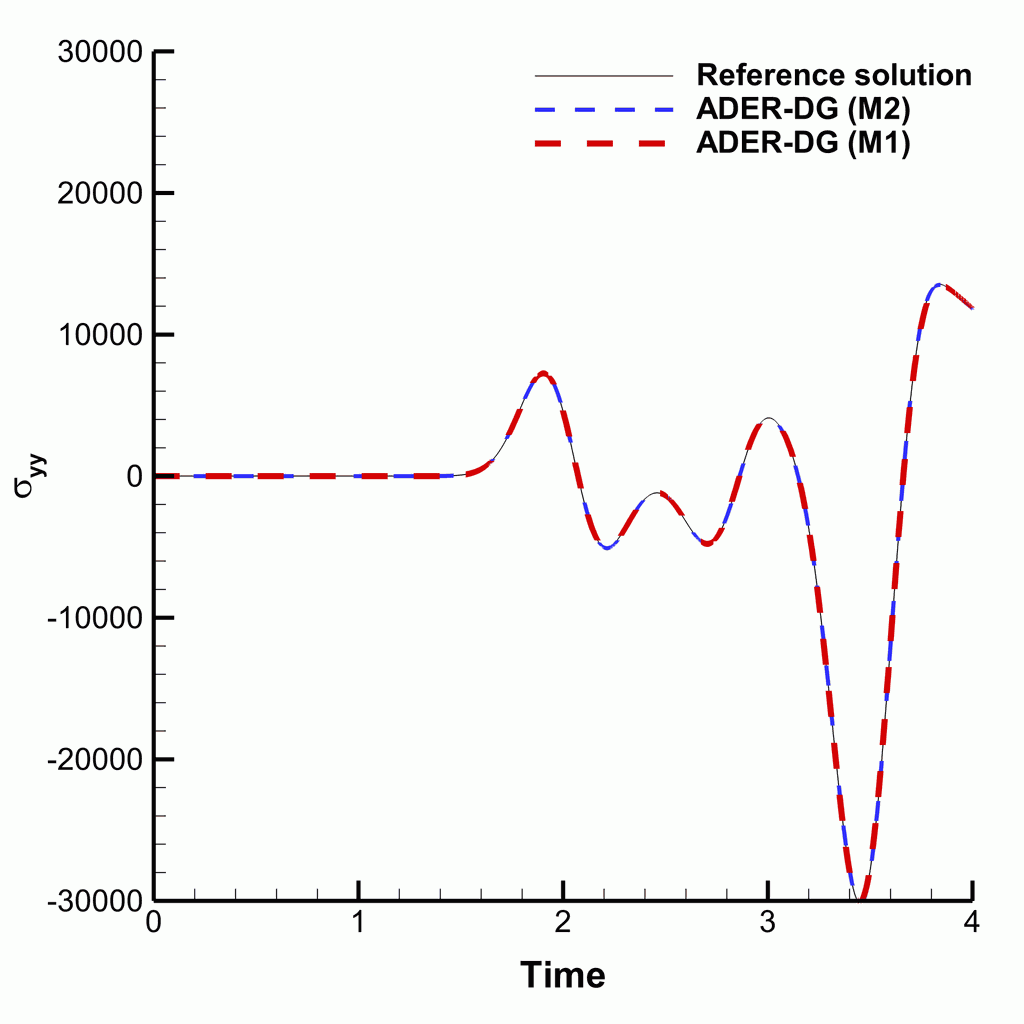}\hfill
	\includegraphics[width=0.32\linewidth]{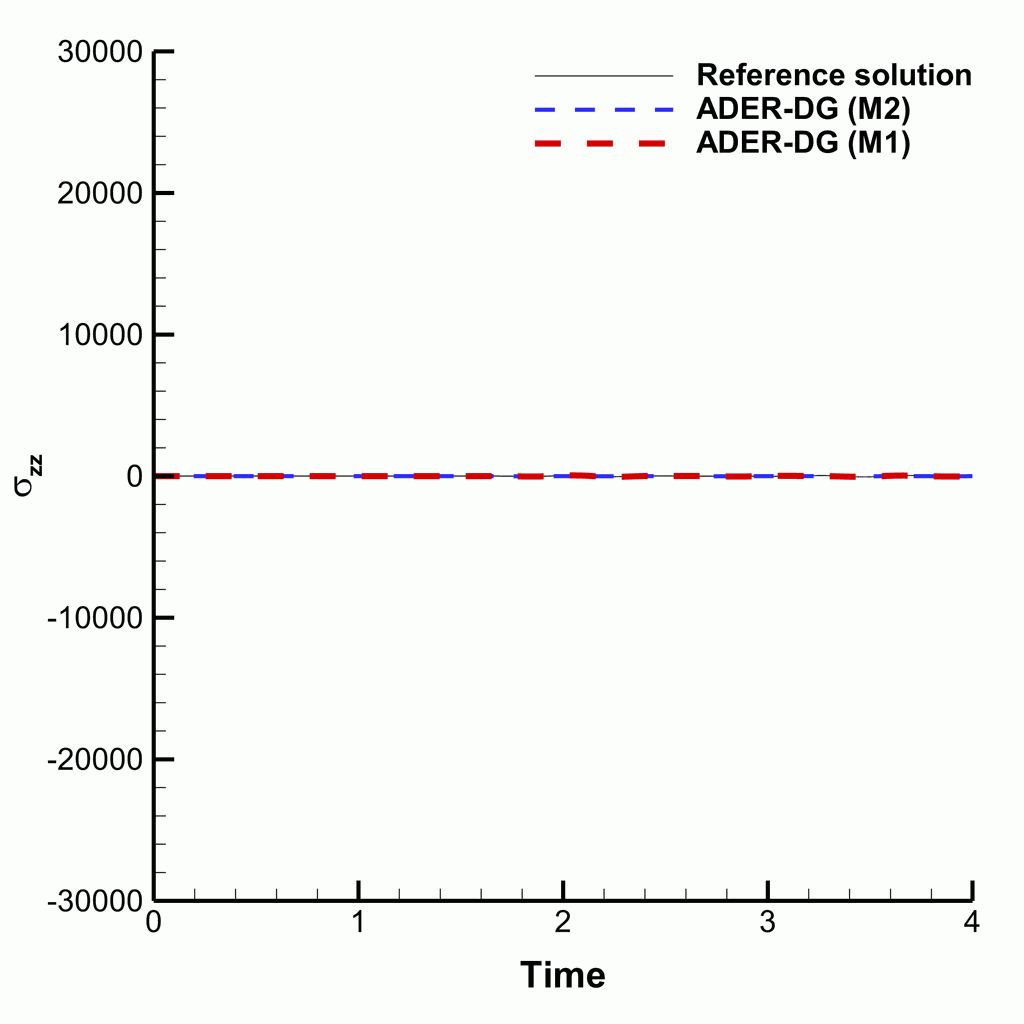}\\	
	\includegraphics[width=0.32\linewidth]{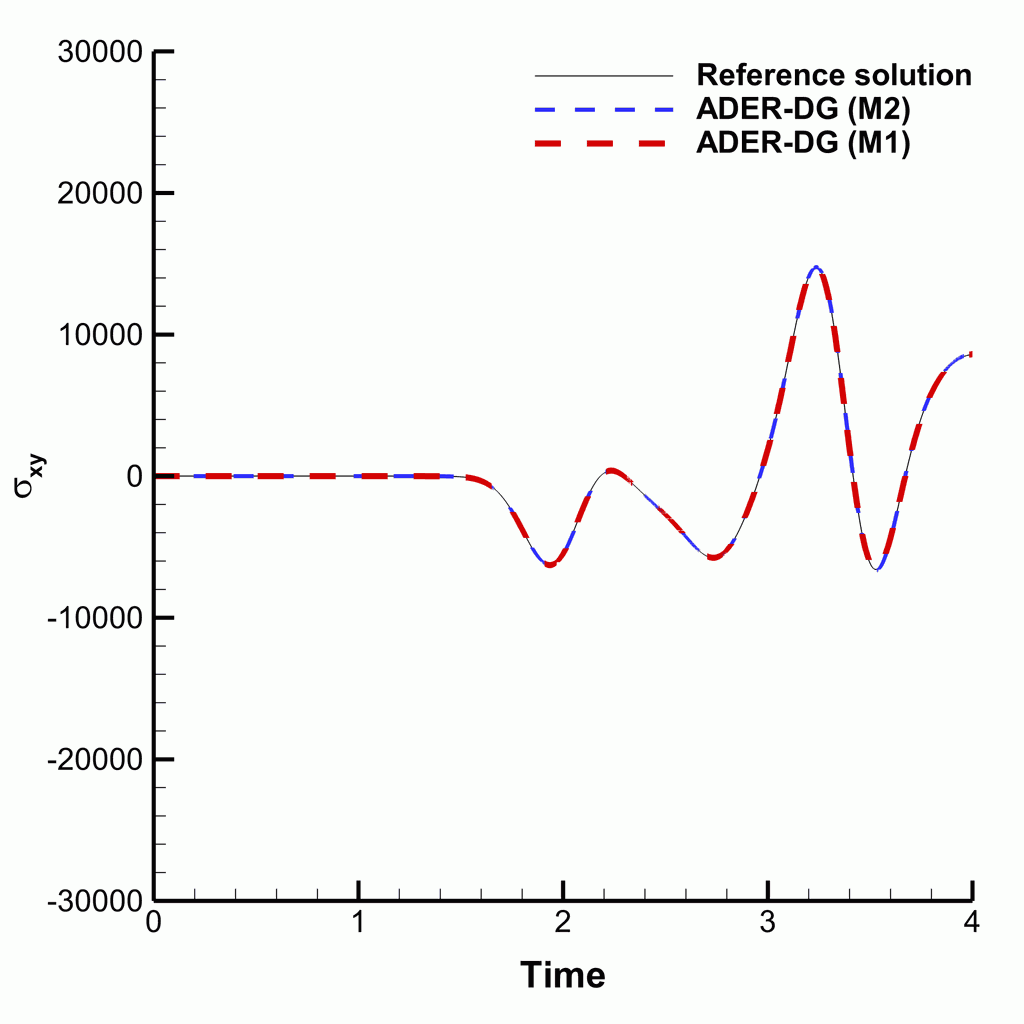}\hfill
	\includegraphics[width=0.32\linewidth]{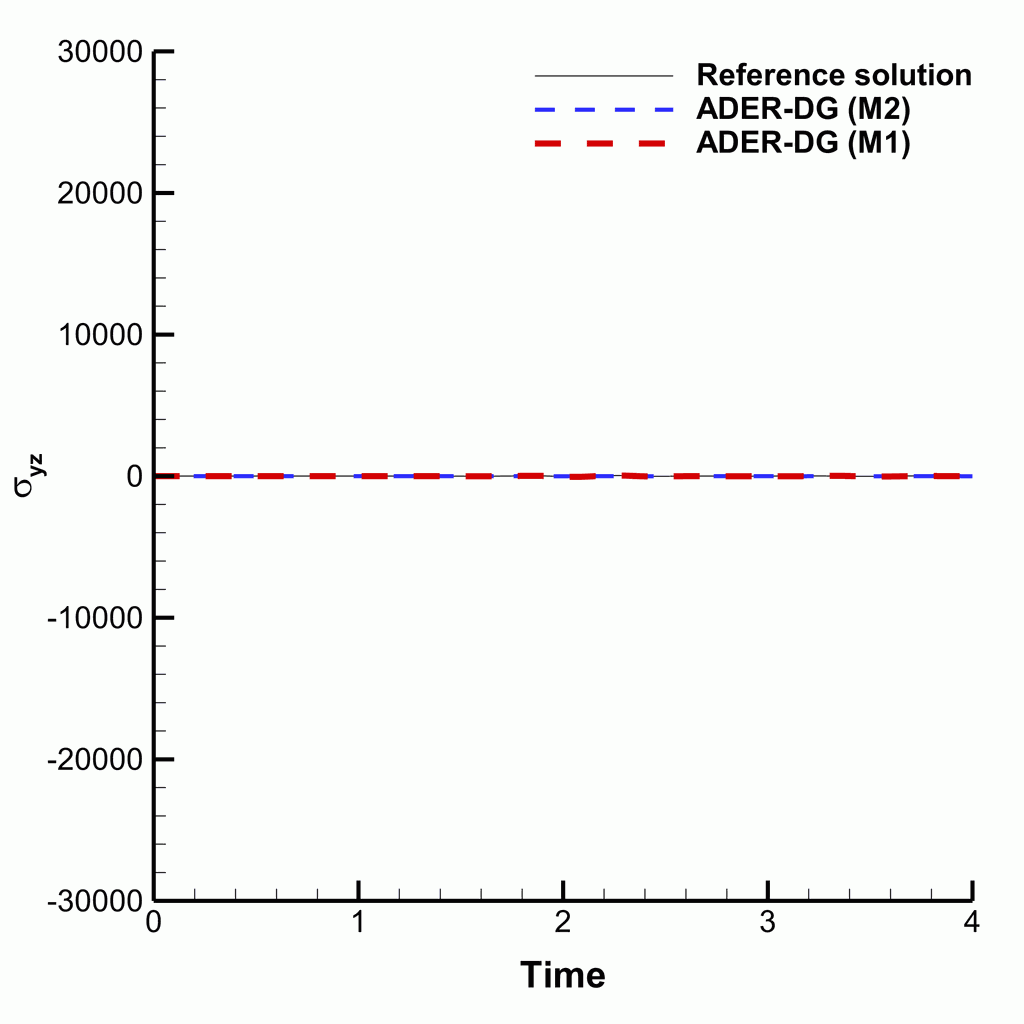}\hfill
	\includegraphics[width=0.32\linewidth]{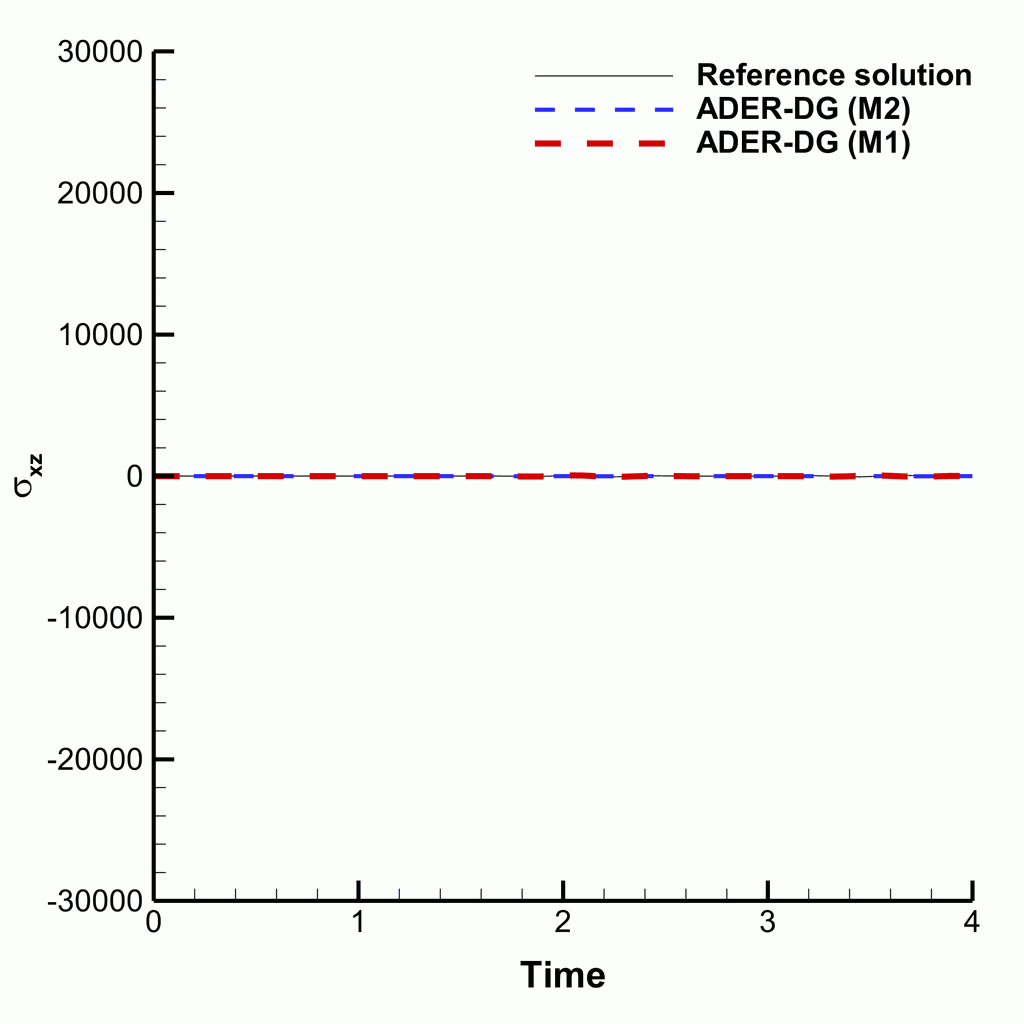}\\	
	\includegraphics[width=0.32\linewidth]{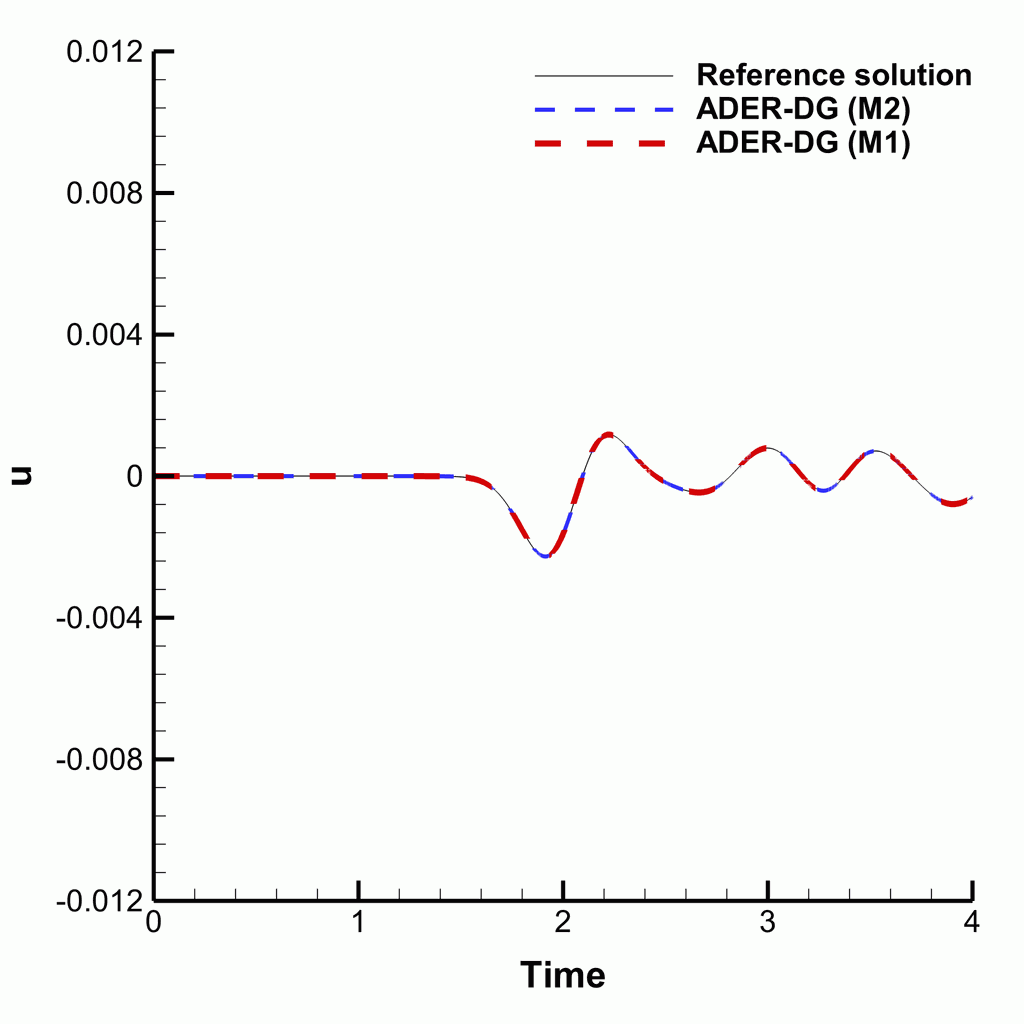}\hfill
	\includegraphics[width=0.32\linewidth]{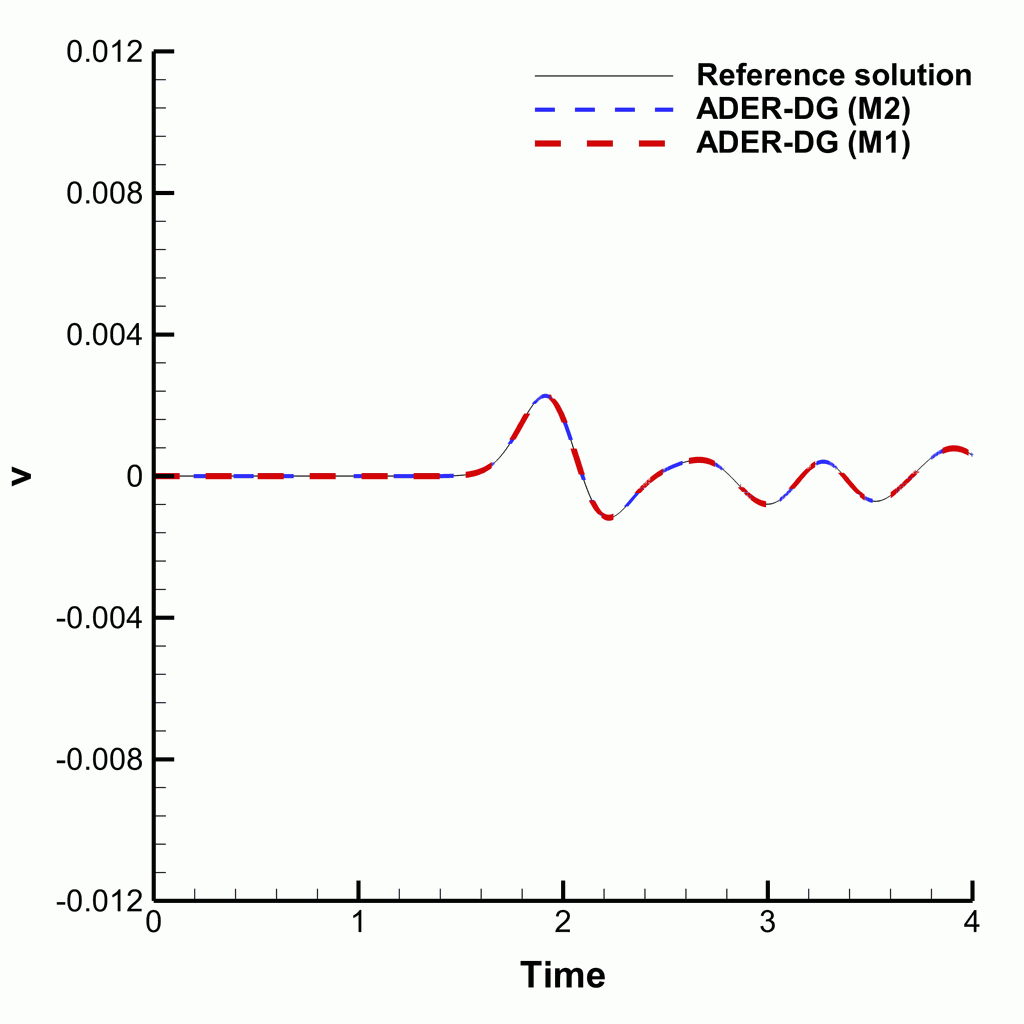}\hfill
	\includegraphics[width=0.32\linewidth]{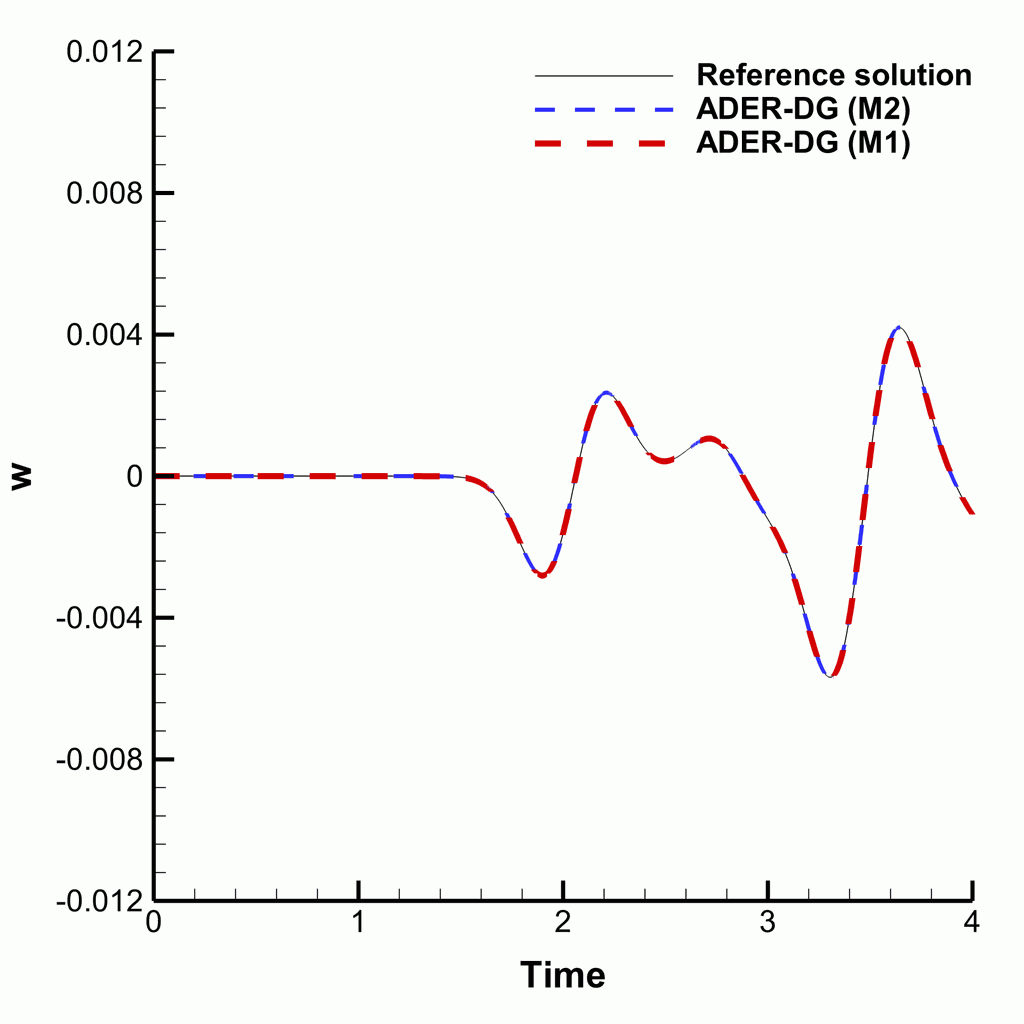}
	
	\caption{3D wave propagation. Seismogram at receiver $\mathbf{x}_{r_{2}}=\left(3000,-3000,5000\right)$ obtained using ADER-DG $\mathcal{O}4$ on two different meshes, M1 made of $30\times 30 \times 30$ (red long dashed line) and M2 using $60\times 60 \times 60$ elements (navy short dashed line), and reference solution computed with a $P_{3}P_{3}$ ADER-DG scheme on an unstructured grid of $195301$ tetrahedra (black line).}
	\label{fig:3DExp_timeevolution_x3000y_3000z5000}
\end{figure}

\begin{figure}
	\centering
	\includegraphics[width=0.32\linewidth]{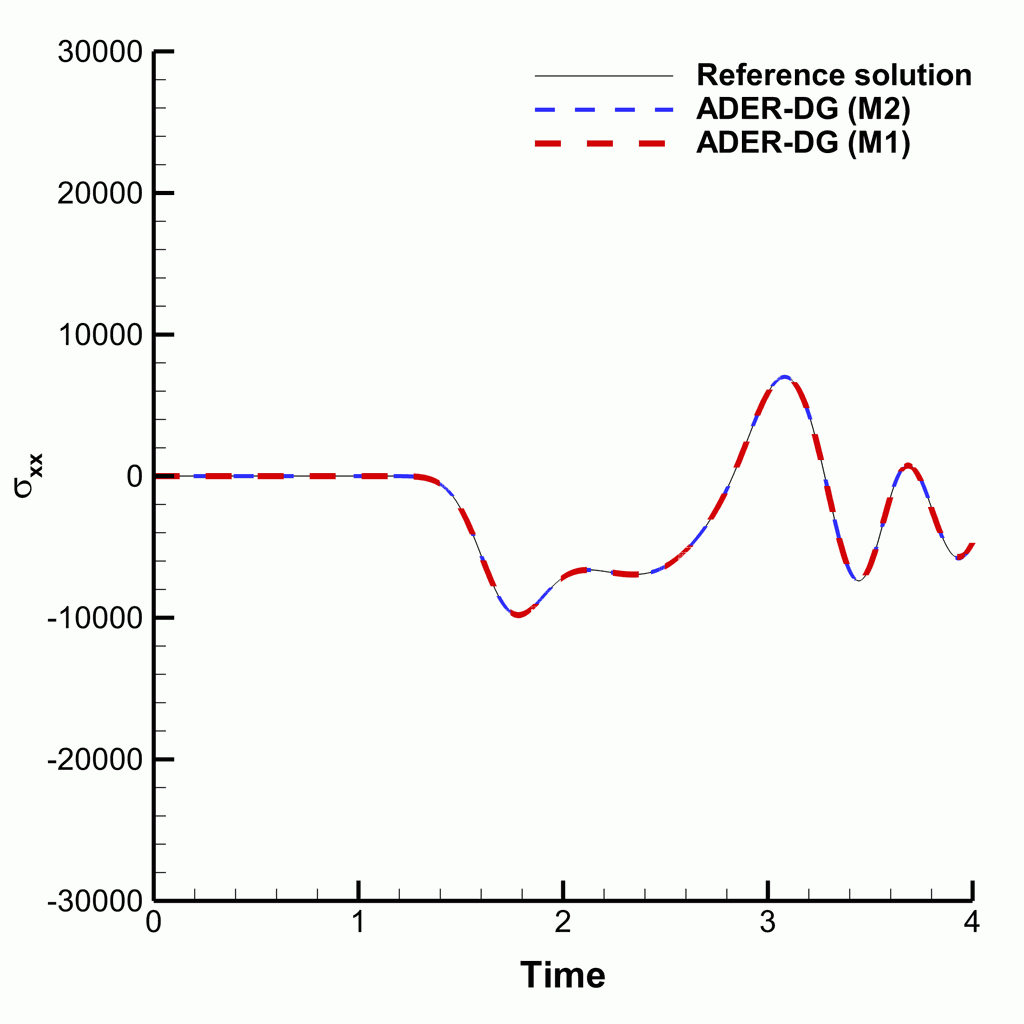}\hfill
	\includegraphics[width=0.32\linewidth]{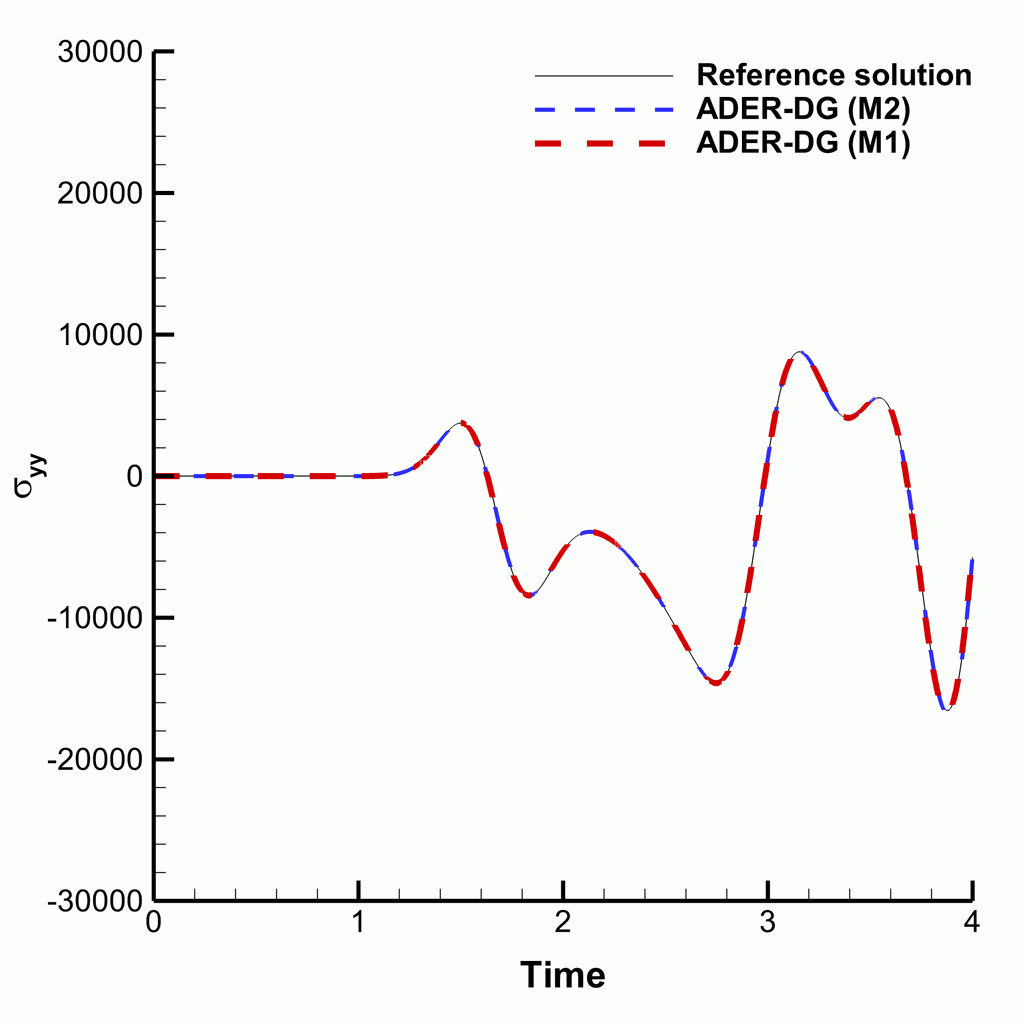}\hfill
	\includegraphics[width=0.32\linewidth]{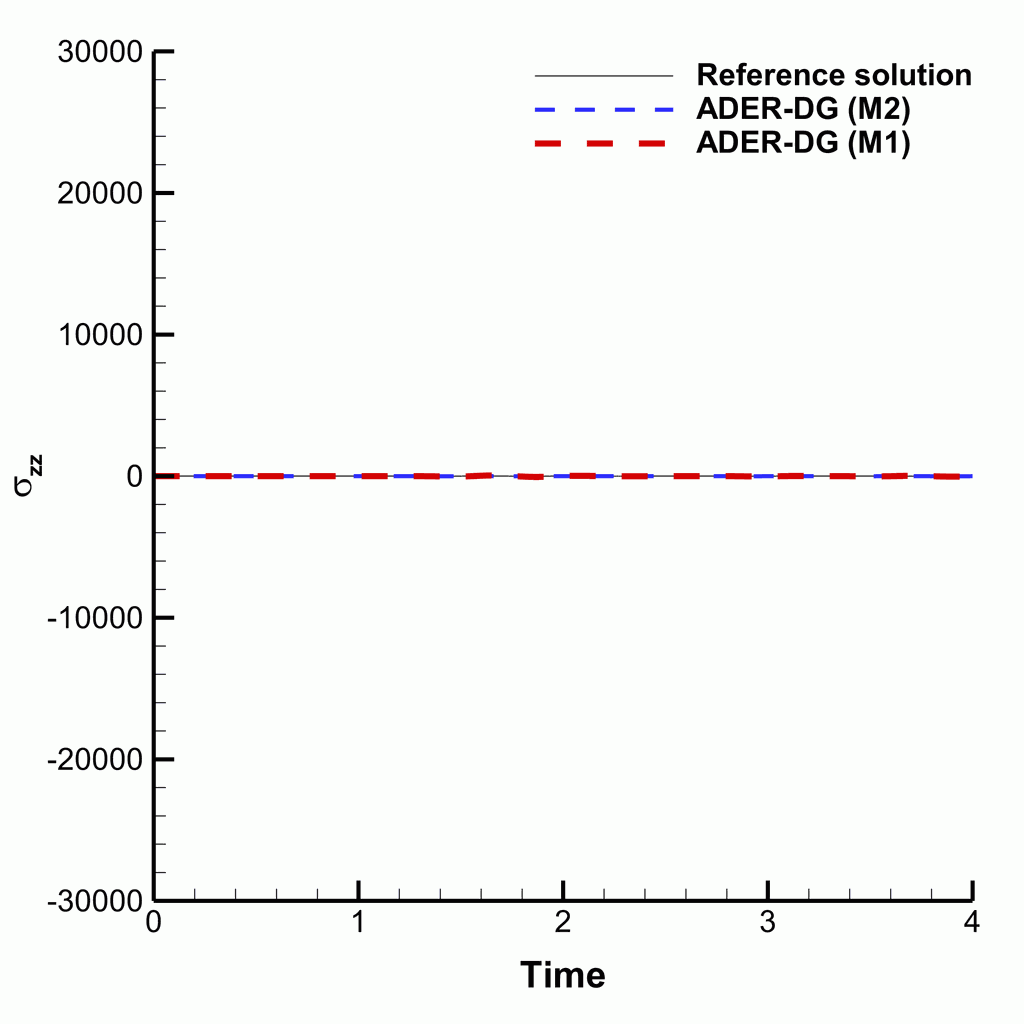}\\	
	\includegraphics[width=0.32\linewidth]{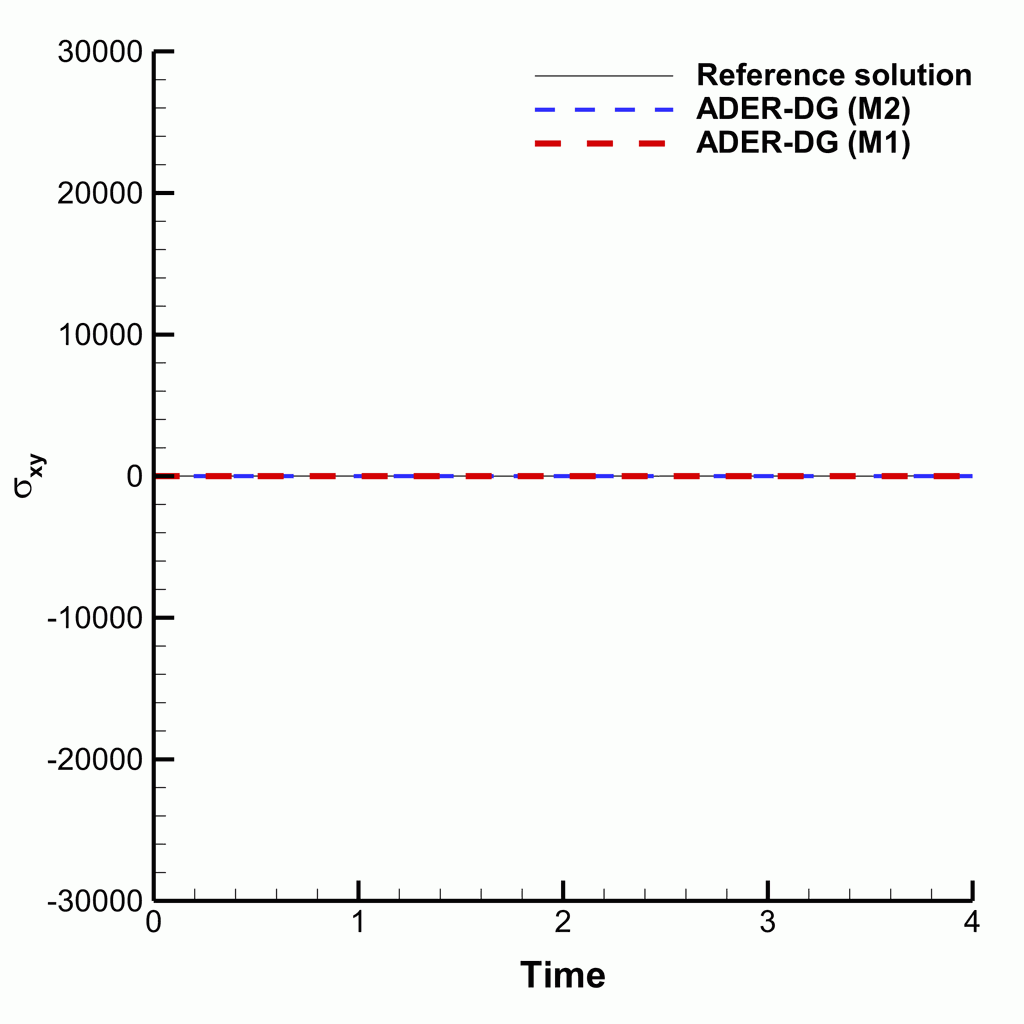}\hfill
	\includegraphics[width=0.32\linewidth]{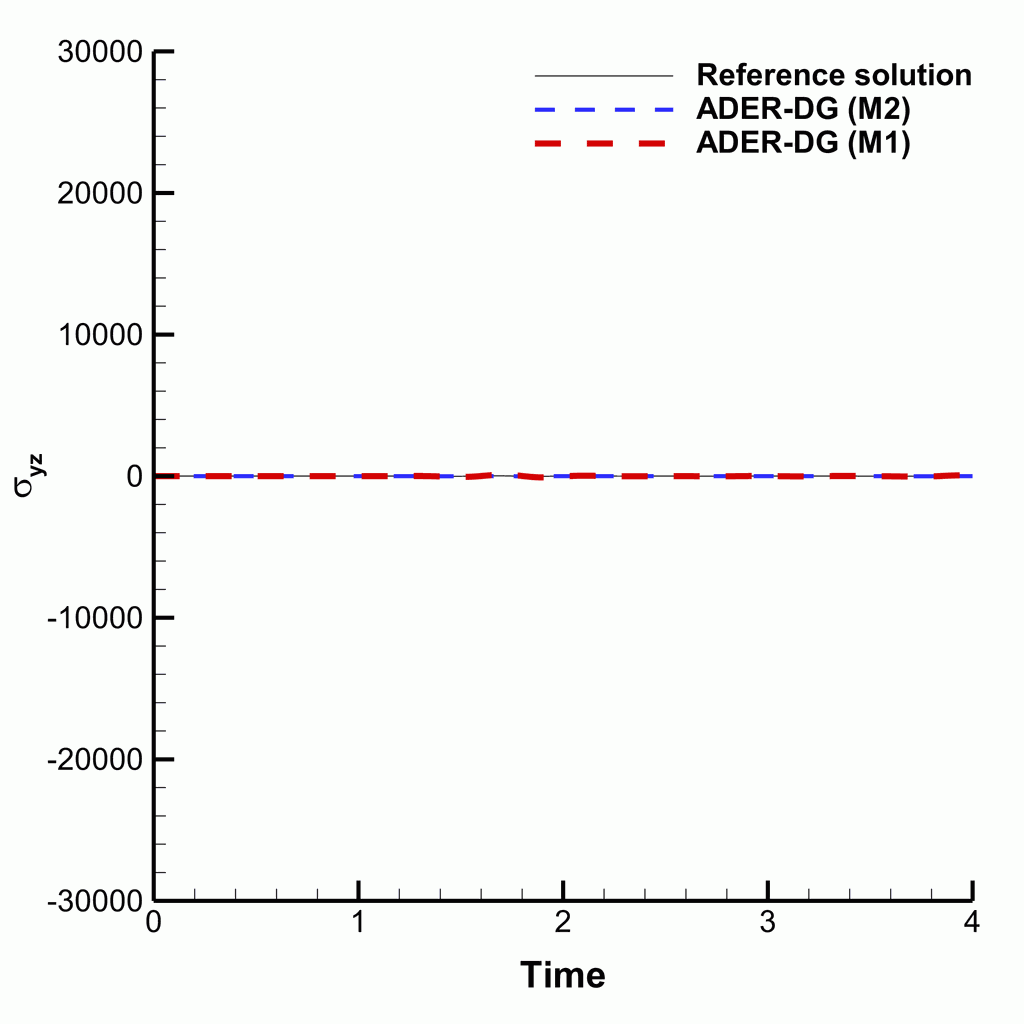}\hfill
	\includegraphics[width=0.32\linewidth]{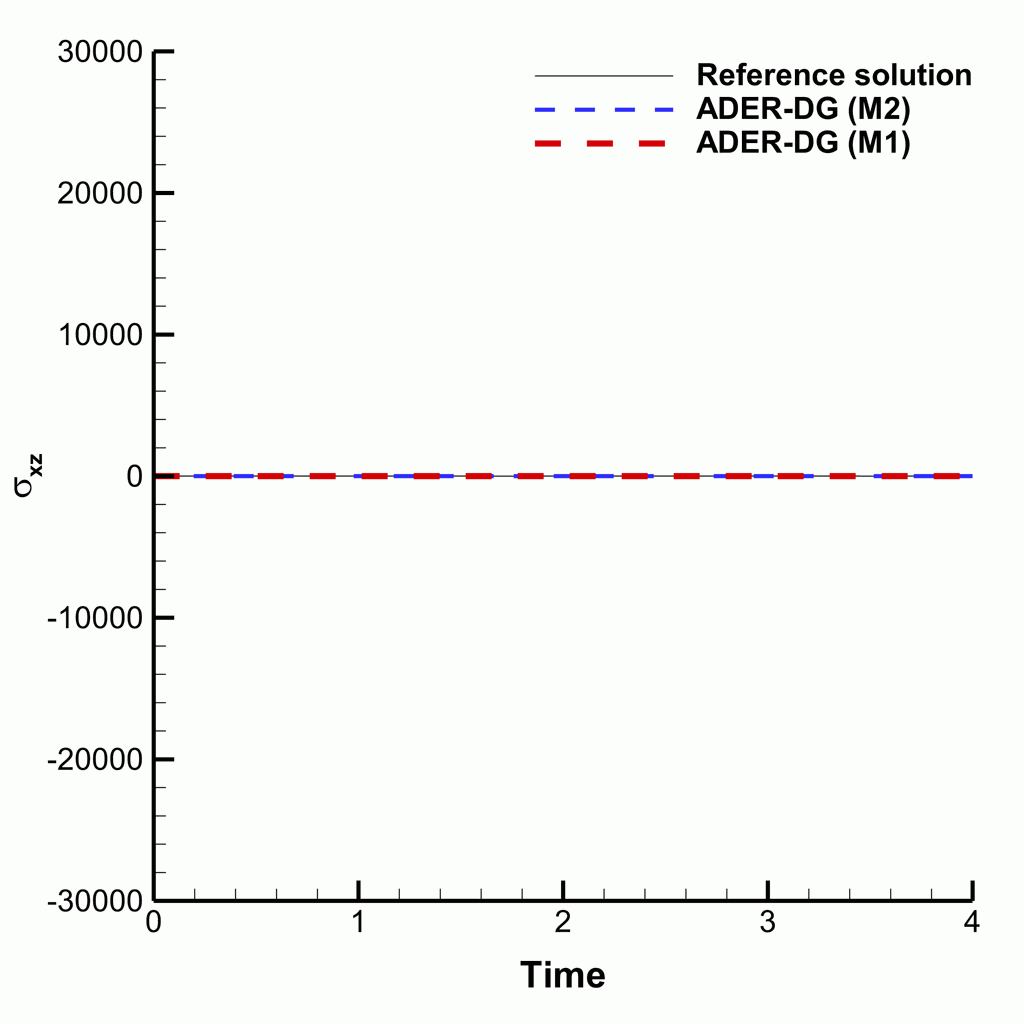}\\	
	\includegraphics[width=0.32\linewidth]{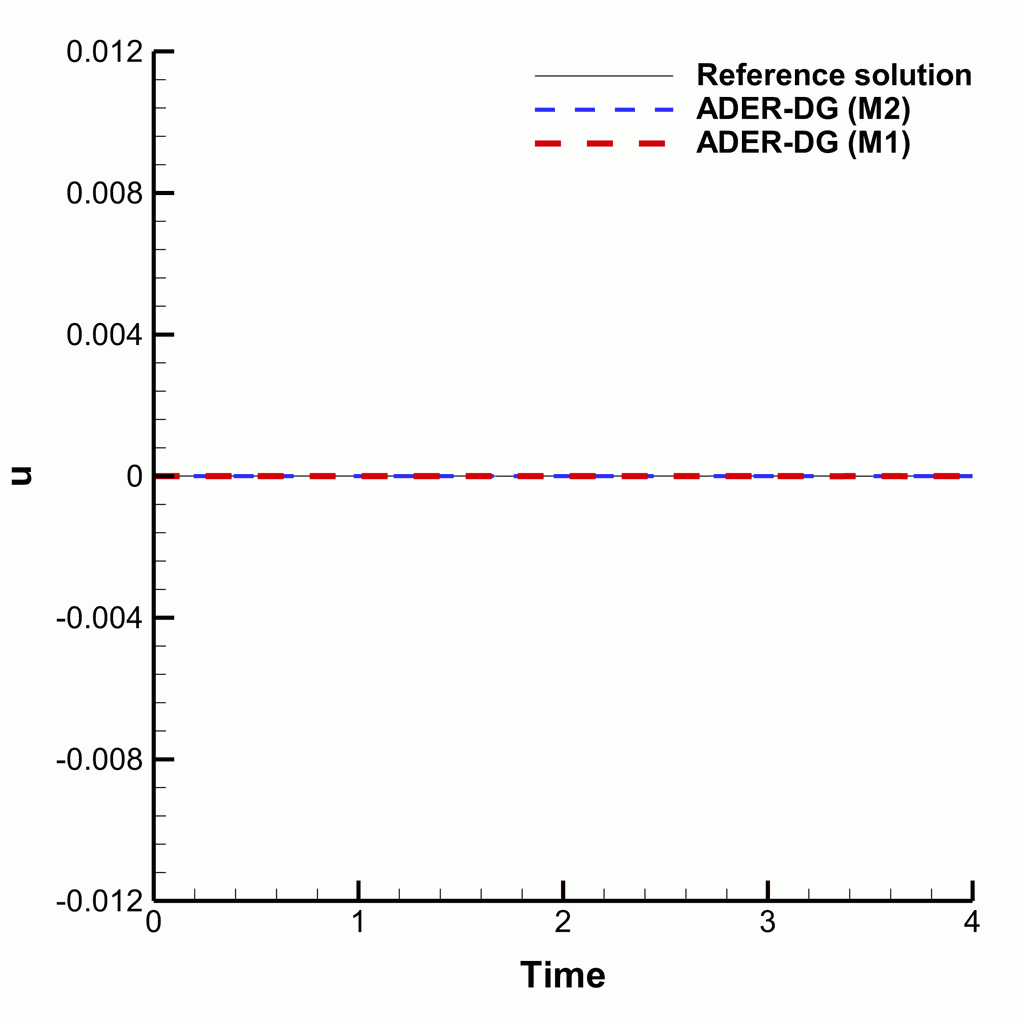}\hfill
	\includegraphics[width=0.32\linewidth]{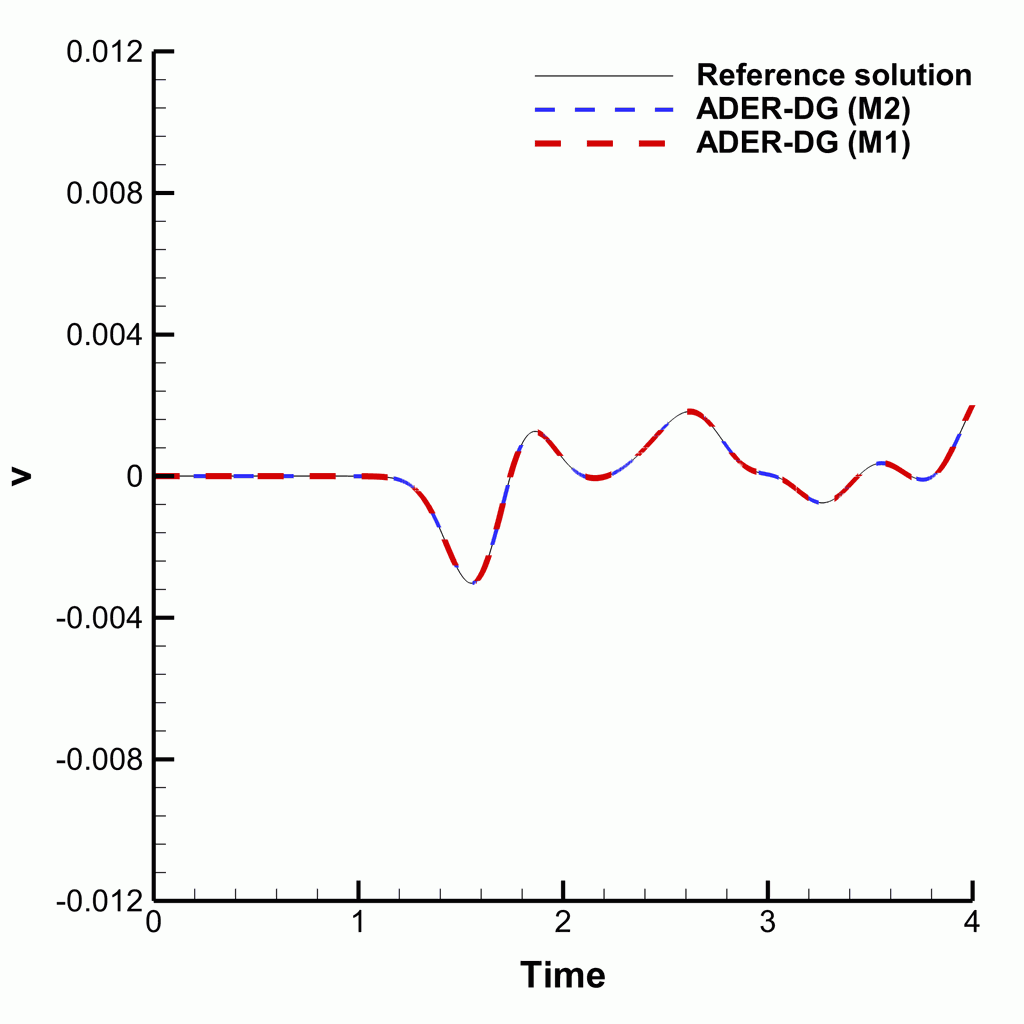}\hfill
	\includegraphics[width=0.32\linewidth]{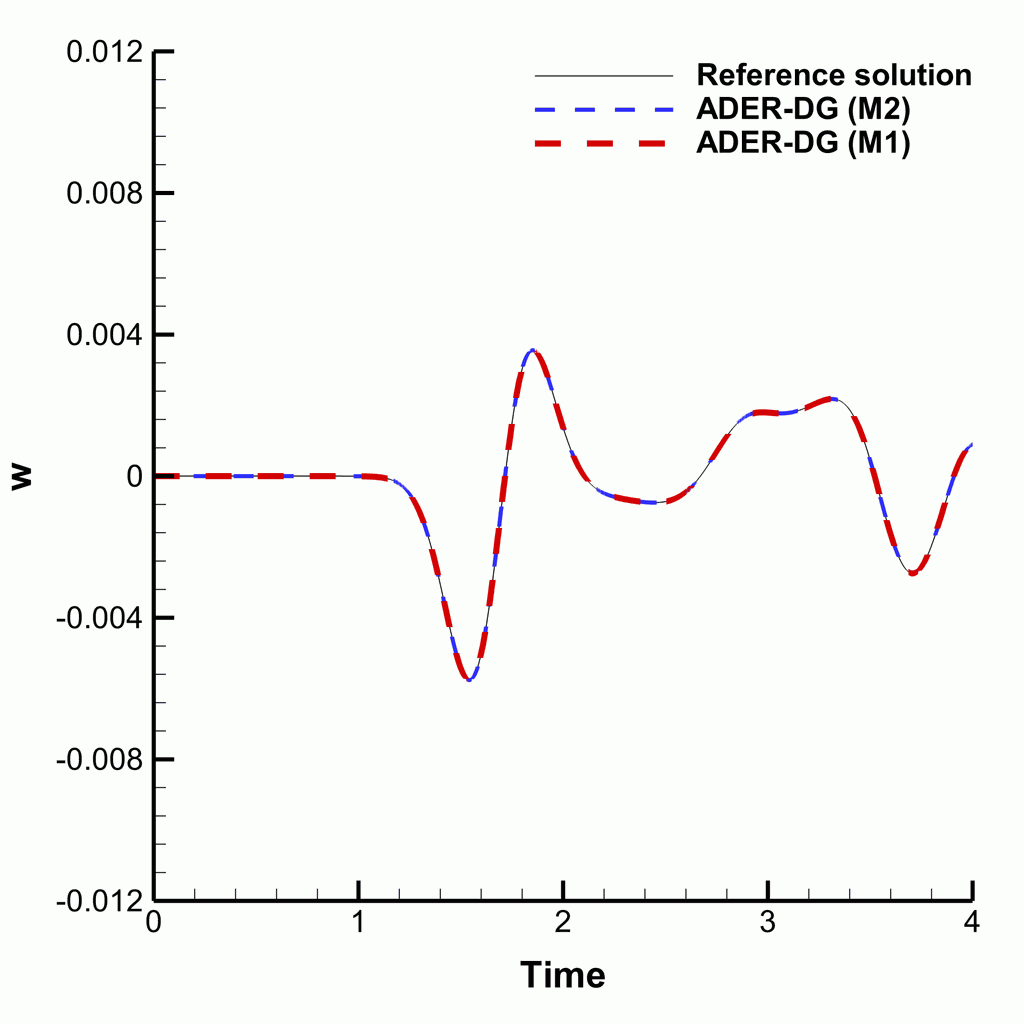}
	
	\caption{3D wave propagation. Seismogram at receiver $\mathbf{x}_{r_{3}}=\left(0,2000,5000\right)$ obtained using ADER-DG $\mathcal{O}4$ on two different meshes, M1 made of $30\times 30 \times 30$ (red long dashed line) and M2 using $60\times 60 \times 60$ elements (navy short dashed line), and reference solution computed with a $P_{3}P_{3}$ ADER-DG scheme on an unstructured grid of $195301$ tetrahedra (black line).}
	\label{fig:3DExp_timeevolution_x0y2000z5000}
\end{figure}

The third test to be analysed corresponds to a 3D wave propagation problem. We consider the computational domain $\Omega=\left[-5000,5000\right]\times\left[-5000,5000\right]\times\left[-5000,5000\right]$ and a homogeneous material with propagation velocities $c_{p}=3200$, $c_{s}=1847.5$ and density $\rho=2200$. Following \cite{TD18LE}, the initial condition for the vertical velocity is given by the Gaussian profile 
\begin{equation}
v_{3} = -0.1 \exp\left( \frac{-r^{2}}{2R^{2}}\right),
\end{equation}
where $R=500$ denotes the initial impulse size and $r$ is the distance with respect to the centre of the impulse, $\mathbf{x}_{0}=\left(0,0,0\right)$. The remaining variables are set to zero. 
The simulation has been run on a Cartesian grid made of $30\times 30 \times 30$ elements, using the fourth order scheme. In Figure \ref{fig:3DExp_w_t2}, the vertical velocity contours obtained at $t=2$ are plotted. To study the wave propagation three receivers have been placed at $\mathbf{x}_{r_{1}}=\left(1000,1000,1000\right)$, $\mathbf{x}_{r_{2}}=\left(3000,-3000,5000\right)$, and $\mathbf{x}_{r_{3}}=\left(0,2000,5000\right)$. The time evolution for the main variables involved is depicted in Figures \ref{fig:3DExp_timeevolution_xyz1000}, \ref{fig:3DExp_timeevolution_x3000y_3000z5000}, and \ref{fig:3DExp_timeevolution_x0y2000z5000}, respectively. 
As expected, due to the location of receiver $r_{1}$, the values obtained in  $\mathbf{x}_{r_{1}}$ for $\sigma_{xx}$, $\sigma_{xz}$, and $u$ coincide with $\sigma_{yy}$, $\sigma_{yz}$, and $v$, respectively.
Similarly, at $\mathbf{x}_{r_{2}}$, $\sigma_{xx}$ matches $\sigma_{yy}$ whereas $u$ and $v$ take opposite values. Finally, since ${r_{3}}$ is located on the $x=0$ plane, we should get zero horizontal velocity.
The solution obtained with a finer mesh of $60\times 60 \times 60$ elements is also included to demonstrate that mesh convergence is attained. Moreover, the obtained results are compared against a reference solution obtained with the unstructured ADER-DG code \texttt{SeisSol} using a $P_{3}P_{3}$ scheme on a mesh made of $195301$ tetrahedra. Again, we can note an excellent agreement between our solution and the numerical reference solution obtained with a community code.   

\subsection{HSGN numerical results}
In the HSGN framework, we employ a solitary wave over a flat bottom test to assess the methodology. Next, a step-shaped bathymetry is considered and the obtained propagation of a soliton wave is compared against the experimental data and the numerical results obtained for the original SGN system. Further analysis on the employed ADER-DG method applied to non-hydrostatic free surface models can be found in \cite{BBBD20}.

\subsubsection{Solitary wave over a flat bottom}\label{sec:SoWFB}
\begin{figure}
	\centering
	\includegraphics[width=0.45\linewidth]{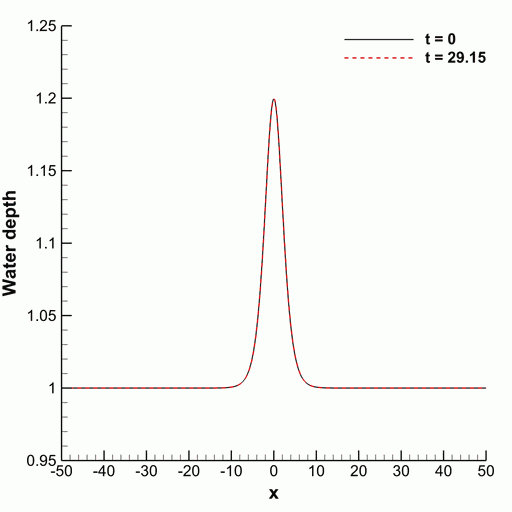}\hspace{0.05\linewidth}
	\includegraphics[width=0.45\linewidth]{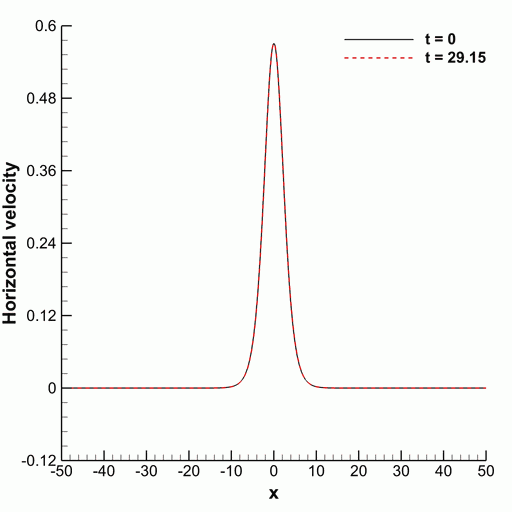}
	
	\includegraphics[width=0.45\linewidth]{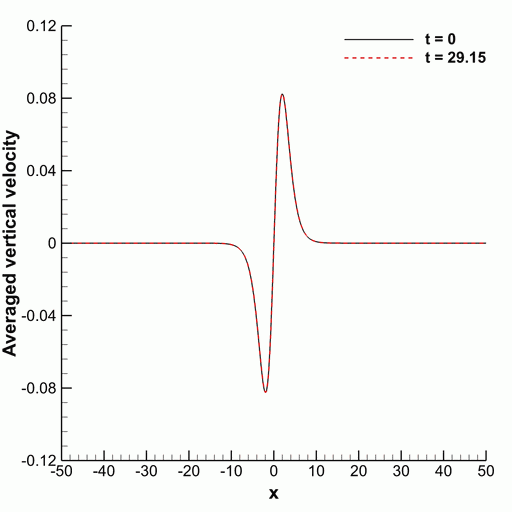}\hspace{0.05\linewidth}
	\includegraphics[width=0.45\linewidth]{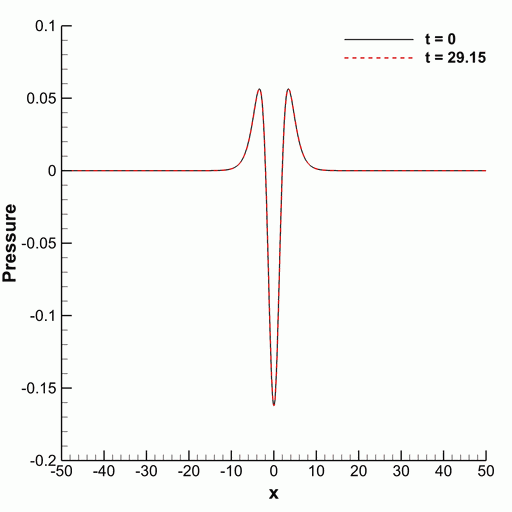}
	
	\caption{Solitary wave over a flat bottom. 1D profile obtained after a complete revolution of the soliton (red dashed line) and initial solution (black line). From left top to right bottom: water depth, $h$, horizontal velocity, $u$, averaged vertical velocity, $w$, pressure, $p$.} 
	\label{fig:SW_p5_m100_t2915}
\end{figure}

To assess the accuracy of the method, we study a solitary wave propagating over a flat bottom, \cite{EDC19,BBBD20}. 
The computational domain is taken to be $\Omega=[-50,50]\times [-1,1]$. As initial condition, we define a soliton of amplitude $A=0.2$ centred at the origin $(x,y)=(0,0)$, an artificial sound velocity $c=20$ and a still water depth $H=1$. Periodic boundary conditions are considered at all boundaries. Let us remark that, for the original non-hyperbolic formulation of the SGN model, an analytical solution is available, see e.g. \cite{JSMarie}. However, this solution does not exactly verify the hyperbolic formulation, thus it should not be employed in a convergence study. Instead, we consider a 1D self-similar solution of the hyperbolic system \eqref{eq:hgn_depth}-\eqref{eq:hgn_pressure} of the form
\begin{equation}
\mathbf{Q}\left( x,t\right) = \mathbf{Q} \left( \zeta \right), \quad \zeta = x - V t,
\end{equation}
$V$ the velocity of the solitary wave, obtained by solving the corresponding nonlinear ODE,
\begin{equation}
\mathbf{Q}^{\prime} = \left( \mathbf{A}\left(\mathbf{Q}\right)  -V \mathbf{I}\right)^{-1} \mathbf{S} \left( \mathbf{Q}\right), \quad \mathbf{A}\left(\mathbf{Q}\right) = \frac{\partial \mathbf{F}}{\partial \mathbf{Q}} + \mathbf{B}\left(\mathbf{Q}\right),
\end{equation}
with initial condition $\mathbf{Q}\left(\zeta_{0}\right) =\left( H_{0},0,0,0,\epsilon\right)$, $\epsilon=10^{-8}$, $H_{0}=H$. The former ODE is solved using a tenth order discontinuous Galerkin scheme, see \cite{ADERNSE}. The solution is used both for the initialization of the soliton and to compute the $L^2$ errors at $t_{\mathrm{end}}=2$. The errors and convergence rates obtained for the density, horizontal velocity, and pressure, using polynomial degrees $N\in\left\lbrace 3,4,5,6,7\right\rbrace$ are depicted in Table~\ref{tab:SW_accuracy}. Overall, the sought order of convergence is achieved but for some particular cases, in which a suboptimal order can be observed on some of the variables. Besides, Figure \ref{fig:SW_p5_m100_t2915} shows the 1D profiles of water depth, horizontal velocity, averaged vertical velocity and pressure obtained after a one complete revolution of the soliton, $t=29.15$. We observe that the results, computed on a mesh made of $100\times 2$ elements using $N=5$, perfectly match the initial condition.

\begin{table}
	\centering
	\begin{tabular}{|c|c||c|c|c||c|c|c||c|}
		\hline
		$N$ & $N_{x}$ & $\epsilon_{L^{2}}\left(h\right) $ & $\epsilon_{L^{2}}\left(u\right) $ & $\epsilon_{L^{2}}\left(p\right) $& $\mathcal{O}_{L^{2}}\left(h\right) $ & $\mathcal{O}_{L^{2}}\left(u\right) $ & $\mathcal{O}_{L^{2}}\left(p\right) $ & Teor.\\
		\hline\hline
		\multirow{5}{*}{$3$} &$80$  & $1.05E-03$ & $8.64E-04$ & $8.33E-03$ & - & - & - & \multirow{5}{*}{$4$}  \\	\hhline{|~|-||-|-|-||-|-|-||~|}
		&$100$ & $3.83E-04$ & $2.84E-04$ & $3.47E-03$ & $4.52$ & $4.99$ & $3.93$ &  $$  \\	\hhline{|~|-||-|-|-||-|-|-||~|}
		&$120$ & $1.60E-04$ & $9.95E-05$ & $1.57E-03$ & $4.81$ & $5.74$ & $4.33$ & $$  \\	\hhline{|~|-||-|-|-||-|-|-||~|}
		&$140$ & $7.75E-05$ & $4.02E-05$ & $7.60E-04$ & $4.69$ & $5.88$ & $4.72$ & $$  \\	\hhline{|~|-||-|-|-||-|-|-||~|}
		&$160$ & $4.32E-05$ & $1.84E-05$ & $4.40E-04$ & $4.37$ & $5.86$ & $4.09$ & $$  \\	\hline\hline
		\multirow{5}{*}{$4$}& $60$ & $3.36E-04$ & $2.55E-04$ & $3.39E-03$ & - & - & - & \multirow{5}{*}{$5$} \\	\hhline{|~|-||-|-|-||-|-|-||~|}
		& $80 $ & $6.49E-05$ & $4.08E-05$ & $7.88E-04$ & $5.72$ & $6.38$ & $5.07$ &   \\	\hhline{|~|-||-|-|-||-|-|-||~|}
		& $100$ & $1.95E-05$ & $1.04E-05$ & $2.56E-04$ & $5.38$ & $6.12$ & $5.03$ &   \\	\hhline{|~|-||-|-|-||-|-|-||~|}
		& $120$ & $8.22E-06$ & $3.55E-06$ & $1.06E-04$ & $4.75$ & $5.89$ & $4.86$ &   \\	\hhline{|~|-||-|-|-||-|-|-||~|}
		& $140$ & $4.10E-06$ & $1.49E-06$ & $5.21E-05$ & $4.51$ & $5.61$ & $4.58$ &   \\	\hline\hline
		\multirow{5}{*}{$5$}& $20$ & $9.48E-03$ & $1.55E-02$ & $8.55E-02$ & - & - & - &  \multirow{5}{*}{$6$} \\	\hhline{|~|-||-|-|-||-|-|-||~|}
		& $40 $ & $3.84E-04$ & $1.86E-04$ & $3.06E-03$ & $4.63$ & $6.39$ & $4.80$ &   \\	\hhline{|~|-||-|-|-||-|-|-||~|}
		& $60 $ & $4.51E-05$ & $1.99E-05$ & $3.70E-04$ & $5.28$ & $5.51$ & $5.22$ &   \\	\hhline{|~|-||-|-|-||-|-|-||~|}
		& $80 $ & $8.41E-06$ & $3.86E-06$ & $6.26E-05$ & $5.84$ & $5.71$ & $6.17$ &   \\	\hhline{|~|-||-|-|-||-|-|-||~|}
		& $100$ & $2.26E-06$ & $9.46E-07$ & $1.74E-05$ & $5.89$ & $6.30$ & $5.72$ &   \\	\hline\hline
		\multirow{5}{*}{$6$}& $30$ & $3.59E-04$ & $3.15E-04$ & $3.17E-03$ & - & - & - &  \multirow{5}{*}{$7$} \\	\hhline{|~|-||-|-|-||-|-|-||~|}
		& $40 $ & $6.43E-05$ & $5.21E-05$ & $5.83E-04$ & $5.98$ & $6.25$ & $5.89$ &   \\	\hhline{|~|-||-|-|-||-|-|-||~|}
		& $50 $ & $1.46E-05$ & $7.03E-06$ & $1.29E-04$ & $6.65$ & $8.98$ & $6.75$ &   \\	\hhline{|~|-||-|-|-||-|-|-||~|}
		& $60 $ & $2.75E-06$ & $1.65E-06$ & $2.93E-05$ & $9.17$ & $7.95$ & $8.15$ &   \\	\hhline{|~|-||-|-|-||-|-|-||~|}
		& $70 $ & $7.71E-07$ & $5.32E-07$ & $9.39E-06$ & $8.24$ & $7.34$ & $7.37$ &   \\	\hline\hline
		\multirow{5}{*}{$7$}& $10$ & $3.10E-02$ & $5.16E-02$ & $1.65E-01$ & - & - & - &  \multirow{5}{*}{$8$} \\	\hhline{|~|-||-|-|-||-|-|-||~|}
		& $20 $ & $8.84E-04$ & $1.30E-03$ & $1.03E-02$ & $5.13$ & $5.31$ & $4.00$ &   \\	\hhline{|~|-||-|-|-||-|-|-||~|}
		& $30 $ & $7.54E-05$ & $5.88E-05$ & $6.56E-04$ & $6.07$ & $7.64$ & $6.80$ &   \\	\hhline{|~|-||-|-|-||-|-|-||~|}
		& $40 $ & $8.60E-06$ & $4.36E-06$ & $8.29E-05$ & $7.55$ & $9.05$ & $7.19$ &   \\	\hhline{|~|-||-|-|-||-|-|-||~|}
		& $50 $ & $8.27E-07$ & $6.61E-07$ & $1.12E-05$ & $10.50$ & $8.45$ & $8.98$ &   \\	\hline
	\end{tabular}
	\caption{Solitary wave over a flat bottom. $L^{2}$ errors and convergence rates obtained for $N\in\left\lbrace 3,4,5,6,7\right\rbrace$ at $t_{\mathrm{end}}=2$.}
	\label{tab:SW_accuracy}
\end{table}

\subsubsection{Solitary wave over a step}
As second test case, we consider the solitary wave propagating over a step benchmark, introduced in \cite{SRT87}. We define the computational domain $\Omega=[-16,16]\times[-1,1]$ and a step shaped obstacle of height $H_{\mathrm{obs}}=0.1$ at $x=0$:
\begin{equation}
z_{b}= 0.05 \left( \mathrm{erf}(8 x)+1\right).
\end{equation}
As initial conditions, we set a soliton of amplitude $A=0.0365$ centered at $(x,y)=(-3,0)$ and a still water depth $H=0.2$. 
Periodic boundary conditions are imposed in $x$ and $y$ directions. 

\begin{figure}
	\centering
	\includegraphics[width=0.4\linewidth]{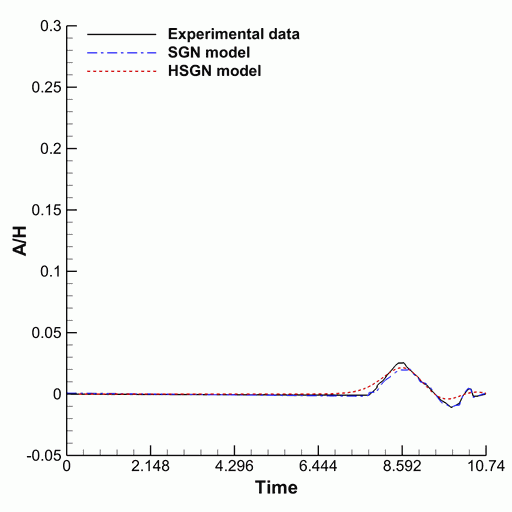}\hspace{0.1\linewidth}
	\includegraphics[width=0.4\linewidth]{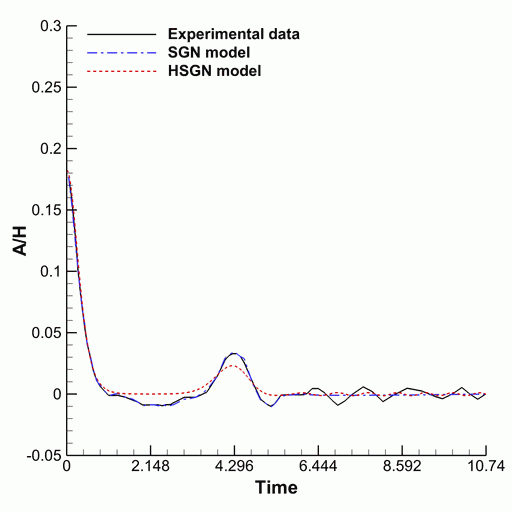}
	
	\includegraphics[width=0.4\linewidth]{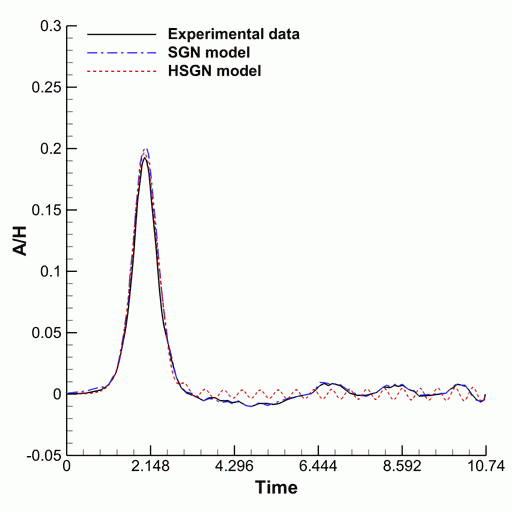}\hspace{0.1\linewidth}	
	\includegraphics[width=0.4\linewidth]{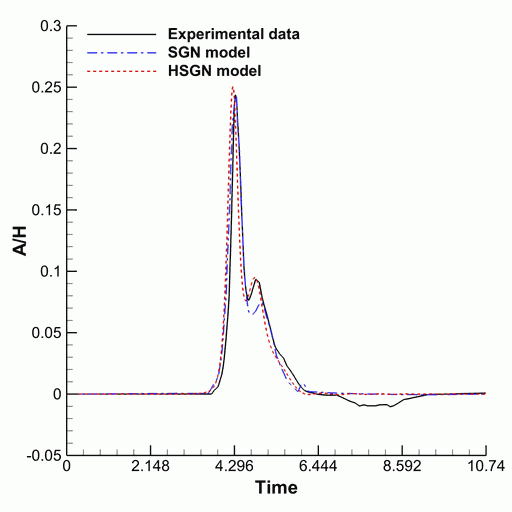}
	
	\includegraphics[width=0.4\linewidth]{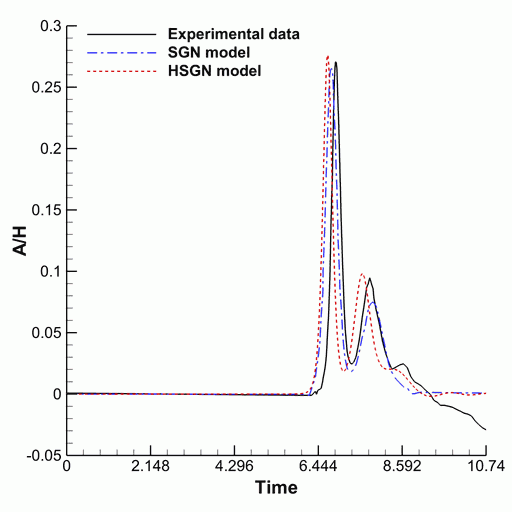}\hspace{0.1\linewidth}
	\includegraphics[width=0.4\linewidth]{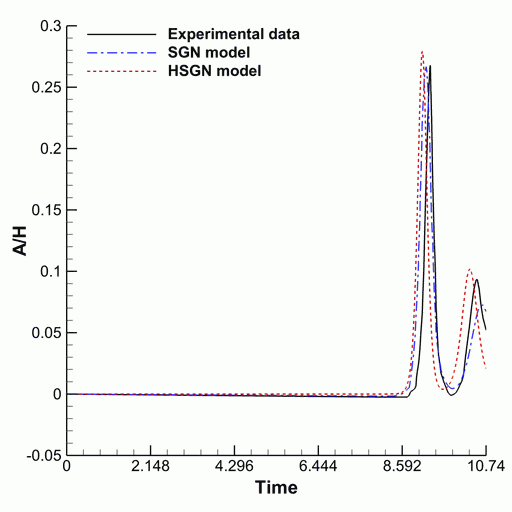}
	
	\caption{Solitary wave over a step. Comparison of the time evolution of numerical results obtained using ADER-DG against the numerical and experimental data given in \cite{SRT87} at  locations $x\in\left\lbrace -9,-3,0,3,6,9 \right\rbrace$ (from left top to right bottom).}
	\label{fig:solstep_comparison}
\end{figure}

\begin{figure}
	\centering
	\includegraphics[width=0.32\linewidth]{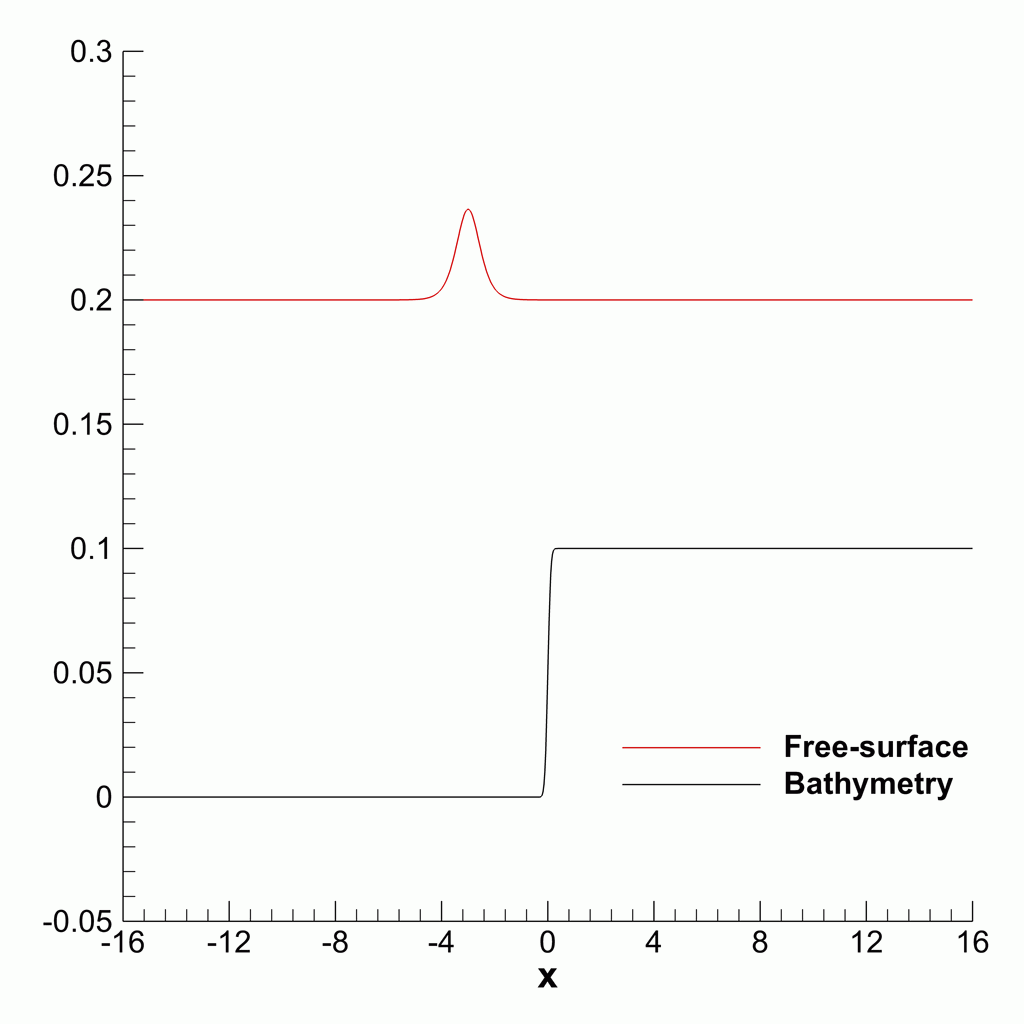}\hfill
	\includegraphics[width=0.32\linewidth]{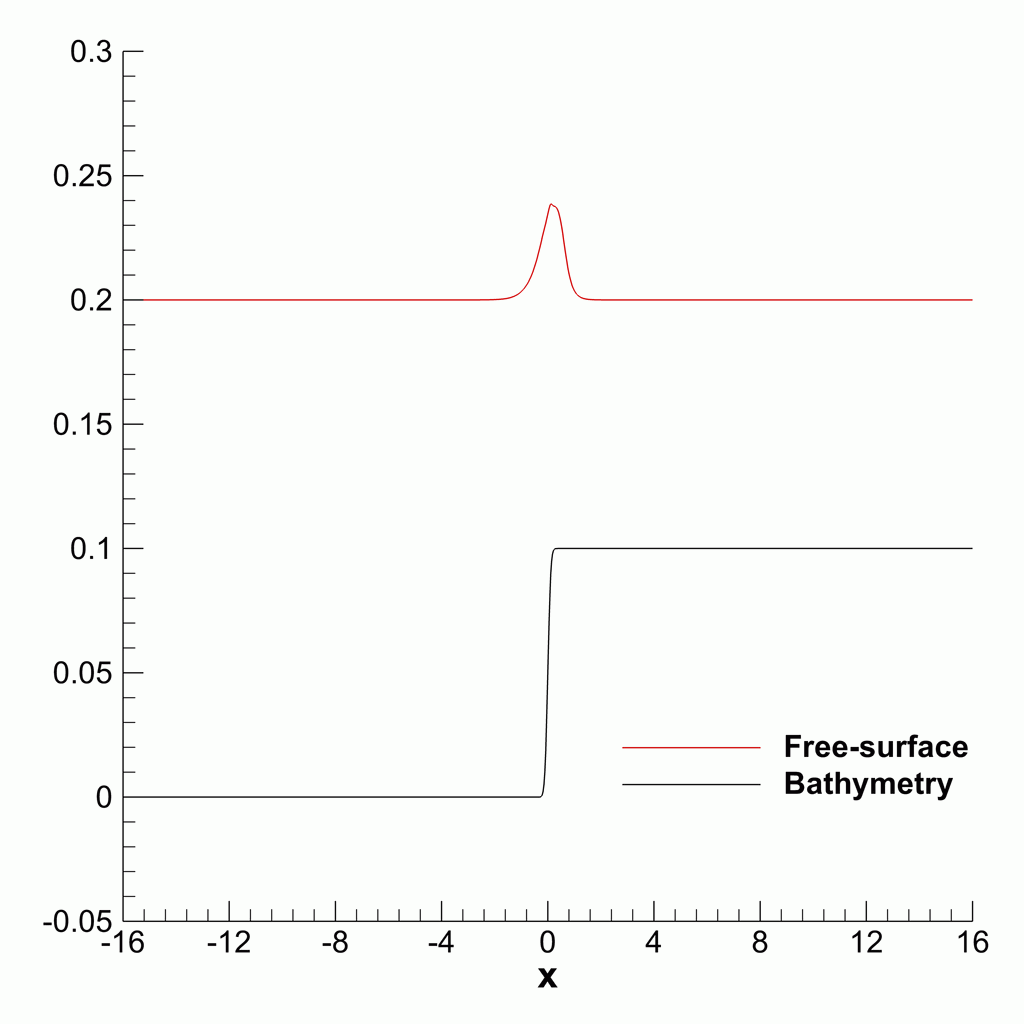}\hfill
	\includegraphics[width=0.32\linewidth]{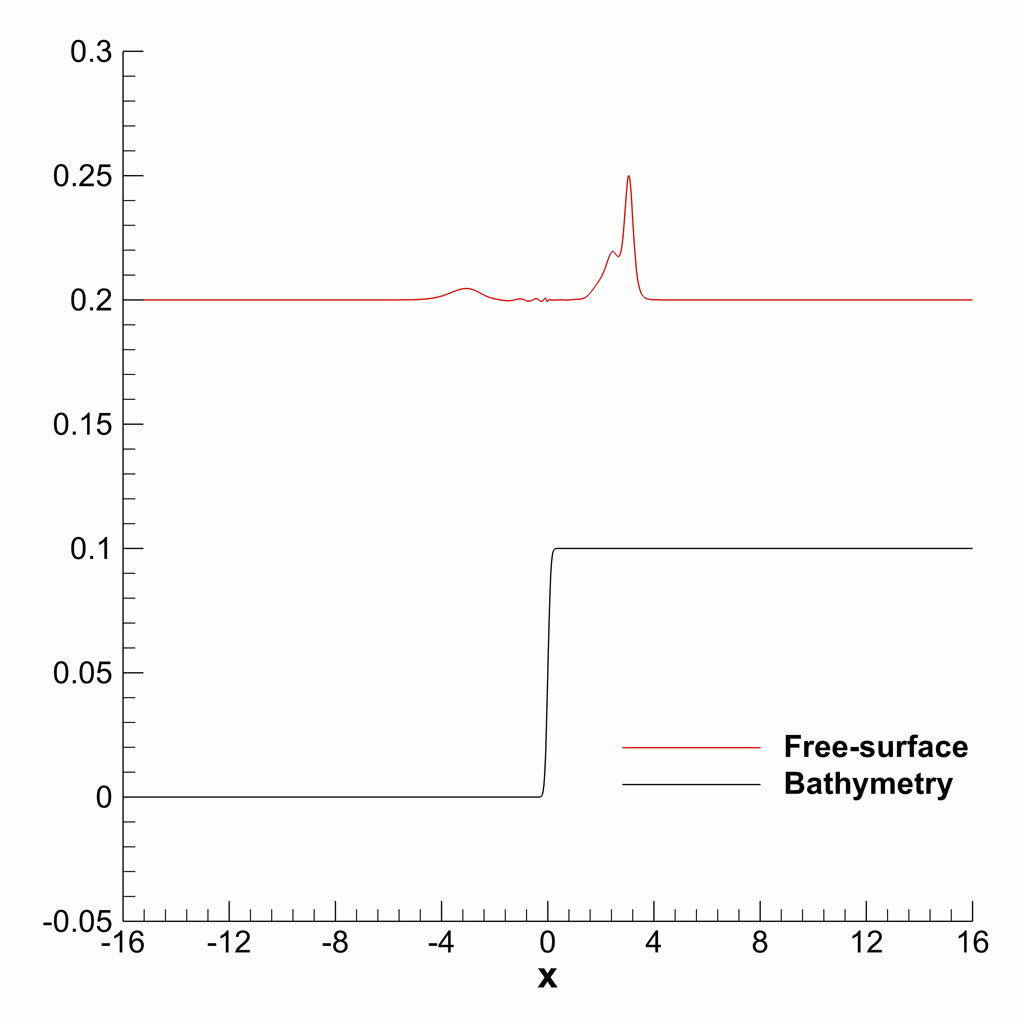}
	
	\includegraphics[width=0.32\linewidth]{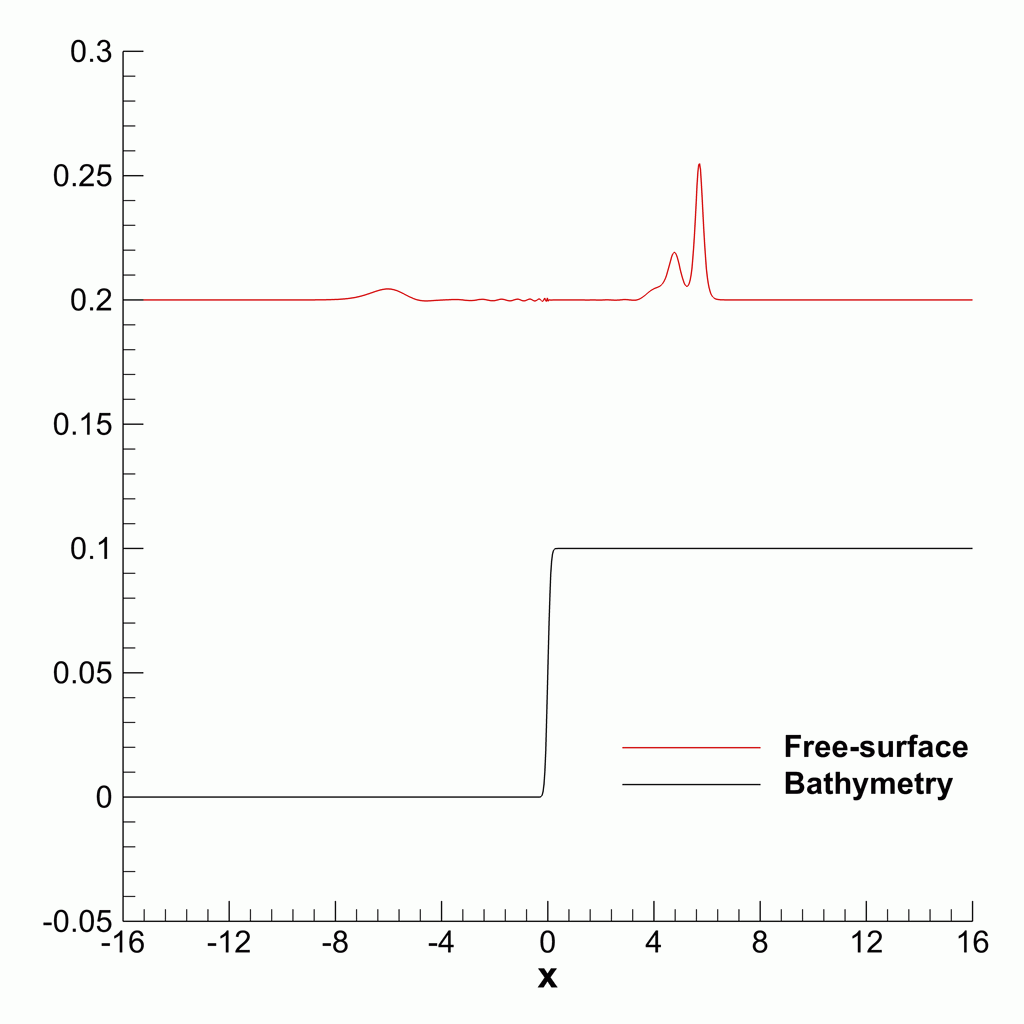}\hfill
	\includegraphics[width=0.32\linewidth]{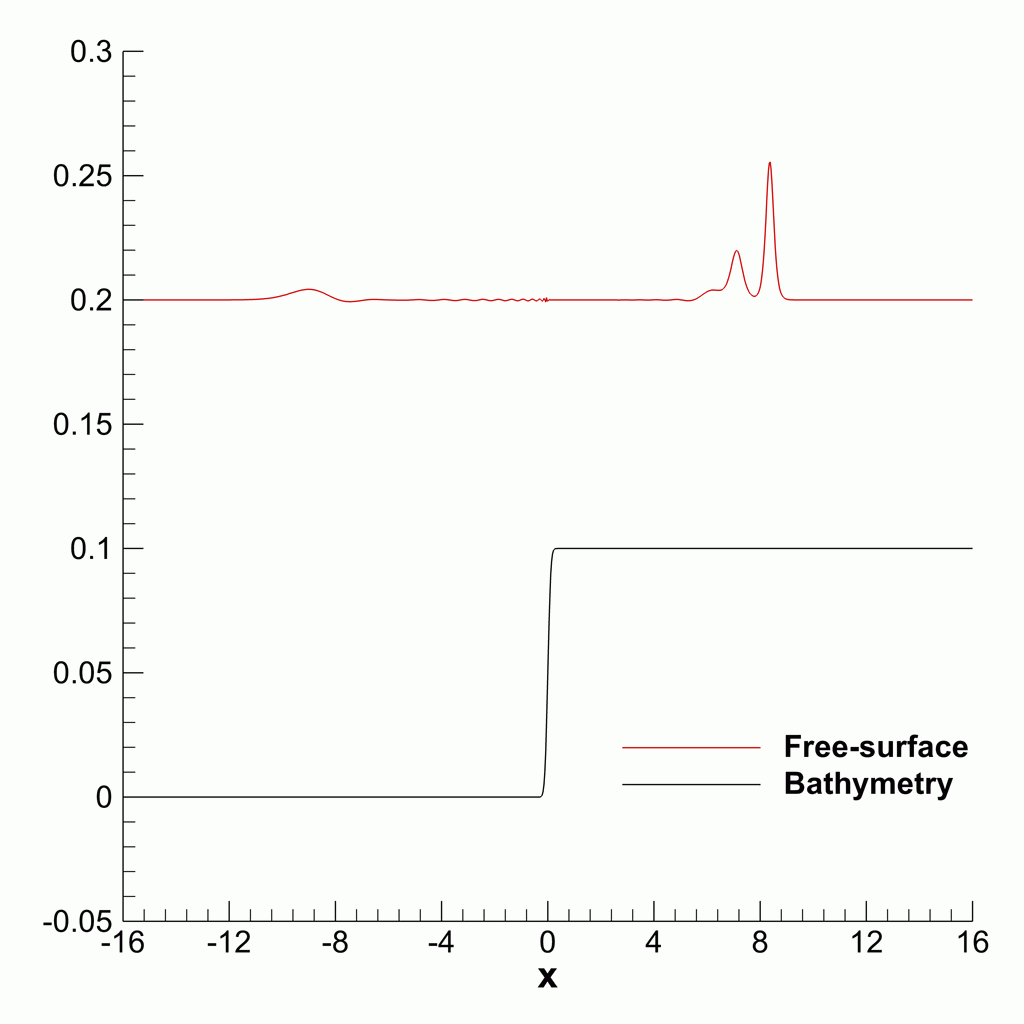}\hfill
	\includegraphics[width=0.32\linewidth]{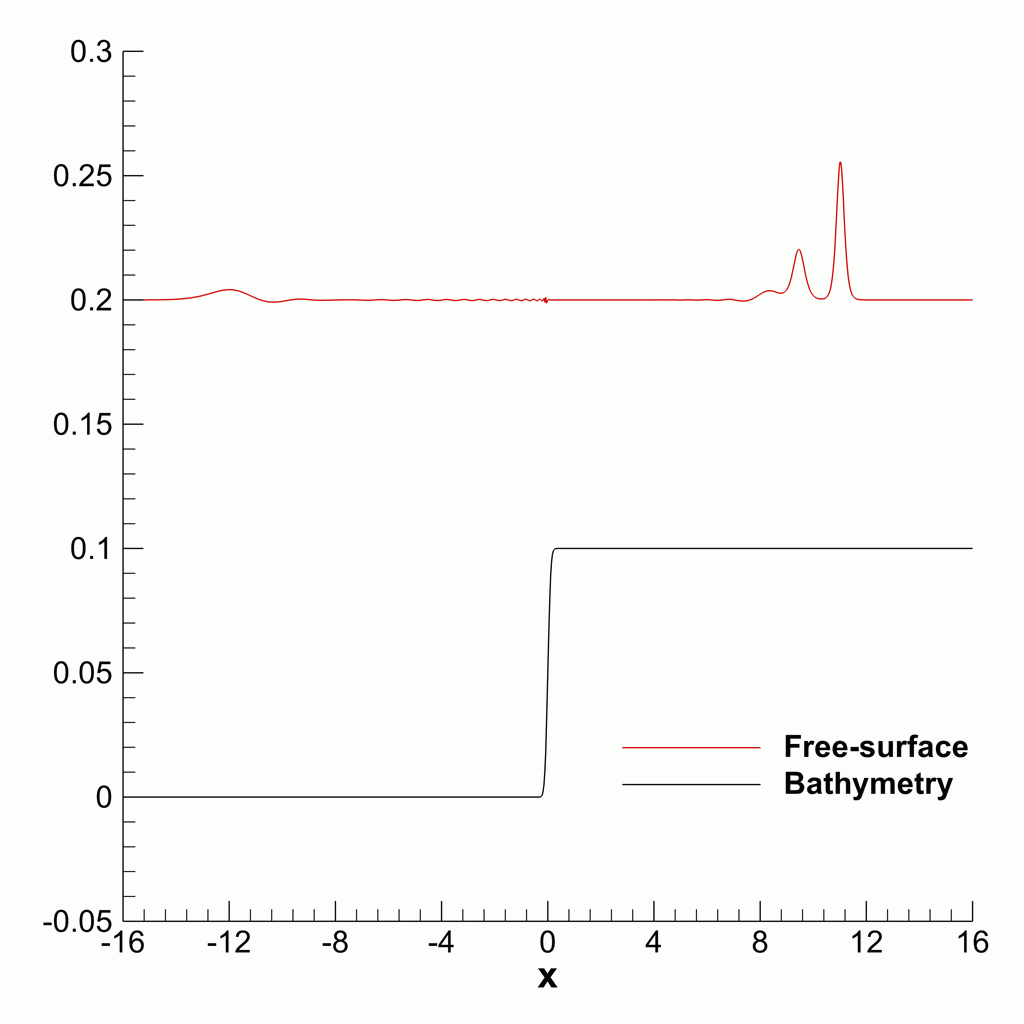}
	
	\caption{Solitary wave over a step. 1D cut along $y=0$ of the free surface and the bottom bathymetry at times \mbox{$t\in \left\lbrace 0,2.148,4.296,6.444,8.592,10.74 \right\rbrace$} (from left top to right bottom).}
	\label{fig:solstep_freesurface}
\end{figure}

The simulation is run using the HSGN model on a mesh made of $2000\times 2$ elements. Figure \ref{fig:solstep_freesurface} depicts a 1D free-surface cut at  \mbox{$t\in \left\lbrace 0,2.148,4.296,6.444,8.592,10.74 \right\rbrace$}. Due to the presence of the step, the amplitude of the soliton starts growing until it splits into two transmitted waves. Moreover, a reflected wave starts to propagate in the opposite direction with respect to the soliton, followed by a train of small dispersive waves. Notice that small spurious oscillations appear in correspondence to obstacle location $x_{obs}=0$. As already pointed out in \cite{BBBD20}, this is due to the fact that the model HSGN is rigorously valid only in the presence of a slowly varying bottom in space, which could create some problems when a strongly varying bottom topography is considered (as in the present test). 
To compare the numerical results obtained with the experimental data provided in  \cite{SRT87}, we compute the ratio between the wave amplitude and the still water depth,
\begin{equation}
\frac{A}{H}= \frac{h-H}{H},
\end{equation}
at seven different locations $x\in\left\lbrace -9,-6,-3,0,3,6,9 \right\rbrace$. In Figure \ref{fig:solstep_comparison}, we observe that the results obtained match pretty well the experimental data, improving the numerical results already presented in \cite{SRT87}. 
In the first three plots, we can observe a good agreement of the amplitude and location of the reflected waves, even if the already mentioned spurious oscillations at the obstacle location can be detected.
Also transmitted waves are properly captured, including a third transmitted wave that is missing in the numerical results in \cite{SRT87}.

\subsection{Coupled model tests}
The last numerical tests aim at showing the behaviour of the proposed methodology for the simulation of the coupled problem. Two different initial water wave profiles are considered: a solitary wave and a train of sinusoidal waves.

\subsubsection{Solitary wave}
\begin{figure}
	\centering	
	\hspace{0.015\linewidth}\includegraphics[width=0.27\linewidth,align=t]{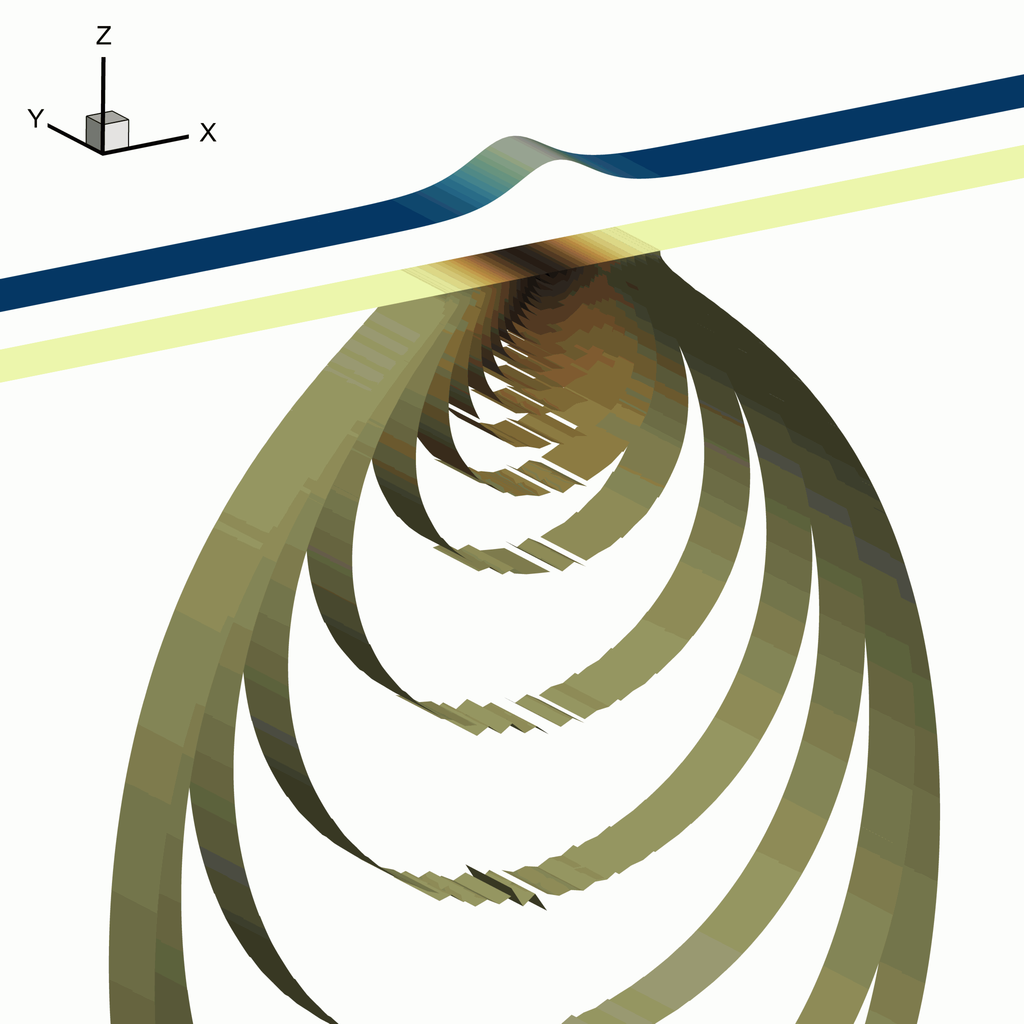}\hspace{0.025\linewidth}
	\includegraphics[width=0.3\linewidth,align=t]{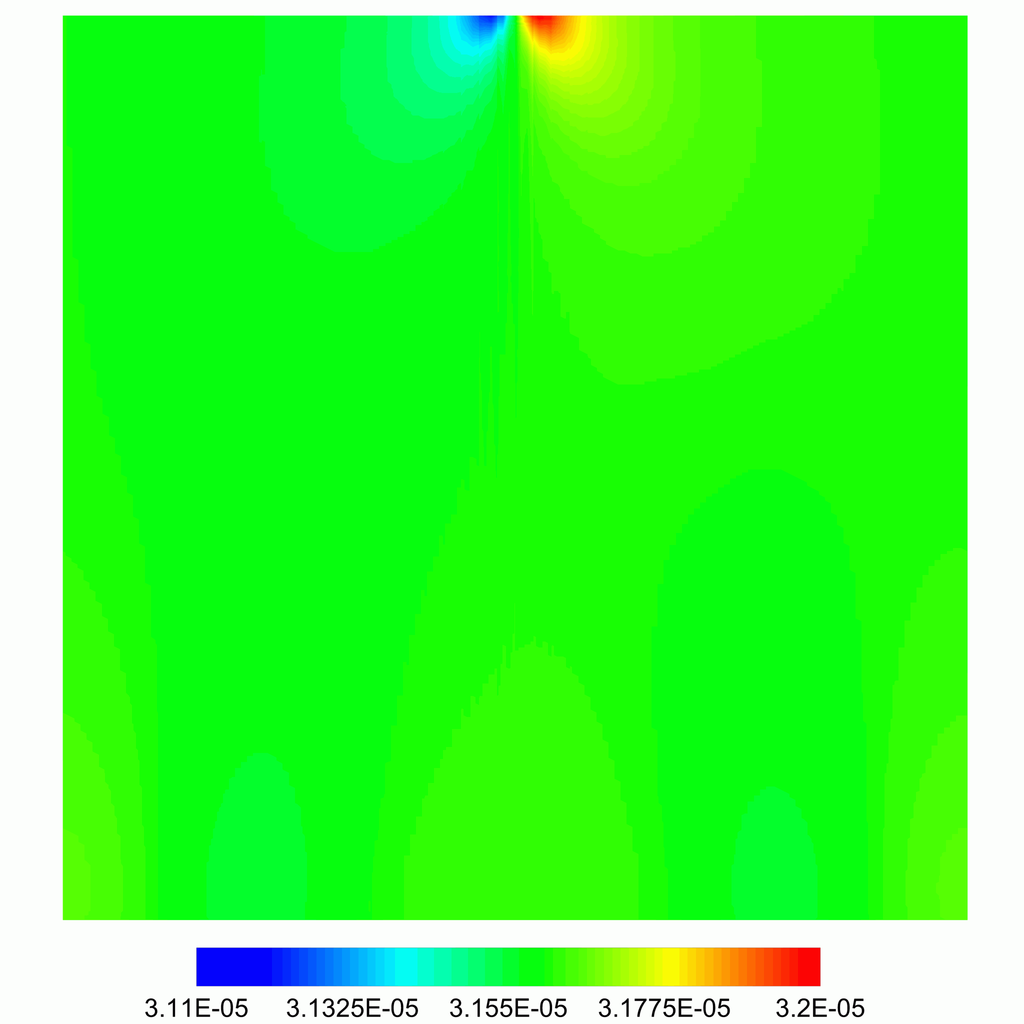}\hspace{0.02\linewidth}
	\includegraphics[width=0.3\linewidth,align=t]{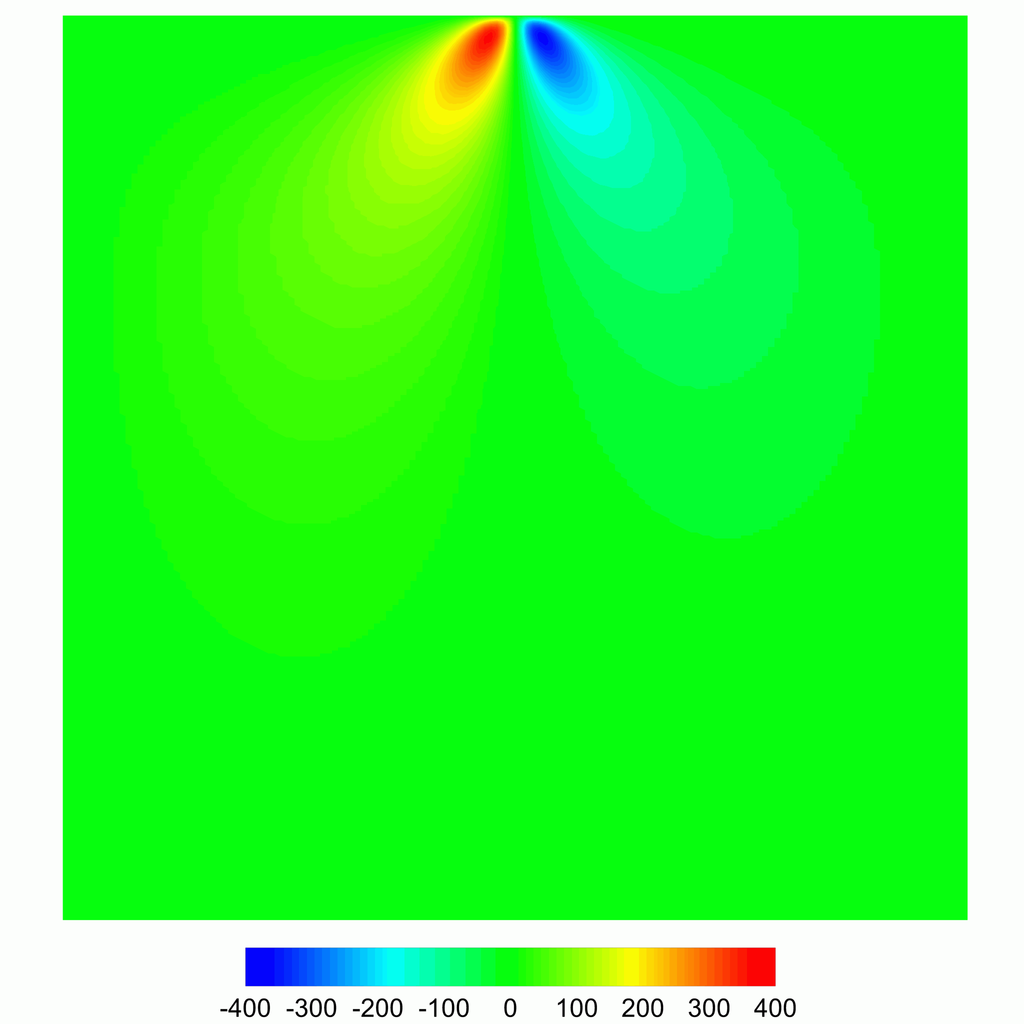}
	
	\vspace{0.05\linewidth}
	\includegraphics[width=0.3\linewidth]{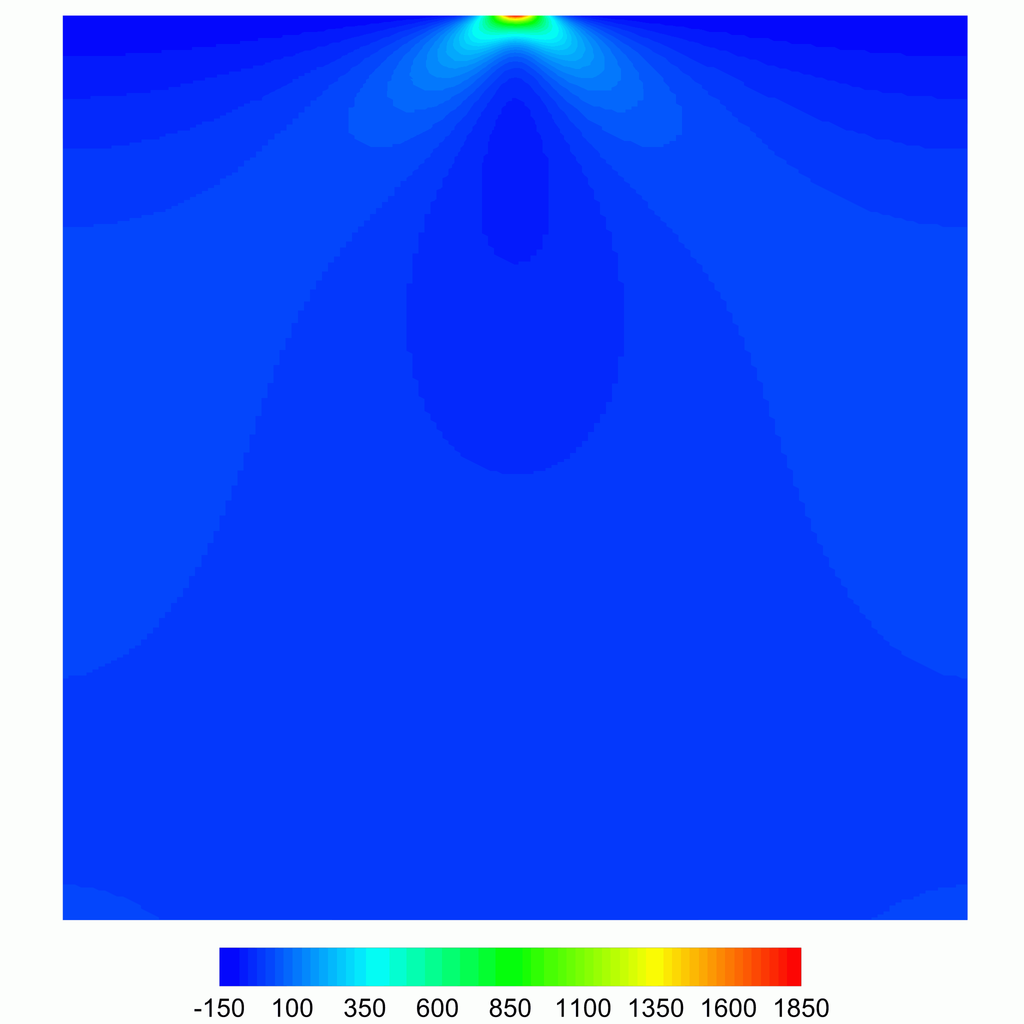}\hspace{0.02\linewidth}	
	\includegraphics[width=0.3\linewidth]{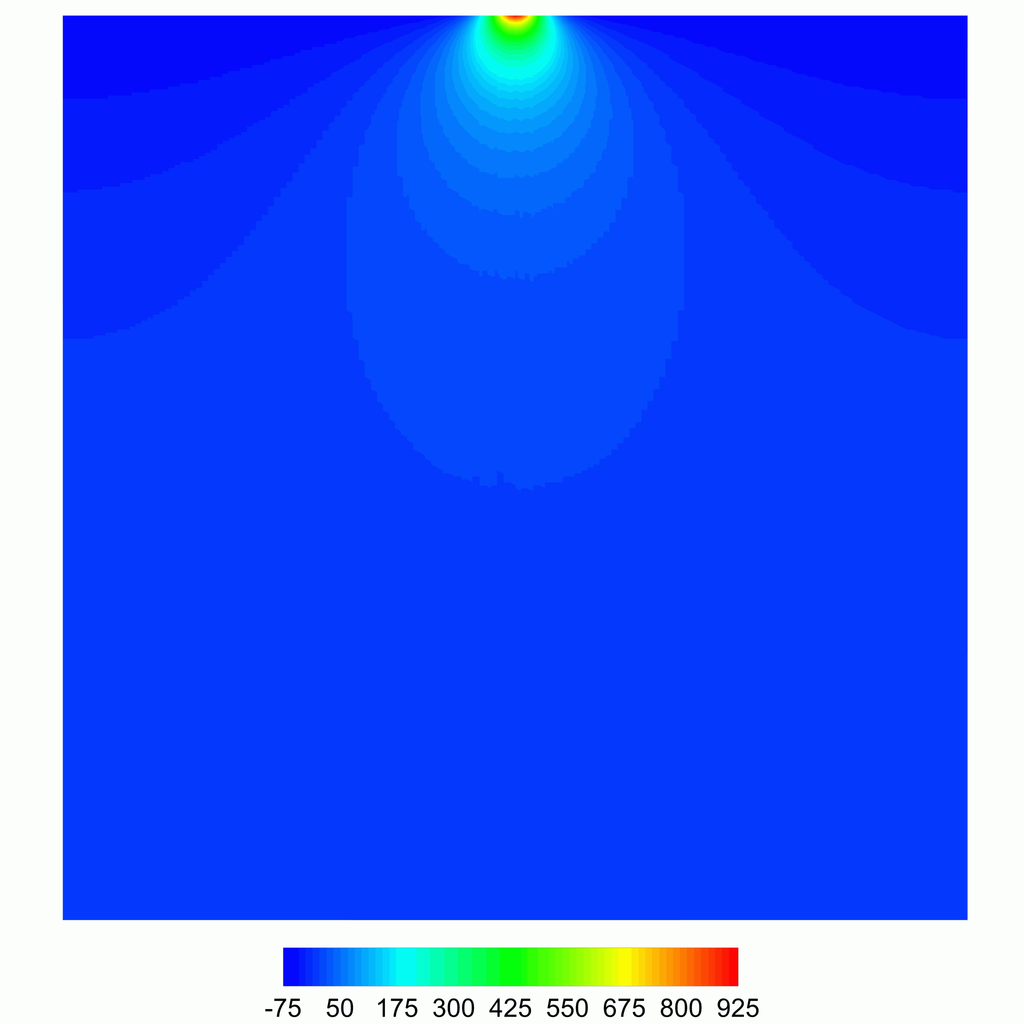}\hspace{0.02\linewidth}
	\includegraphics[width=0.3\linewidth]{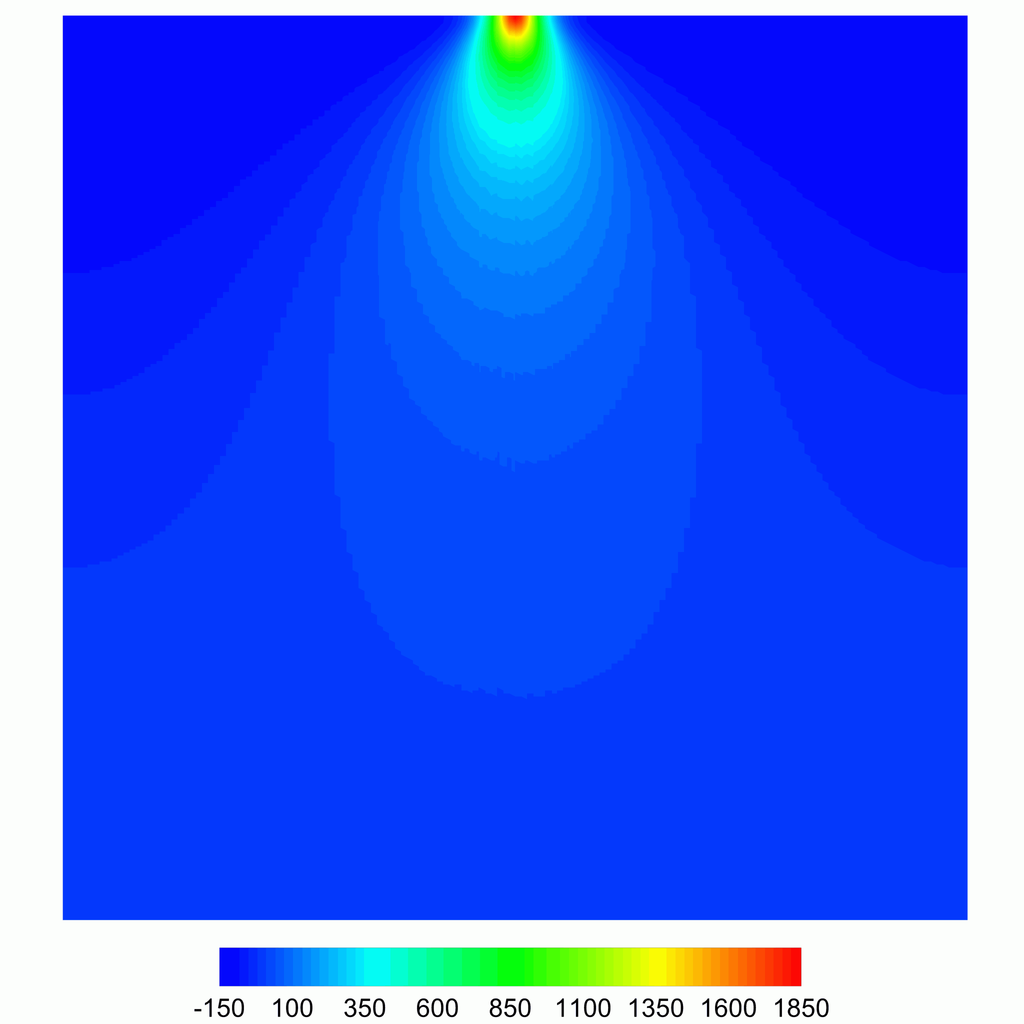}
	
	\caption{Solitary wave. Solution obtained at time $t=58.3$. Left top: 3D extrusion of the water depth, $h$, isosurfaces of $\sigma_{zz}$ and contour plot of $\sigma_{zz}$ at slice $z=0$. From middle top to bottom right: contour plots at plane $y=0$ of $w$, $\sigma_{xz}$, $\sigma_{xx}$, $\sigma_{yy}$, and $\sigma_{zz}$, respectively.}

	\label{fig:SoWC_3Dplot}
\end{figure}

\begin{figure}
	\centering	
	\includegraphics[width=0.4\linewidth]{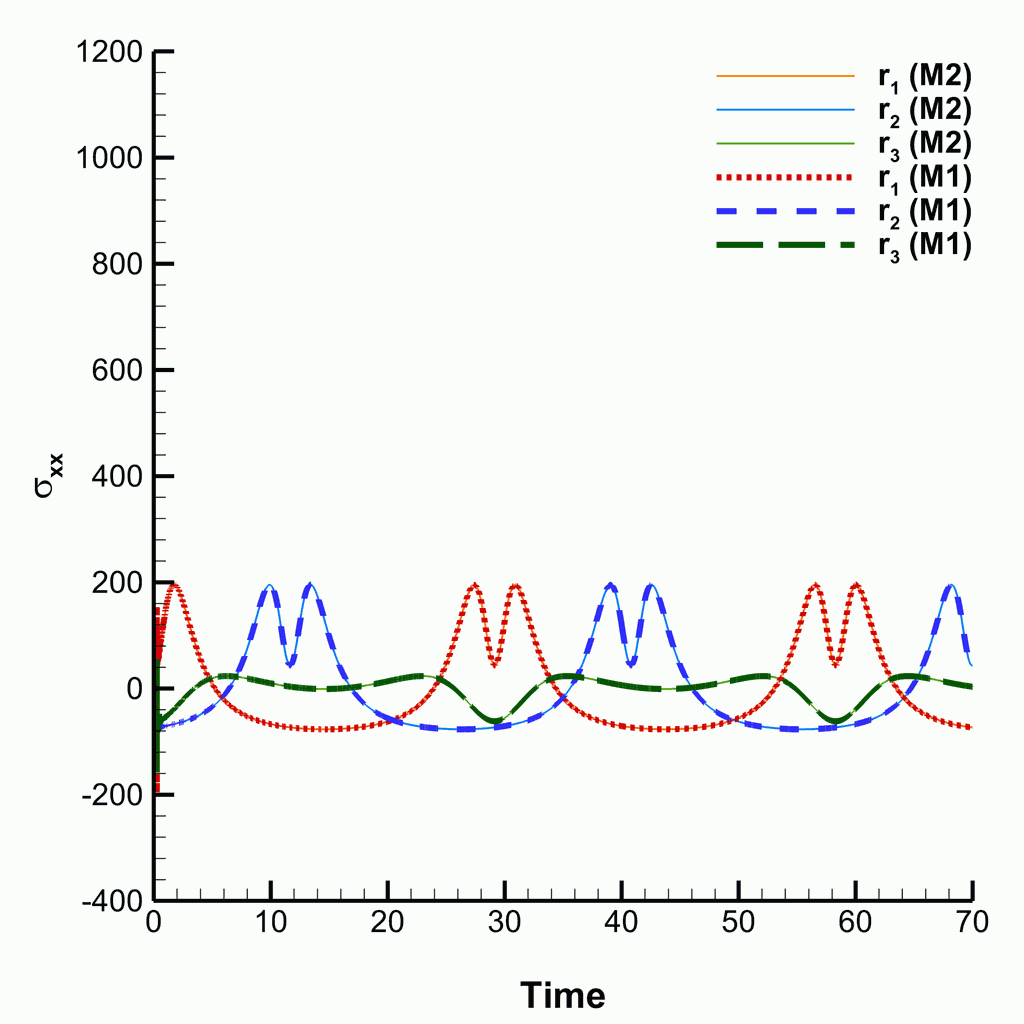}\hspace{0.05\linewidth}
	\includegraphics[width=0.4\linewidth]{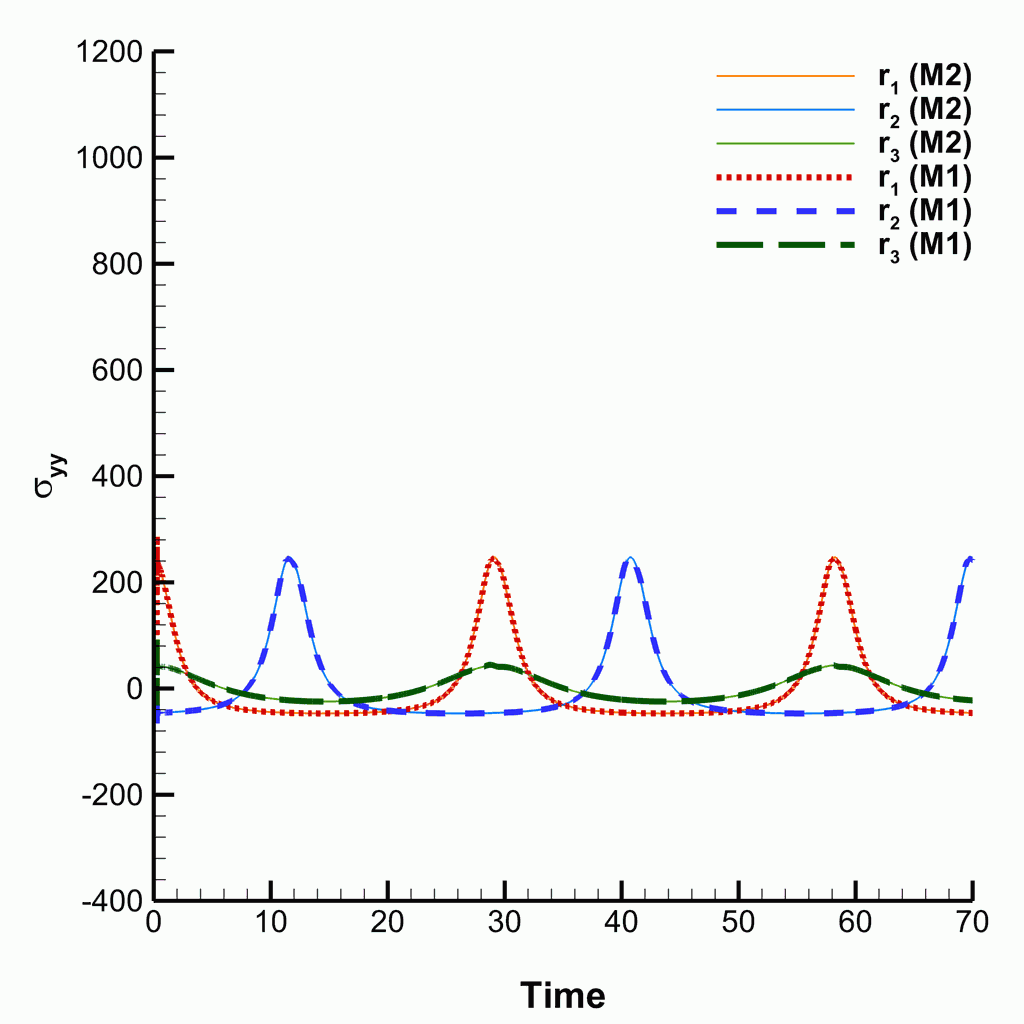}
	
	\vspace{0.05\linewidth}
	\includegraphics[width=0.4\linewidth]{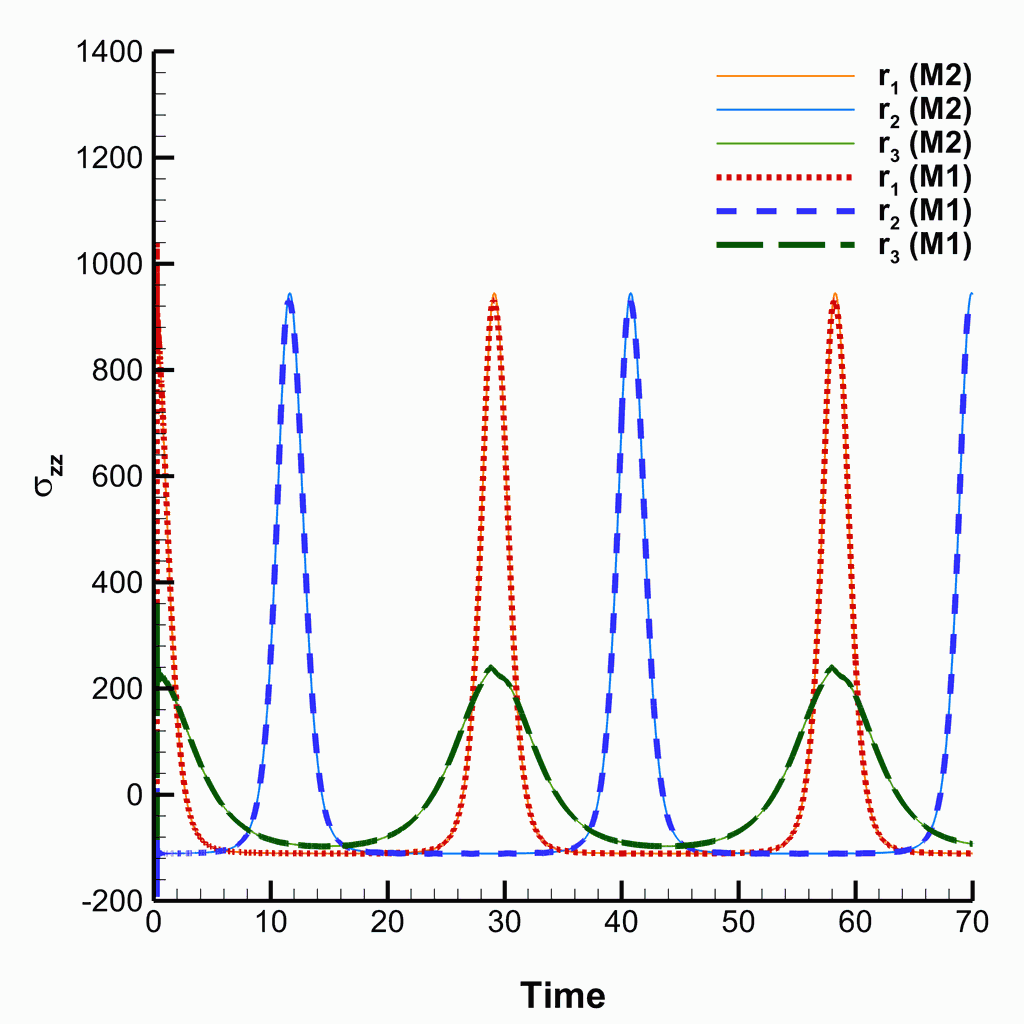}\hspace{0.05\linewidth}	
	\includegraphics[width=0.4\linewidth]{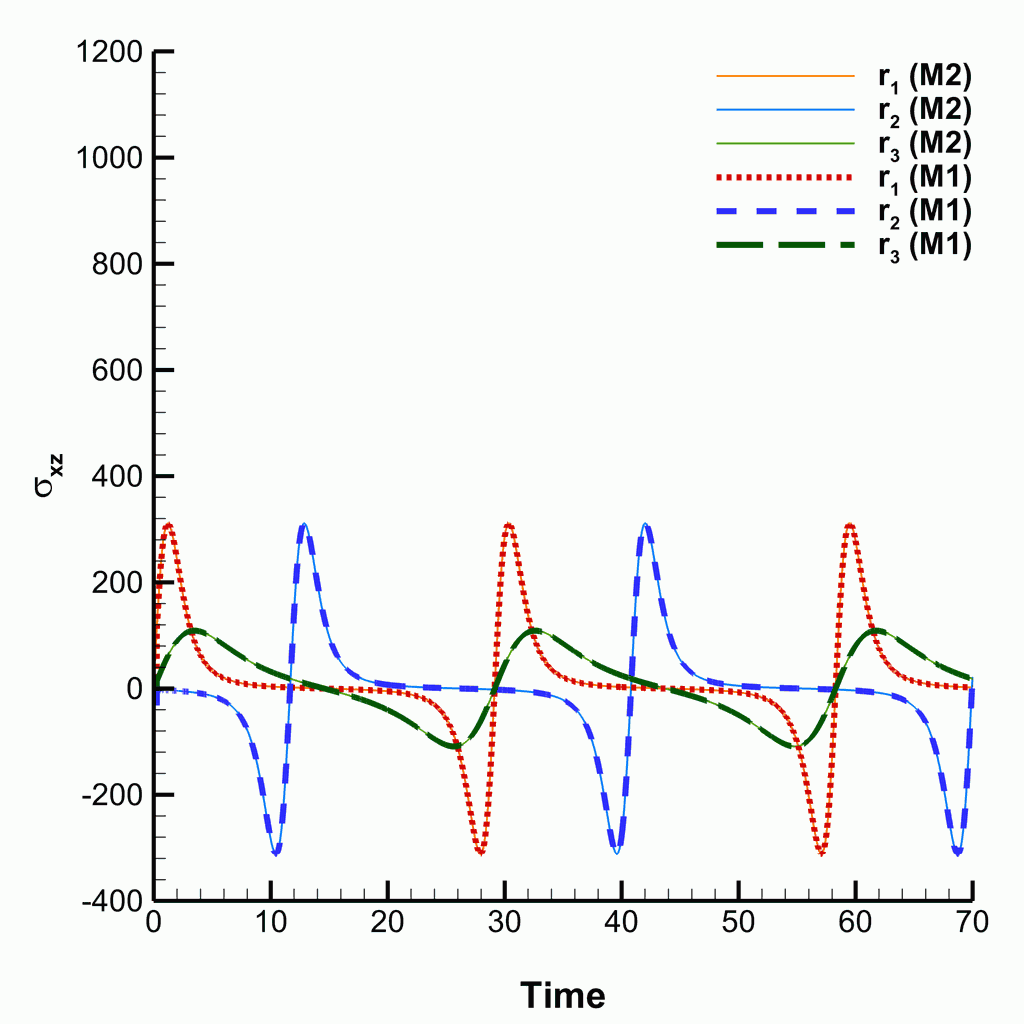}
	
	\caption{Solitary wave. Seismograms at receivers $\mathbf{x}_{r_{1}}=\left(0,0,-5\right)$ (M1: red dotted line, M2: orange line), $\mathbf{x}_{r_{2}}=\left(40,0,-5\right)$ (M1: blue dashed line, M2: cyan line), and $\mathbf{x}_{r_{3}}=\left(0,0,-20\right)$ (M1: green long dashed line; M2: light green line) obtained using ADER-DG $\mathcal{O}4$. From left top to right bottom: $\sigma_{xx}$, $\sigma_{yy}$, $\sigma_{zz}$, and $\sigma_{xz}$.}
	
	\label{fig:SoWC_pickpoints}
\end{figure}

\begin{figure}
	\centering	
	\includegraphics[width=0.4\linewidth]{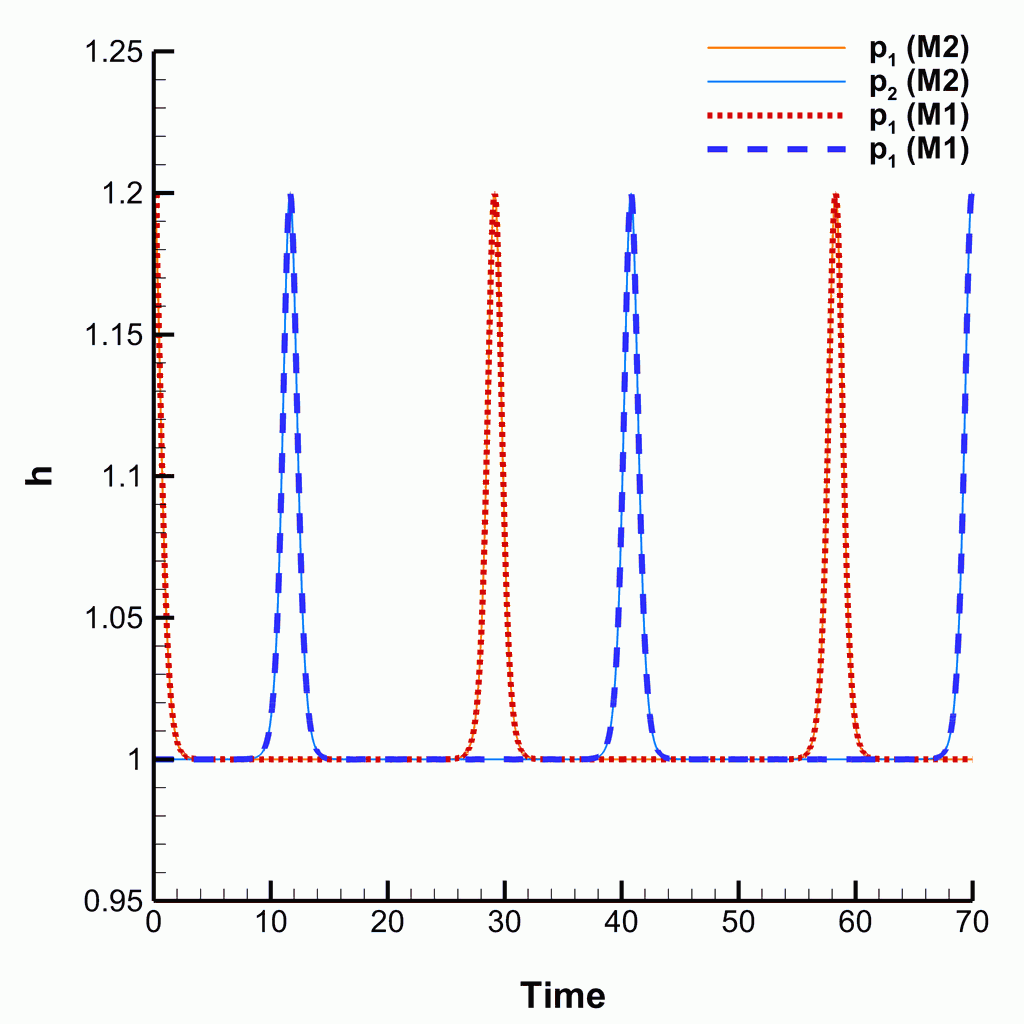}\hspace{0.05\linewidth}
	\includegraphics[width=0.4\linewidth]{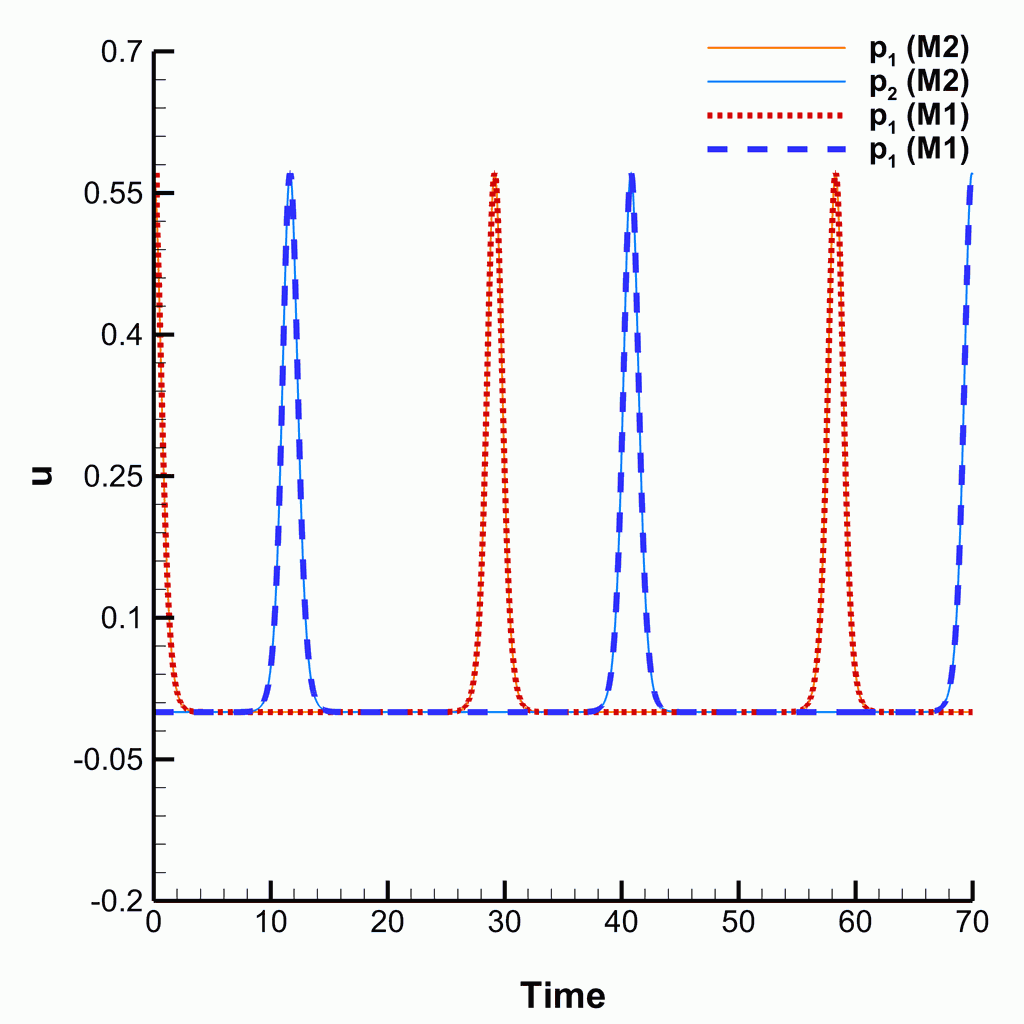}
	
	\vspace{0.05\linewidth}
	\includegraphics[width=0.4\linewidth]{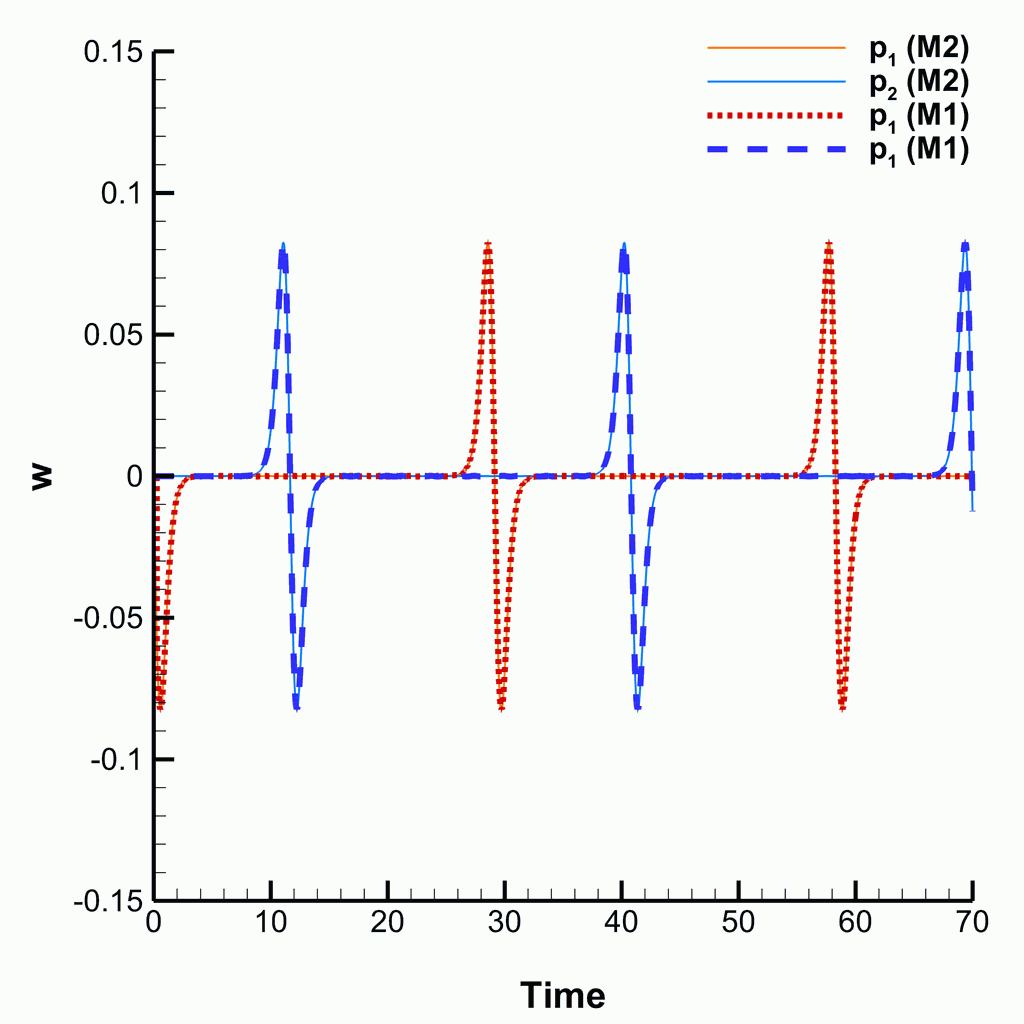}\hspace{0.05\linewidth}
	\includegraphics[width=0.4\linewidth]{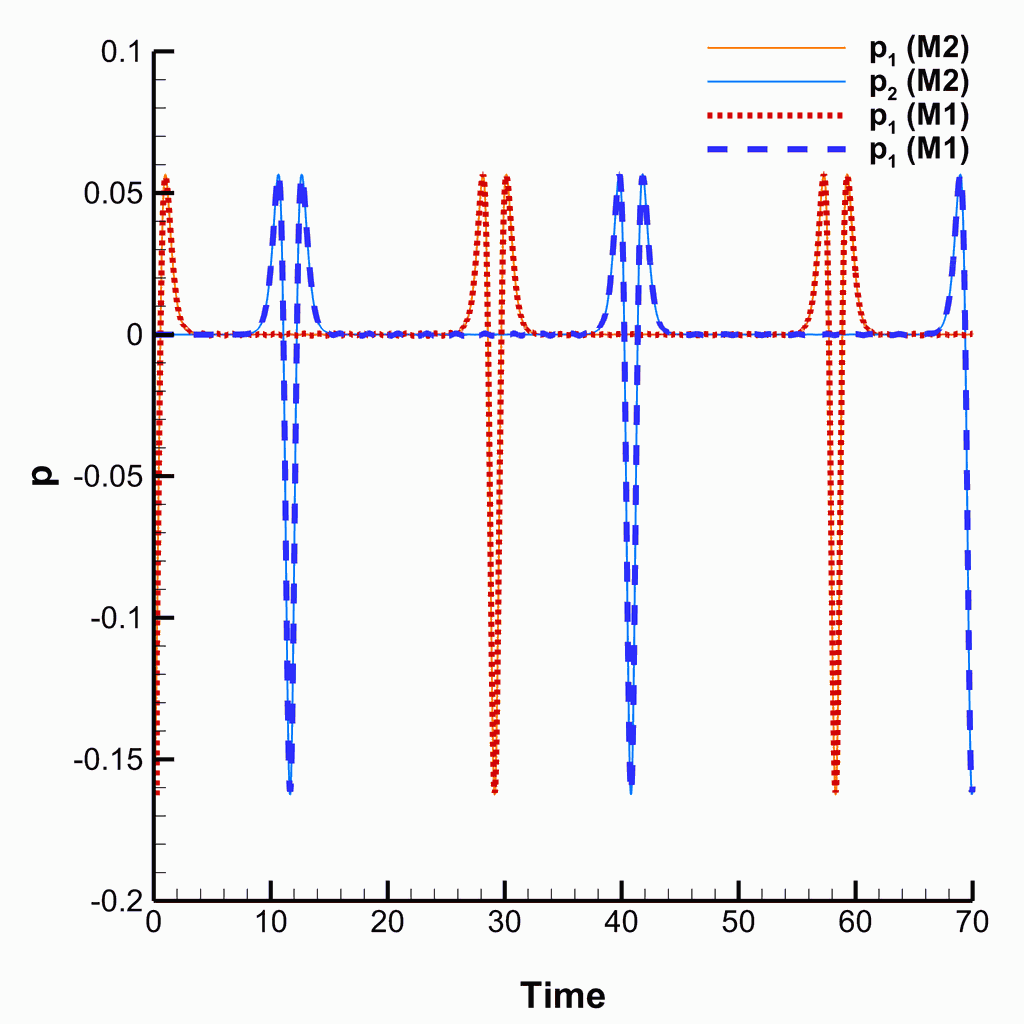}

	\caption{Solitary wave. Time evolution of the main variables of the HSGN model at pick points  $\mathbf{x}_{p_{1}}=\left(0,0\right)$ (M1 red line) and $\mathbf{x}_{p_{2}}=\left(40,0 \right)$ (M1 blue line) obtained using ADER-DG $\mathcal{O}4$. From left top to right bottom: $h$, $u$, $w$, and $p$.}
	
	\label{fig:SoWC_2Dpickpoints}
\end{figure}

We first study the seismic waves generated in the solid domain by the propagation of a soliton on the water surface. 
We consider the computational domains $\Omega_{w}=[-50,50]\times [-2,2]$ for the fluid and $\Omega_e=[-50,50]\times [-2,2]\times [-100,0]$ for the solid.
As initial condition for the HSGN model we employ the planar solitary wave over a flat bottom  already analysed in Section \ref{sec:SoWFB}. The linear elasticity model is initiated with zero values. The simulation is run using a mesh of $50\times 2 \times 100$ elements on $\Omega_e$ and $100\times 4$ on $\Omega_{w}$, (M1). 
Periodic boundary conditions are set in $x$ and $y$ directions whereas a free surface boundary condition is defined on the bottom of $\Omega_e$.
The simulation is run up to time $t_{\mathrm{end}}=70$, which corresponds to $2.4$ revolutions of the soliton.  Note that the simulation can also be seen as a train of solitons, which do not interact among them thanks to the great length of the domain.  The solution obtained at $t=58.3$ is depicted in Figure \ref{fig:SoWC_3Dplot}. The first plot of the figure shows a good correspondence between the position of the water surface wave and the seismic wave propagating on the sea bottom, which validates the coupling methodology.
To better analyse the results obtained, we place three receivers at  $\mathbf{x}_{r_{1}}=\left(0,0,-5\right)$, $\mathbf{x}_{r_{2}}=\left(40,0,-5\right)$, and $\mathbf{x}_{r_{3}}=\left(0,0,-20\right)$.
Figure \ref{fig:SoWC_pickpoints} shows the time evolution of the main stress variables until the final simulation time, $t_{\mathrm{end}}=70$. The spurious oscillations generated at the beginning of the simulation, due to the homogeneous initial condition used, quickly disappear, leading to the smooth wave profile arising in response to the soliton. 
As expected, we observe that the magnitude of the stress tensor decreases as we get far from the sea bottom surface, while the wave front is reached at the same time instants. On the other hand, comparison of the seismograms obtained at $\mathbf{x}_{r_{1}}$ and $\mathbf{x}_{r_{2}}$ proves the conservation of the stress magnitude, as the seismic wave advances in the horizontal direction. We also include the results obtained with a finer grid (M2), made of $70\times 4 \times 140$ elements on 
$\Omega_e$ and $140\times 8$ on $\Omega_{w}$. Finally, in Figure \ref{fig:SoWC_2Dpickpoints} we show the time evolution of the main variables of the non-hydrostatic model.

\subsubsection{Sinusoidal wave}
\begin{figure}
	\centering	
	\includegraphics[width=\linewidth]{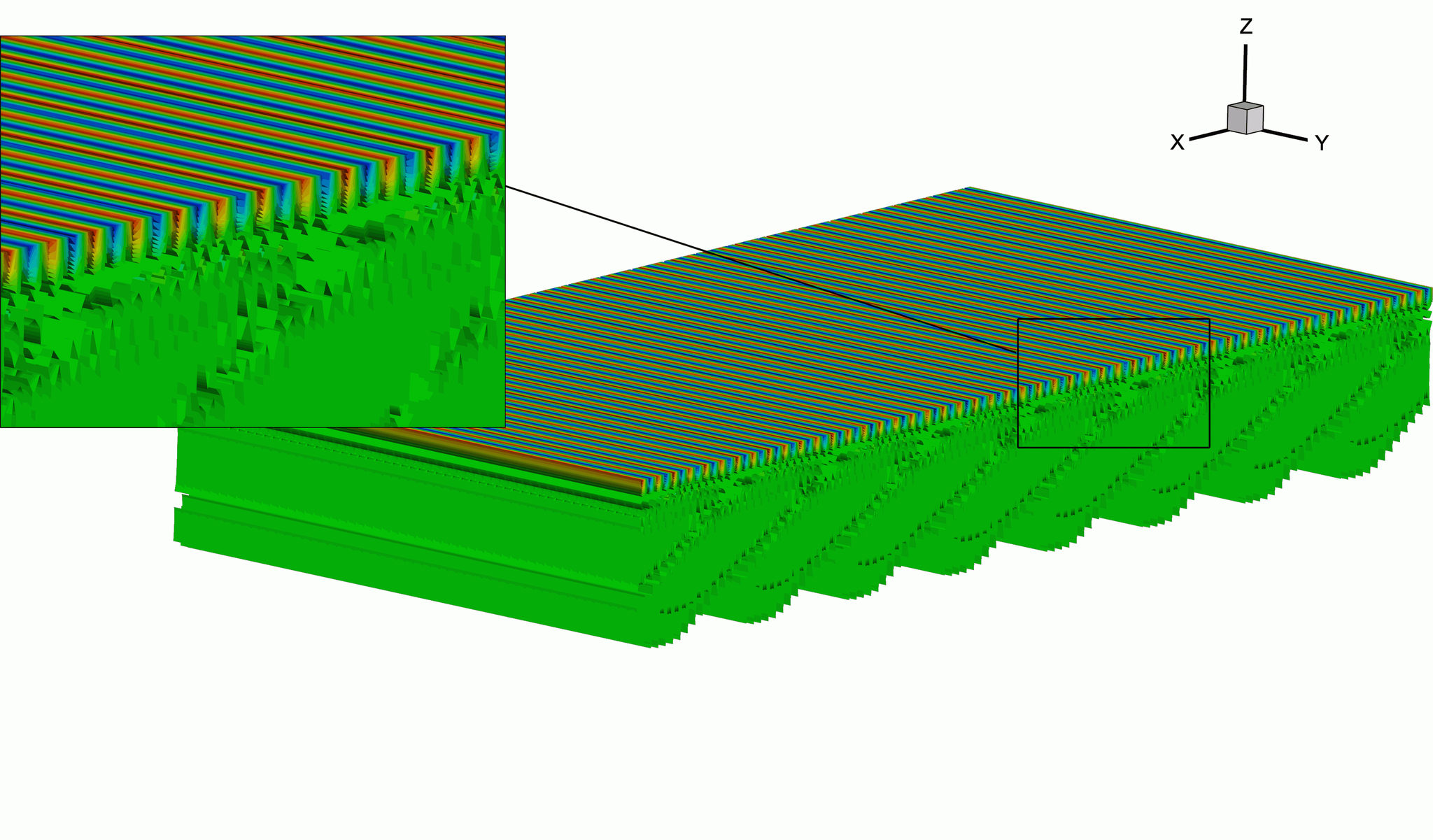}
	
	\caption{Sinusoidal wave. Isosurfaces of $\sigma_{zz}$, $\sigma_{zz}\in\left\lbrace -9, -8, -7, -6, -5, -4,-3,-2,-1,-10^{-2},10^{-2},\right.$ $\left.1,2,3,4,5,6,7,8,9 \right\rbrace$, at time $t=50$.}
	
	\label{fig:SiW_szz3D_subfig}
\end{figure}

\begin{figure}
	\centering	
	\includegraphics[width=0.32\linewidth]{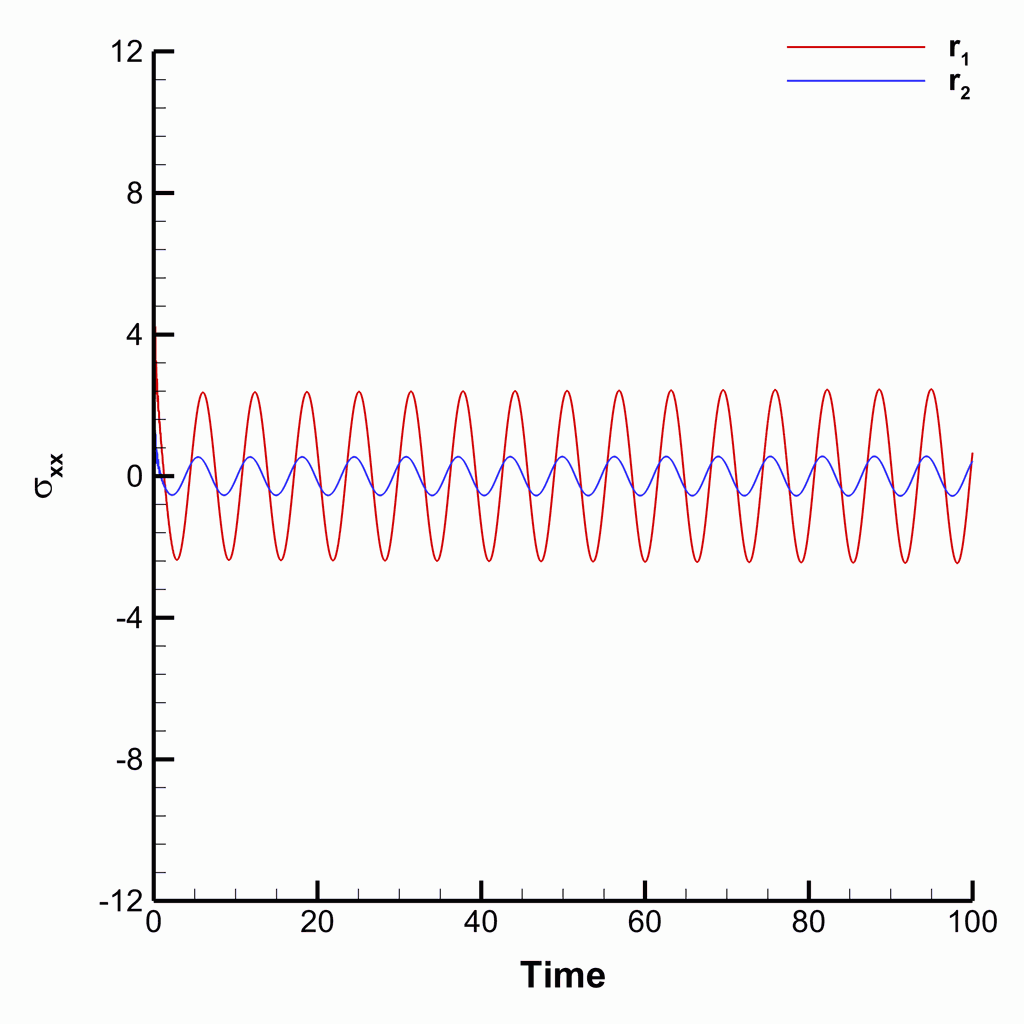}\hfill
	\includegraphics[width=0.32\linewidth]{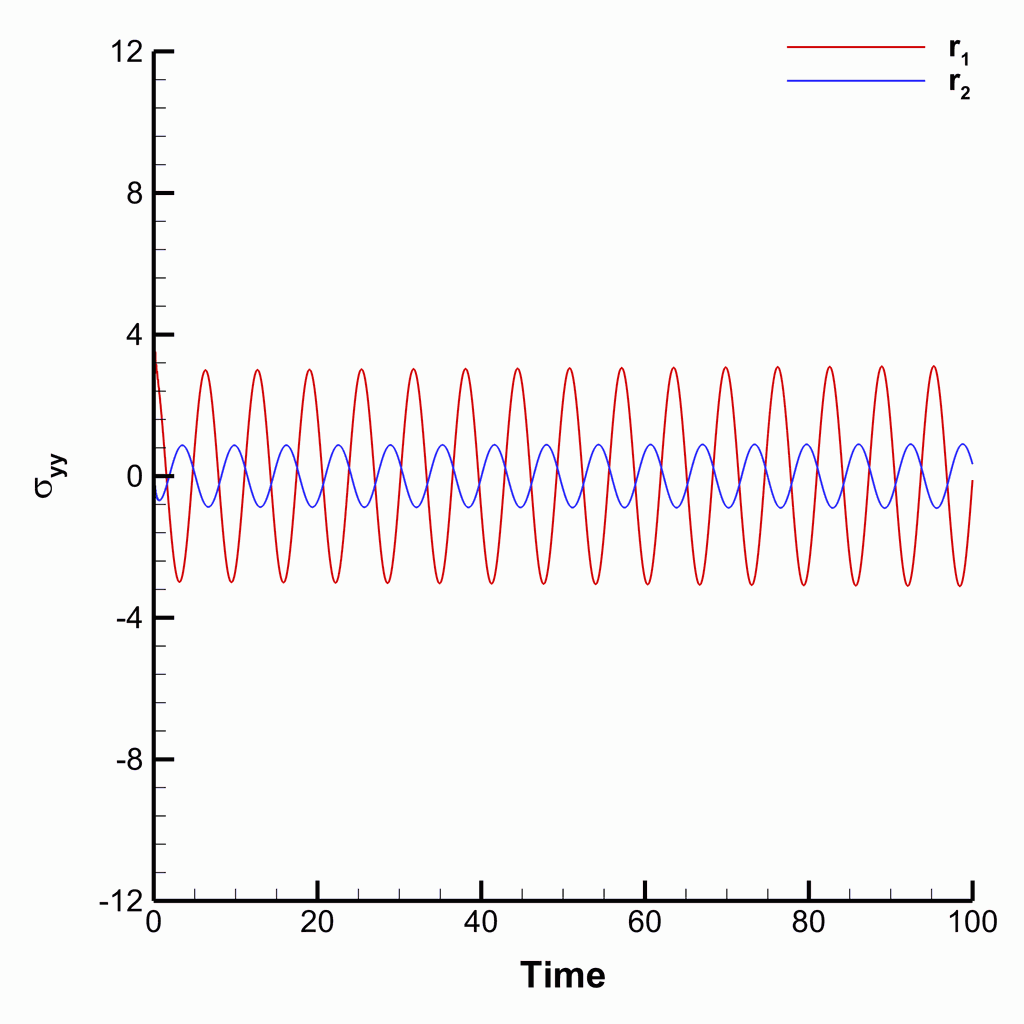}\hfill
	\includegraphics[width=0.32\linewidth]{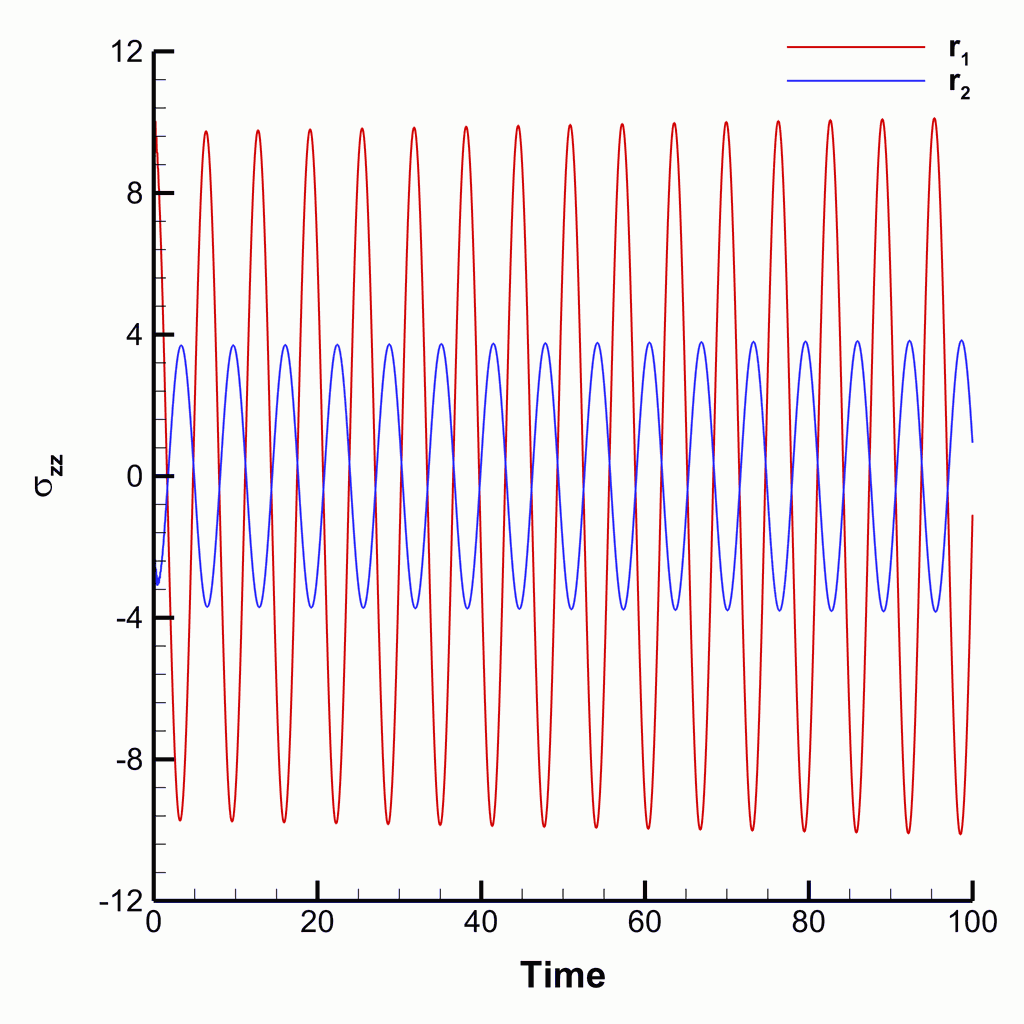}
	
	\vspace{0.05\linewidth}
	\includegraphics[width=0.32\linewidth]{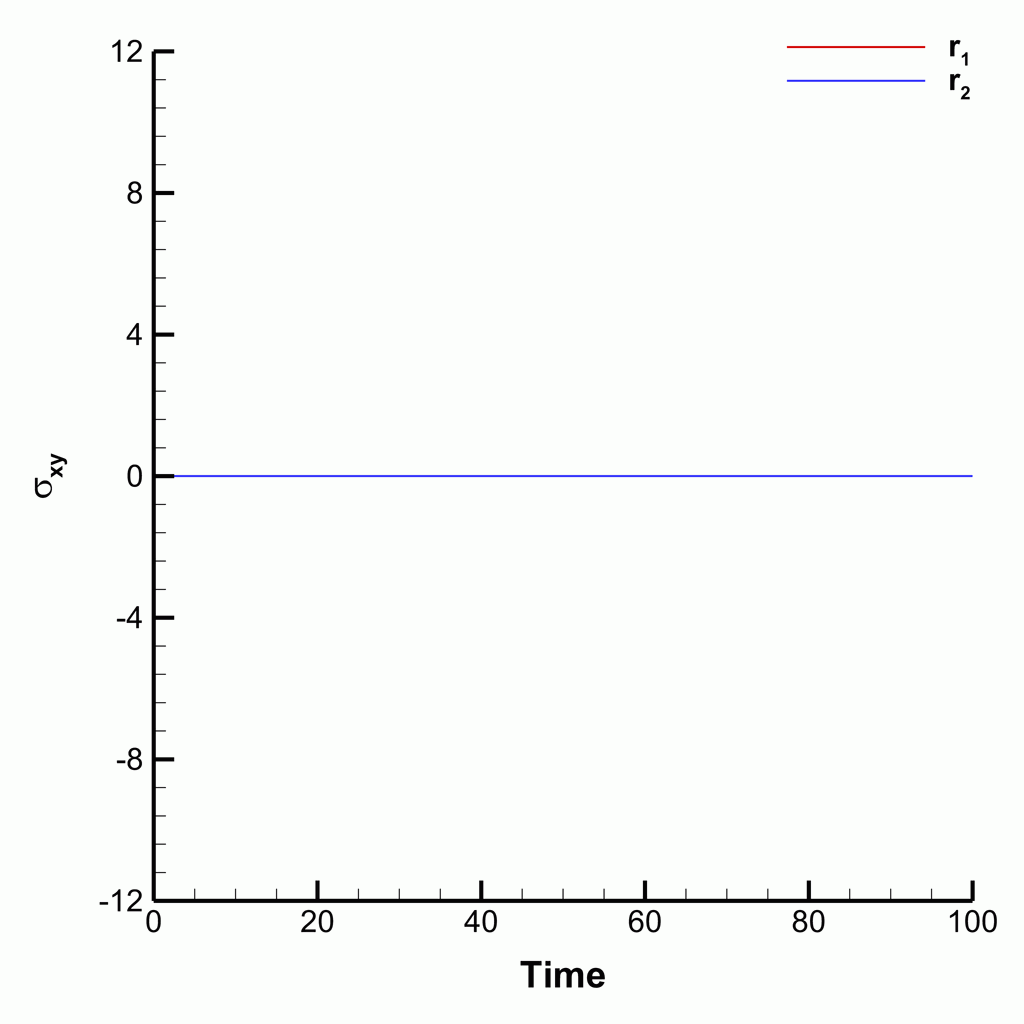}\hfill
	\includegraphics[width=0.32\linewidth]{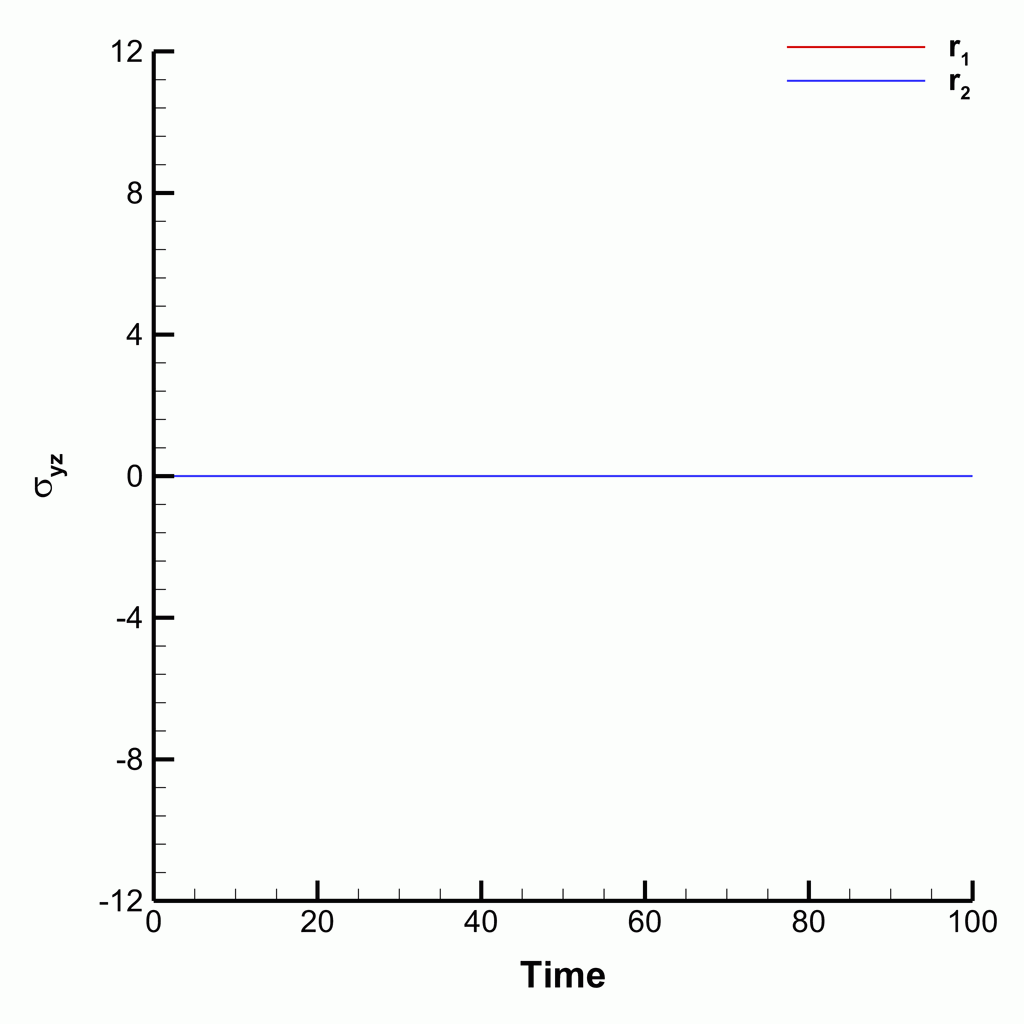}\hfill
	\includegraphics[width=0.32\linewidth]{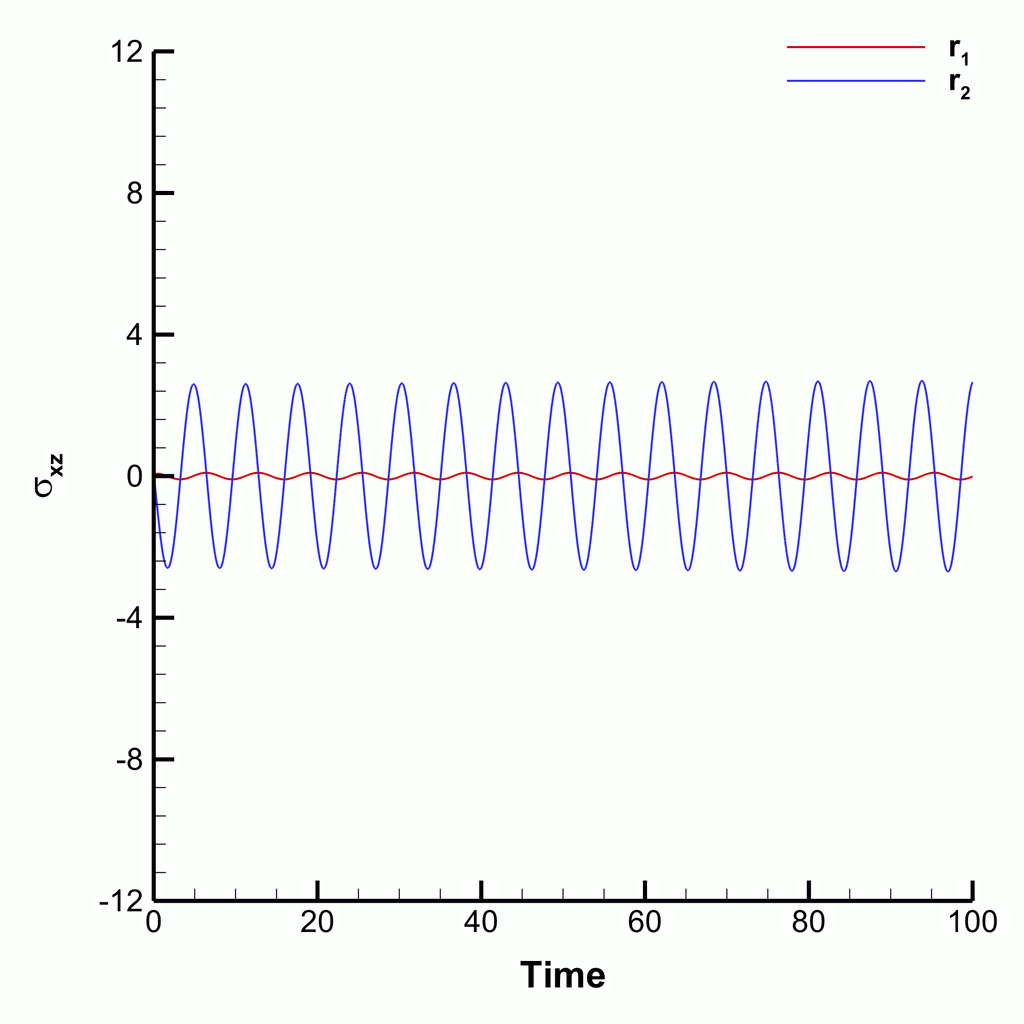}
	
	\includegraphics[width=0.32\linewidth]{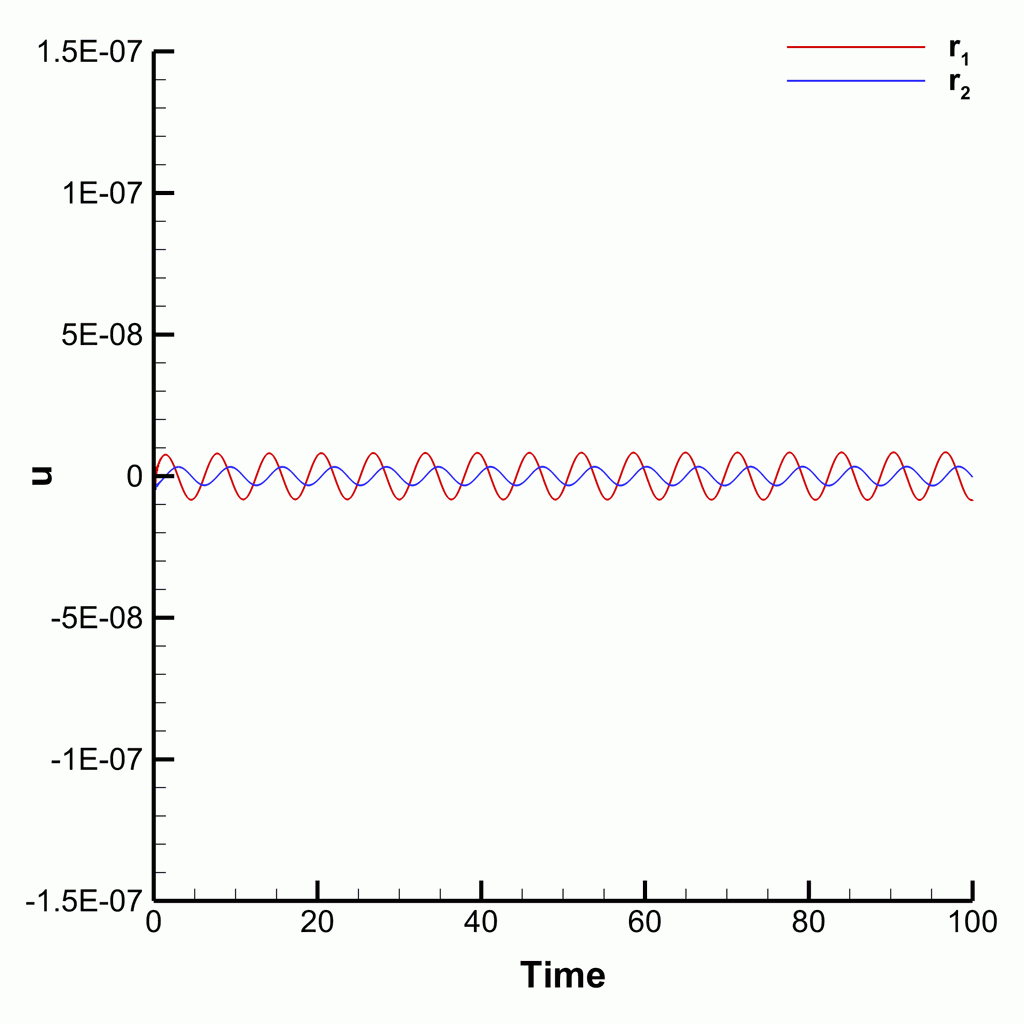}\hfill
	\includegraphics[width=0.32\linewidth]{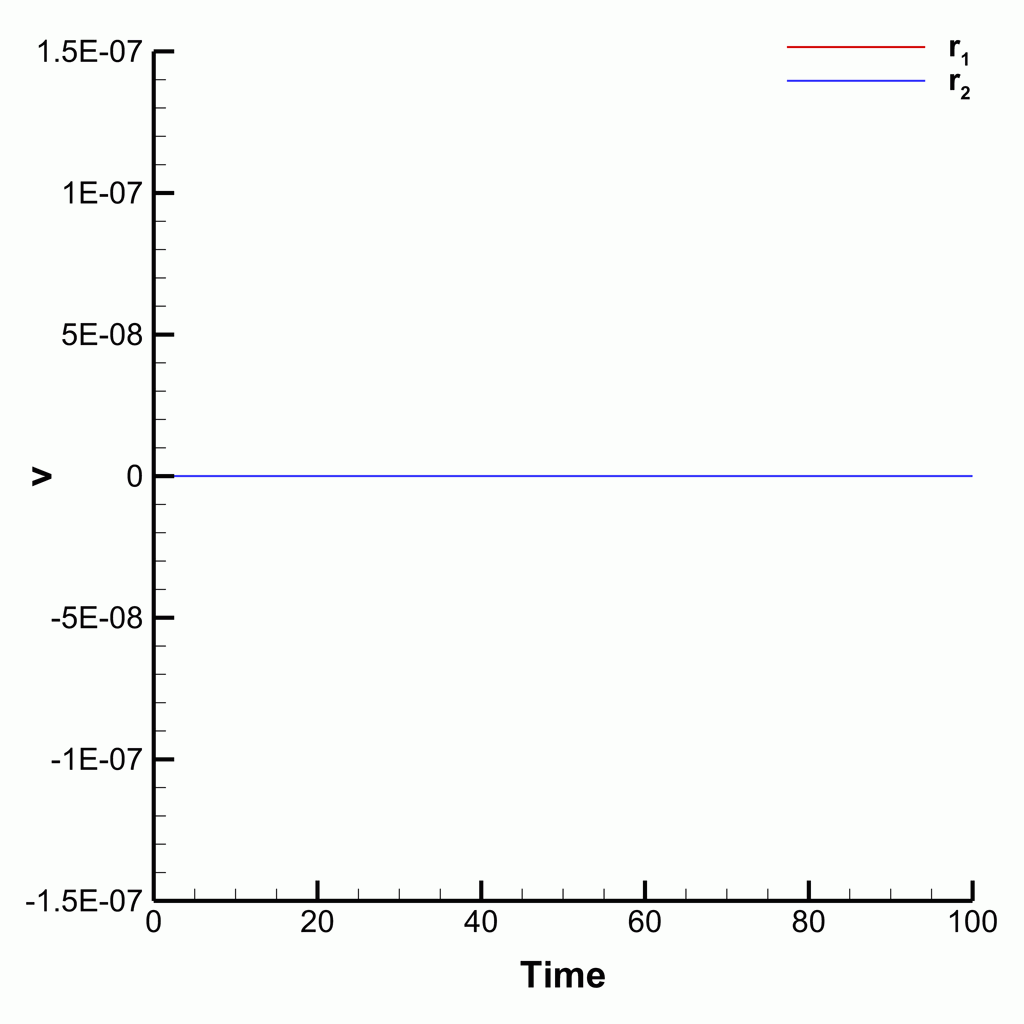}\hfill
	\includegraphics[width=0.32\linewidth]{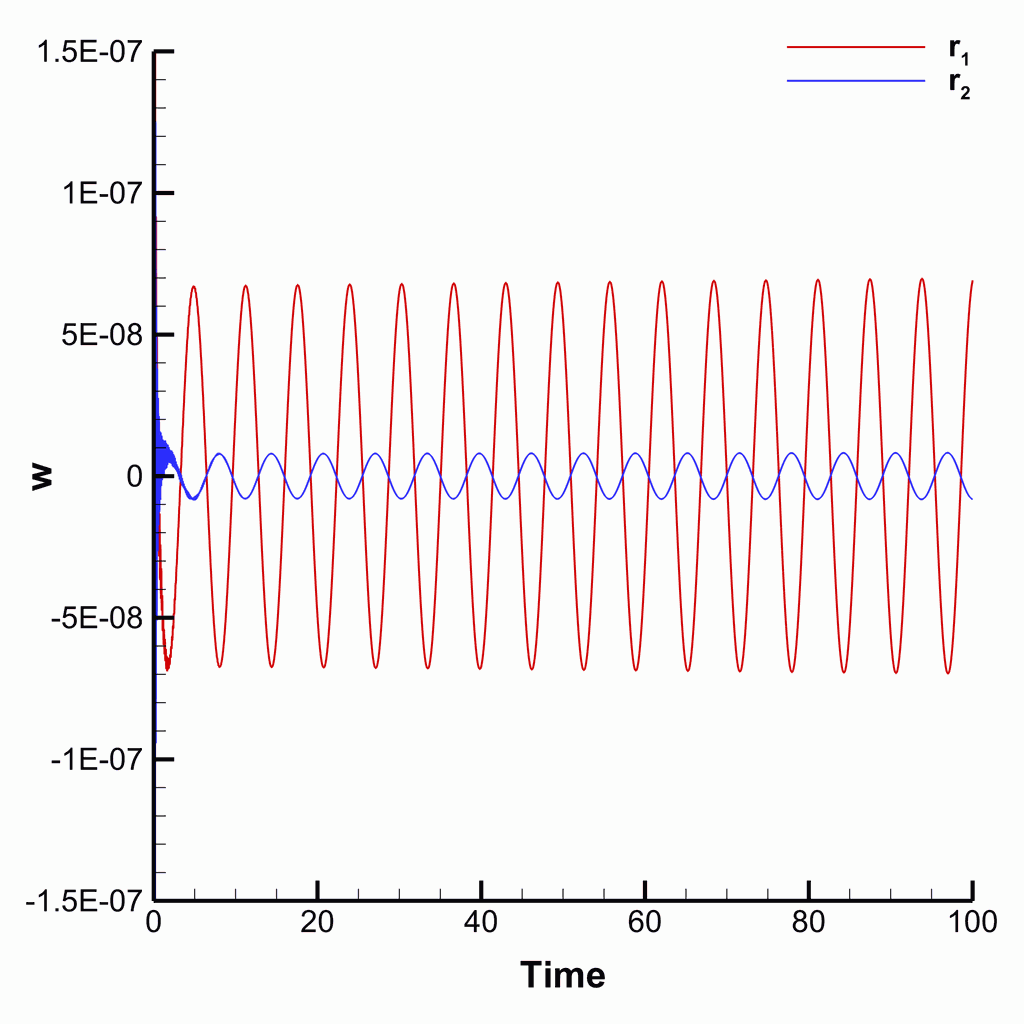}
	
	\caption{Sinusoidal wave. Seismograms at receivers $\mathbf{x}_{r_{1}}=\left(0,0,0\right)$ (red line) and $\mathbf{x}_{r_{2}}=\left(500, 500,-100 \right)$ (blue line) obtained using ADER-DG $\mathcal{O}4$. From left top to right bottom: $\sigma_{xx}$, $\sigma_{yy}$, $\sigma_{zz}$, $\sigma_{xz}$, $u$, $v$ and $w$.}
	
	\label{fig:SiW_pickpoints}
\end{figure}

\begin{figure}
	\centering	
	\includegraphics[width=0.4\linewidth]{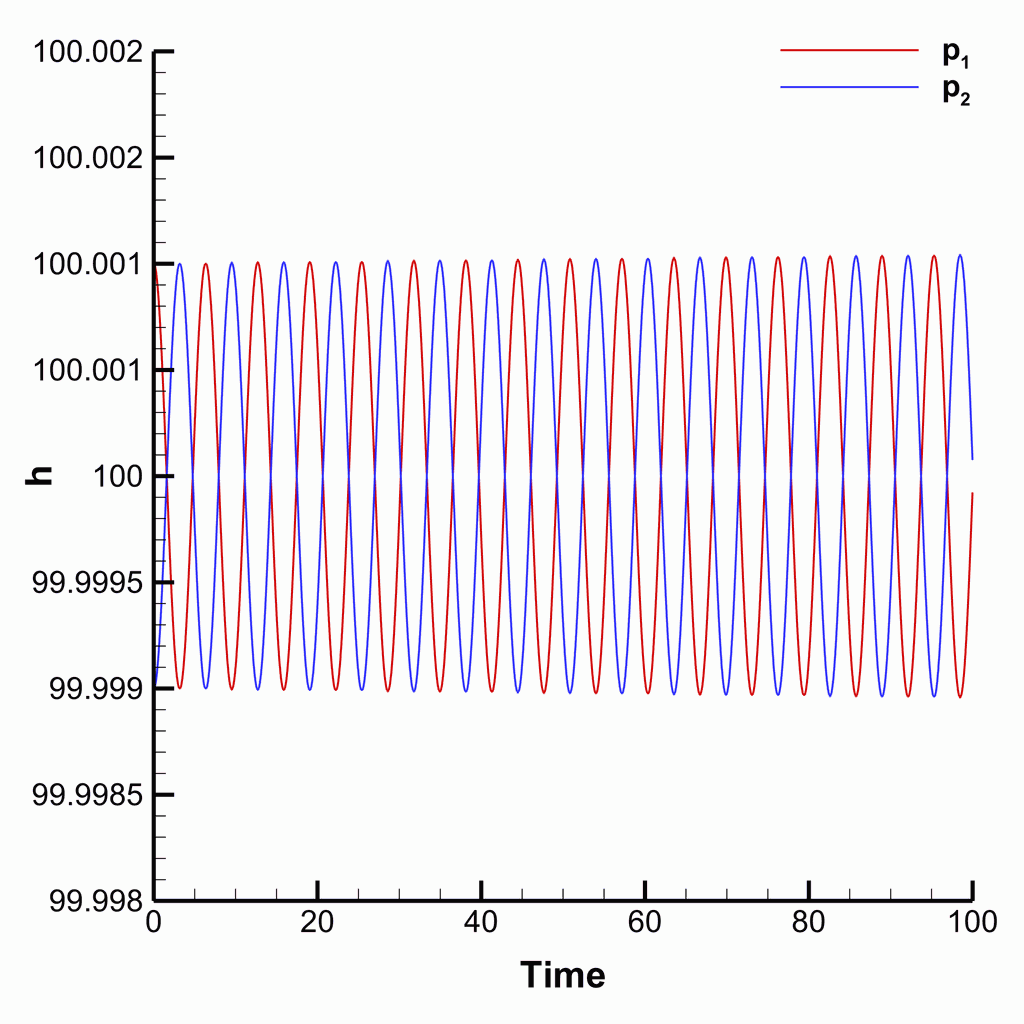}\hspace{0.05\linewidth}
	\includegraphics[width=0.4\linewidth]{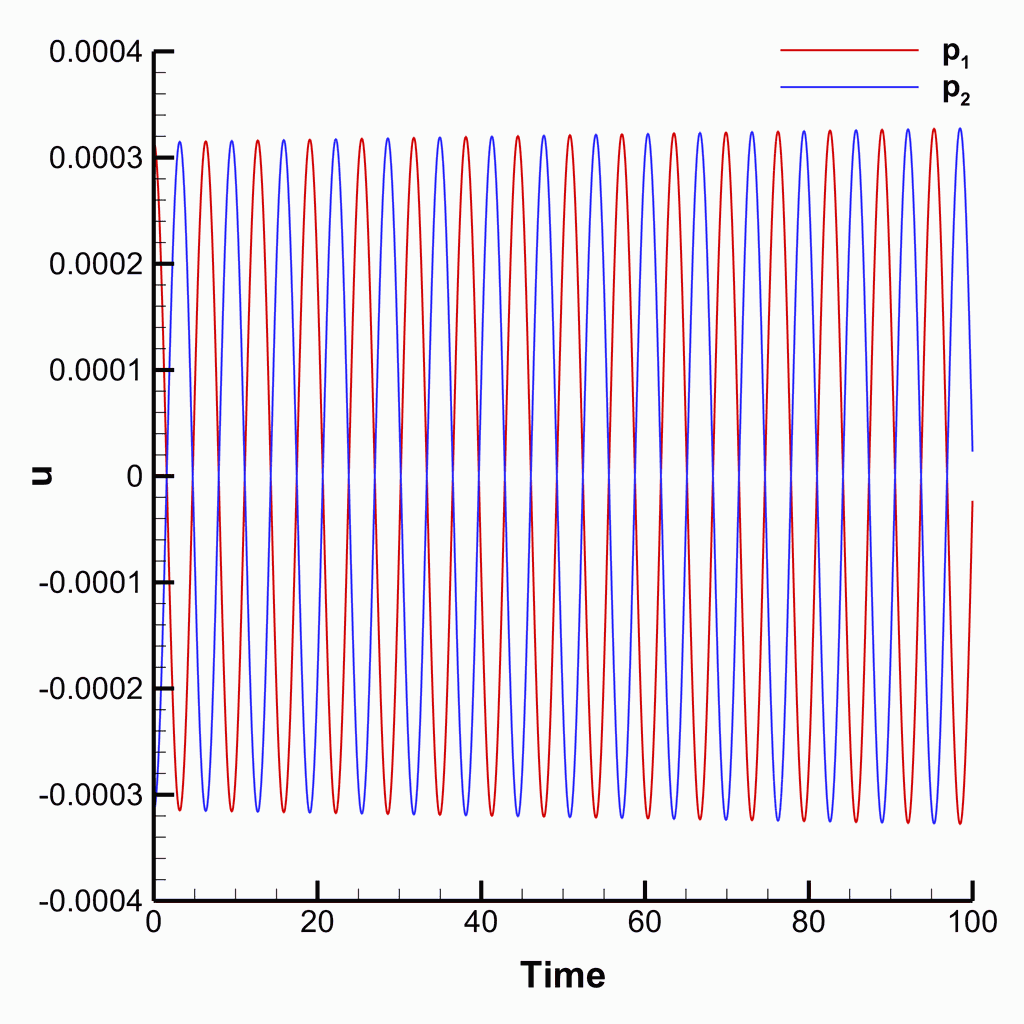}
	
	\vspace{0.05\linewidth}
	\includegraphics[width=0.4\linewidth]{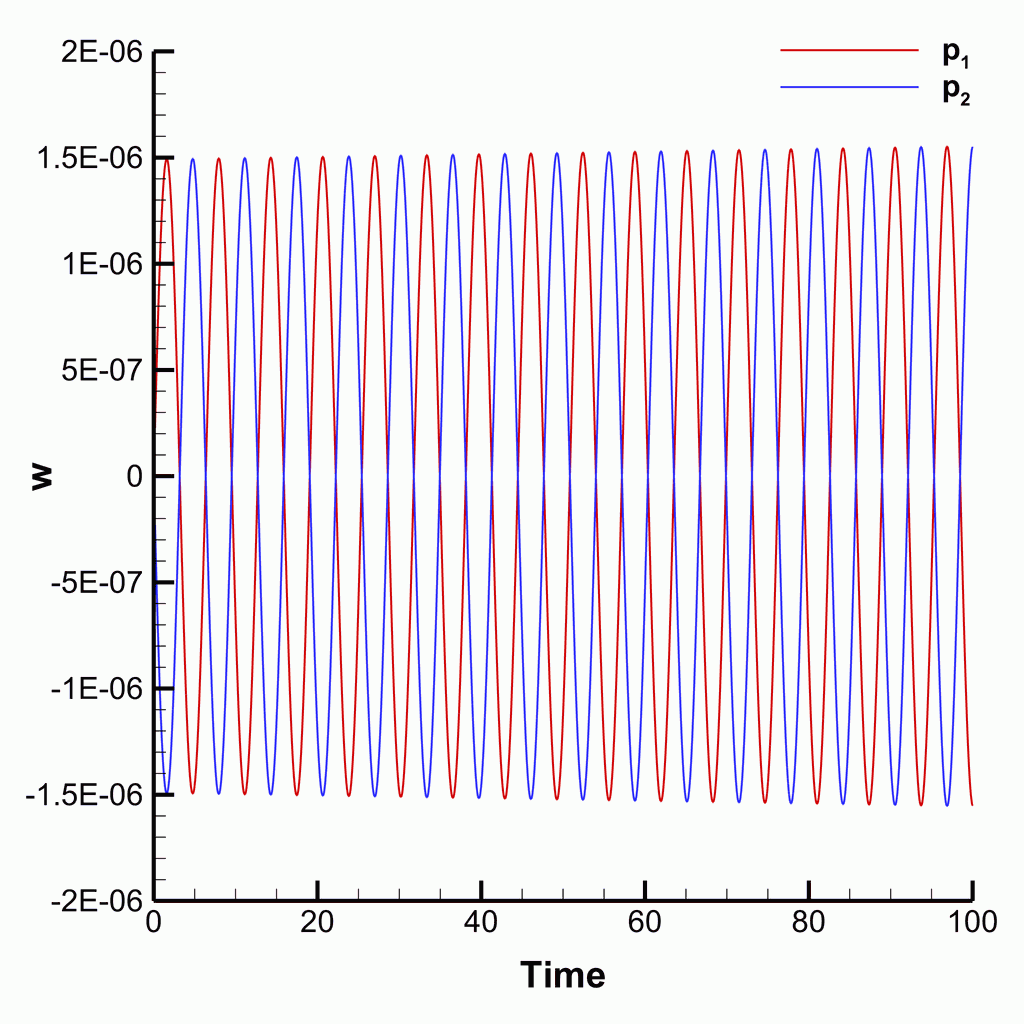}\hspace{0.05\linewidth}
	\includegraphics[width=0.4\linewidth]{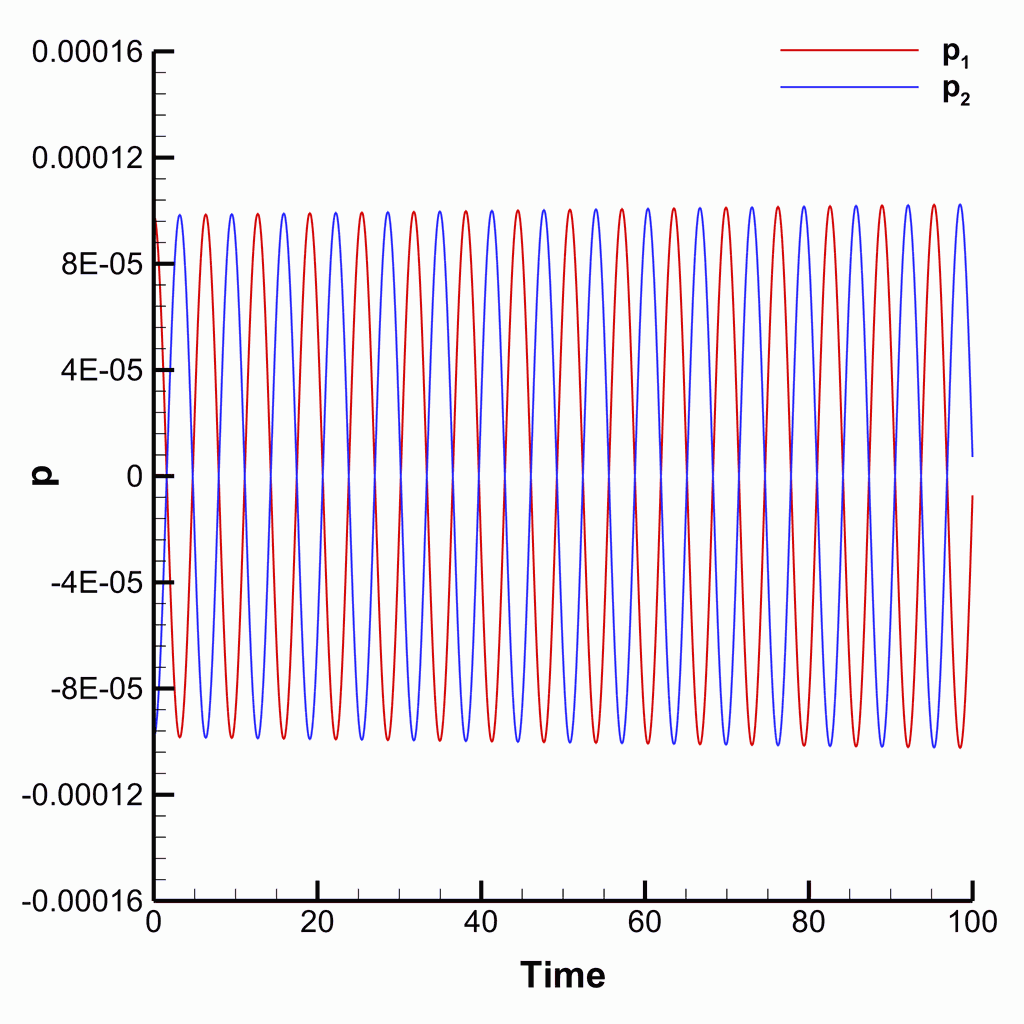}

	\caption{Sinusoidal wave. Time evolution of the main variables of the HSGN model at pick points  $\mathbf{x}_{p_{1}}=\left(0,0\right)$ (red line) and $\mathbf{x}_{p_{2}}=\left(500, 500 \right)$ (blue line) obtained using ADER-DG $\mathcal{O}4$. From left top to right bottom: $h$, $u$, $w$, and $p$.}
	
	\label{fig:SiW_pickpoints2D}
\end{figure}

Finally, we propose a sinusoidal wave propagation test.
We consider the computational domains $\Omega_{w}=[-5000,5000]\times [-2500,2500]$ and $\Omega_e=[-5000,5000]\times [-2500,2500]\times [-20000,0]$. 
The initial condition for the HSGN model is
computed taking into account the dispersion relation of the linearized system,
\begin{equation}
\frac{\partial \mathbf{Q}}{\partial t} + \mathbf{A}\left(\mathbf{Q}\right)\frac{\partial \mathbf{Q}}{\partial x} = \mathbf{S}\left(\mathbf{Q}\right), \quad \mathbf{A}\left(\mathbf{Q}\right) = \frac{\partial \mathbf{F}}{\partial \mathbf{Q}} + \mathbf{B}\left(\mathbf{Q}\right).
\end{equation}
That is, we first decompose the conservative variables into a stationary part plus the contribution of small time dependent fluctuations, i.e. $\mathbf{Q}(x,t)=\mathbf{Q}_{0}(x)+\mathbf{Q}^{\prime}(x,t)$. Next, we assume  $\mathbf{Q}^{\prime}=\hat{\mathbf{Q}}e^{i\left(kx-\omega t\right)}$, where $k$ denotes the wave number and $\omega$ is the angular frequency, obtaining the eigenvalue problem
\begin{equation}
\left[ \mathbf{A}\left(\mathbf{Q}_{0}\right)+i \mathbf{E}\left(\mathbf{Q}_{0}\right)\right] \mathbf{Q}^{\prime} = \omega \mathbf{I} \mathbf{Q}^{\prime}, \quad \mathbf{E}\left(\mathbf{Q}_{0}\right)= \frac{\partial \mathbf{S}}{\partial \mathbf{Q}}\left(\mathbf{Q}_{0}\right).\label{eq:linearizedeigenvalueproblem}
\end{equation}
Once \eqref{eq:linearizedeigenvalueproblem} is solved, we select one real non zero eigenvalue, 
\begin{equation}
\lambda = \dfrac{ \displaystyle{\sqrt{2 c_{0}^2 \gamma+ g H_{0}^3 k^2 
			+ c_{0}^2 H_{0}^2 k^2
			+ \displaystyle{\sqrt{g^2 H_{0}^6 k^4 
					+ 2 g c_{0}^2 H_{0}^5 k^4 
					- 4 g c_{0}^2 H_{0}^3 \gamma k^2
					+ c_{0}^4 H_{0}^4 k^4 
					+ 4 c_{0}^4 H_{0}^2 \gamma k^2 
					+ 4 c_{0}^4 \gamma^2} }}}}{\sqrt{2} H_{0}},
\end{equation}
and its corresponding eigenvector, and we compute the real part of the associated $\mathbf{Q}^{\ast}$,
\begin{equation}
\mathbf{r}_{l} = \begin{pmatrix}
1\\
\dfrac{\lambda}{k}\\
0\\
-\dfrac{H_{0}\lambda (c_{0}^{2} k^{2} + g H_{0} k^{2}- \lambda^{2}) \tan(\lambda t - k x) }{2 c_{0}^{2} k^{2}}\\[10pt]
-\dfrac{H_{0}^{2} \lambda^{2} (c_{0}^{2} k^{2} + g H_{0} k^{2} - \lambda^{2})}{2 \gamma c_{0}^{2} k^{2}}
\end{pmatrix}.\label{eq:eigsol_lp}
\end{equation}
Thus, adding the corresponding lake at rest solution,
\begin{equation}
\mathbf{Q}_{0} = \left(H_{0}, 0, 0, 0, 0 \right)
\end{equation}
to \eqref{eq:eigsol_lp} multiplied  by the sinusoidal function
\begin{equation}
f(x,y,0) =  10^{-3} \cos\left( \frac{ 2\pi x}{s}\right),
\end{equation}
yields the initial condition
\begin{equation}
\mathbf{Q}^{HSGN}(\mathbf{x},0) = \mathbf{Q}_{0}^{\mathrm{HSGN}} + f(x,y,0) \mathbf{r}_{l},
\end{equation}
with $H_{0}=100$ the still water depth, $s=200$ the wave length, $k=\frac{2\pi}{s}$, $g=9.81$ the gravity, and  $c=\sqrt{g H_{0}}$ the celerity. Zero homogeneous initial conditions are considered for the linear elasticity model. Likewise in the previous test case, we define periodic boundary conditions on $x$ and $y$ directions and a free surface boundary condition on the bottom of $\Omega_e$. The solid domain is meshed using $108\times 54 \times 100$ hexahedral elements and the upper surface is refined with refinement factor $5$ in each spatial direction to get the two-dimensional grid on $\Omega_w$.
The simulation is run until $t_{\mathrm{end}}=100$ using a fourth order ADER-DG scheme in both domains. The contour plots, on $y=0$, of the main unknowns of the linear elasticity model are reported in Figure \ref{fig:SiW_pickpoints} for $t=50$. Moreover, 
Figure \ref{fig:SiW_szz3D_subfig} shows the 3D isosurfaces of $\sigma_{zz}$. 
The seismogram recorded at receivers $\mathbf{x}_{r_{1}}=\left(0,0,0\right)$ and $\mathbf{x}_{r_{2}}=\left(500, 500,-100 \right)$ is depicted in Figure \ref{fig:SiW_pickpoints}. 
We observe the expected sinusoidal signal propagating in the solid domain. Like in the previous test case the apparent misbehaviour observed at the beginning of the simulation is caused by transients due to the initialization with zero of all the variables in the linear elastic wave propagation model.  
To check the correct displacement of the waves we also include the time evolution of the water waves at $\mathbf{x}_{p_{1}}=\left(0,0\right)$ and $\mathbf{x}_{p_{2}}=\left(500, 500 \right)$, see Figure \ref{fig:SiW_pickpoints2D}.

\begin{figure}
	\centering	
	
	\includegraphics[width=0.88\linewidth]{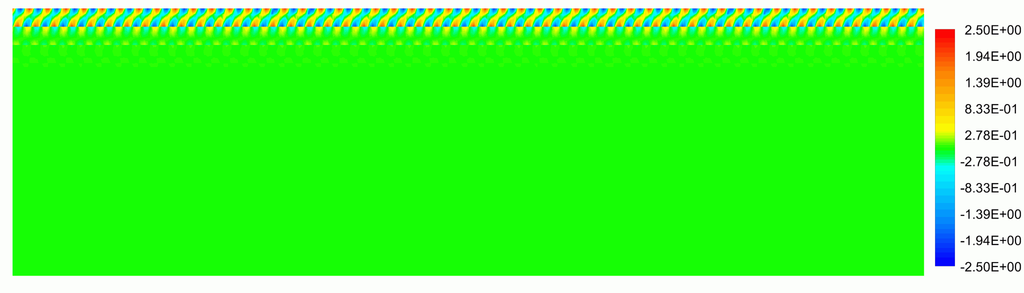}
	
	\includegraphics[width=0.88\linewidth]{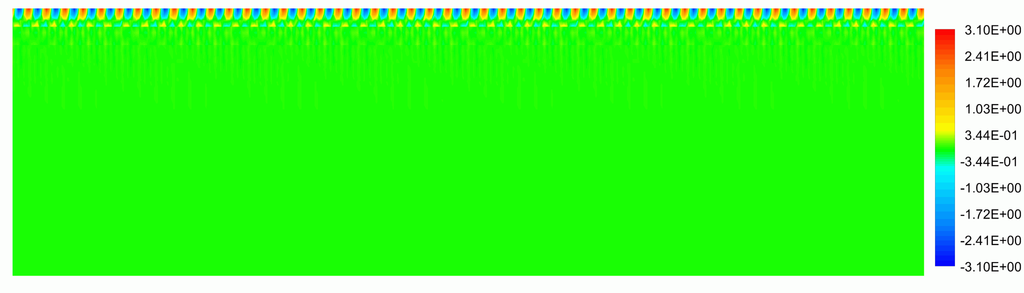}
	
	\includegraphics[width=0.88\linewidth]{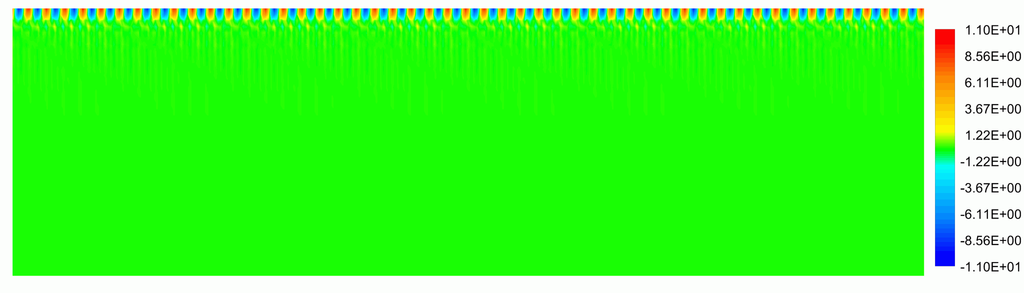}
	
	\includegraphics[width=0.88\linewidth]{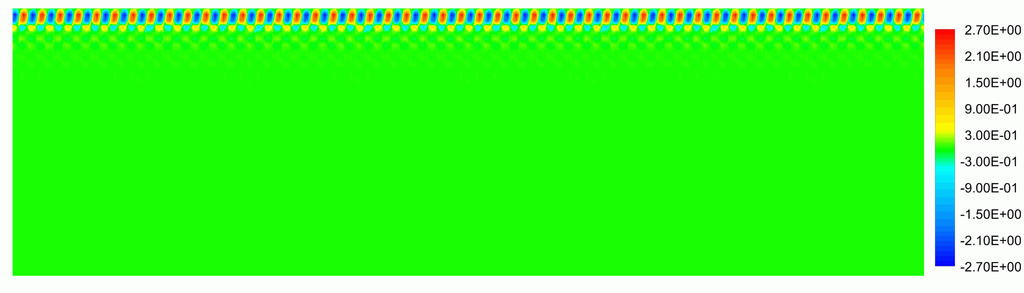}
	
	\includegraphics[width=0.88\linewidth]{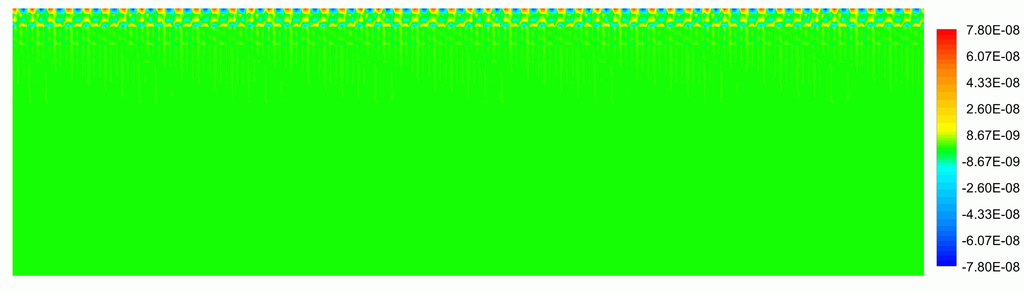}
	
	\caption{Sinusoidal wave. Solution obtained at time $t=50$. From left top to bottom right: contour plots at plane $y=0$, $z>-3000$, of $\sigma_{xx}$, $\sigma_{yy}$, $\sigma_{zz}$, $\sigma_{xz}$, and $w$, respectively.}
	
	\label{fig:SiW_3Dplot}
\end{figure}

\section{Conclusions}\label{sec:conclusions}
In this work we have presented a high order explicit ADER-DG method for the resolution of a one-way coupled system, describing seismic waves in the sea bottom generated by free surface water waves.  Classical formulations of non-hydrostatic dispersive systems and linear elasticity equations involve high order time and space derivatives, yielding to severe time step restrictions for explicit numerical schemes. To overcome this issue, we propose the use of first order hyperbolic reformulations of the original systems, which lead to classical CFL restrictions with $\Delta t \propto \Delta x$. In particular, we have considered the hyperbolic reformulation of the Serre-Green-Naghdi model for non-hydrostatic free surface flows and the first order velocity-stress formulation of linear elasticity for seismic wave propagation. Moreover, the use of hyperbolic models allows an easier coupling of the equations and a unified discretization based on one and the same method, i.e. employing the well known discontinuous Galerkin finite element method. 
High order of accuracy in space and time has been achieved using the ADER-DG methodology based on performing a local reconstruction of the data at each cell at the aid of a space-time predictor and the subsequent correction of the obtained approximation by considering the intercell flux within a classical space DG scheme.
The use of a common methodology to solve both PDE systems eases the coupling between the models, which has been done by imposing the normal stress on the upper boundary of the solid domain taking into account the water column heigh computed using the HSGN model. The large discrepancy between the wave lengths in the two media is addressed by considering non-conforming Cartesian grids in the two domains. Firstly, a three-dimensional mesh for the solid domain is designed. Then, the face mesh obtained on the upper boundary of the solid domain is refined to get a mesh for the fluid domain.
A careful assessment of the developed methodology has been performed. Several benchmarks for the linear elasticity equations and for non-hydrostatic dispersive free-surface flows are studied, showing excellent agreement of the obtained results with available reference data. Finally, two new tests of the coupled problem, considering solitary and sinusoidal free surface waves have been presented and allow to successfully validate the proposed approach. 

Within this work we have assumed to have a smooth bathymetry, so that the mild bottom approximation, used for the derivation of the HSGN model, holds. However, practical applications might also involve non mild bottom topographies. Therefore, future research would study the use of non-hydrostatic free-surface models for arbitrary bottom that may enlarge the applicability of the developed methodology. Besides, the simulation of non-hydrostatic flows is done using a dispersive shallow water type model, attending to the reduced depth of the water in comparison with the horizontal dimensions of the considered domain. An alternative approach that might be studied, when smaller differences between the spatial dimensions are involved, is the coupling of the linear elasticity model with the fully three dimensional free surface Navier-Stokes equations. 

\clearpage

\section*{Acknowledgements}
This work was financially supported by INdAM (\textit{Istituto Nazionale di Alta Matematica}, Italy) 
under two Post-doctoral grants of the research project \textit{Progetto premiale FOE 2014-SIES}; M.D. acknowledges partial support by the European Union's Horizon 2020 Research and Innovation Programme under 
the project \textit{ExaHyPE}, grant no. 671698 (call FETHPC-1-2014).  
The authors acknowledge funding from the Italian Ministry of Education, University 
and Research (MIUR) in the frame of the Departments of Excellence Initiative 2018--2022 
attributed to DICAM of the University of Trento (grant L. 232/2016) and in the frame of the 
PRIN 2017 project \textit{Innovative numerical methods for evolutionary partial differential equations and  applications}. Furthermore, M.D. has also received funding from the University of Trento via the Strategic Initiative \textit{Modeling and Simulation}.

\bibliographystyle{plain}
\bibliography{./References}

\end{document}